\documentclass[a4paper,12pt]{amsart}
\usepackage[frenchb,english]{babel}
\usepackage[latin1]{inputenc}
\usepackage{pdfsync}
\usepackage[T1]{fontenc}   
\usepackage{amsmath}
\usepackage{amssymb}
\usepackage{amsxtra}
\usepackage{amscd}
\usepackage{amsthm}
\usepackage{amsfonts}
\usepackage{eucal}
\usepackage[all]{xy}
\usepackage{graphicx}
\usepackage{comment}
\usepackage{epsfig}
\usepackage{psfrag}
\usepackage{mathrsfs}
\usepackage{amscd}
\usepackage{rotating}
\usepackage{lscape}
\usepackage{amsbsy}
\usepackage{verbatim}
\usepackage{moreverb}
\usepackage{url}
\usepackage{color}
\usepackage{xcolor}

\newcommand{\nc}{\newcommand}
\nc{\renc}{\renewcommand}

\nc\restr[2]{{ 
  \left.\kern-\nulldelimiterspace    #1  
  \vphantom{\big|}  
  \right|_{#2}  
  }}

\newtheorem{cor}[subsubsection]{Corollary}
\newtheorem{lem}[subsubsection]{Lemma}
\newtheorem{prop}[subsubsection]{Proposition}

\newtheorem{construction}[subsubsection]{Construction}
\newtheorem{nota}[subsubsection]{Notation}

\newtheorem{thm}[subsubsection]{Theorem}

\renc{\sec}{\section}
\nc{\ssec}{\subsection}
\nc{\sssec}{\subsubsection}

\theoremstyle{definition}
\newtheorem{defi}[subsubsection]{Definition}
\newtheorem{example}[subsubsection]{Example}

\theoremstyle{remark}
\newtheorem{rem}[subsubsection]{Remark}

\numberwithin{equation}{section}

\nc{\on}{\operatorname}

\nc\wt{\widetilde}
\nc\wh{\widehat}
\nc\ol{\ov}
\nc{\oc}[1]{{\overset{\circ}{#1}}}
\nc{\ov}[1]{{\overline{#1}}}
\nc{\isor}[1]{{\xrightarrow[\raisebox{0.25 em}{\smash{\ensuremath{\sim}}}]{#1}}}
\nc{\modmod}{/ \! \! /}

\nc{\mc}{\mathcal}
\nc{\mf}{\mathfrak}
\nc{\mr}{\mathrm}
\nc{\mb}{\mathbb}
\nc{\mbf}{\mathbf}
\nc{\ms}{\mathscr}

\nc{\R}{{\mathbb R}}
\nc{\Z}{{\mathbb Z}}
\nc{\N}{{\mathbb N}}
\nc{\C}{{\mathbb C}}
\nc{\Q}{{\mathbb Q}}

\nc{\Fq}{{\mathbb F}_q}
\nc{\Fl}{{\mathbb F}_\ell}
\nc{\Fqbar}{\ol{{\mathbb F}_q}}
\nc{\Flbar}{\ol{{\mathbb F}_\ell}}
\nc{\Zl}{{\mathbb Z}_\ell}
\nc{\Zlbar}{\ol{{\mathbb Z}_\ell}}
\nc{\Ql}{{\mathbb Q}_\ell}
\nc{\Qlbar}{\ol{{\mathbb Q}_\ell}}
\nc{\hl}{\overset{\leftarrow}h{}}
\nc{\hr}{\overset{\rightarrow}h{}}
\nc{\Gr}{{\on{Gr}}}
\nc{\Hecke}{\on{Hecke}}
 \nc{\Hom}{\on{Hom}}
 \nc{\Coker}{\on{Coker}}
 \nc{\Ker}{\on{Ker}}
 \nc{\Lie}{\on{Lie}}
\nc{\Loc}{\on{Loc}}
\nc{\Pic}{\on{Pic}}
\nc{\Bun}{\on{Bun}}
\nc{\IC}{\on{IC}}
\nc{\Aut}{\on{Aut}}
\nc{\Perv}{\on{Perv}}
\nc{\pos}{{\on{pos}}}
\nc{\Sym}{\on{Sym}}

\nc{\ta} {{}^\tau}
\nc {\tu}[1]{{}^{\tau^{#1}}\!}

\nc{\Id}{\on{Id}}
\nc{\Fil}{\on{Fil}}
\nc{\pr}{\on{pr}}
\nc{\Res}{\on{Res}}
\nc{\cusp}{\on{cusp}}
\nc{\Frob}{\on{Frob}}
\nc{\diag}{\Delta}
\nc{\gr}{\on{gr}}
\nc{\Inj}{\on{Inj}}
\nc{\Bl}{\on{Bl}}
\nc{\dem}{\noindent {\bf Proof. }}
\nc{\cqfd}{{\ }\hfill $\square$ \vskip 1mm}
\nc{\s}[1]{\langle #1 \rangle}
\nc{\Cht}{\on{Cht}}
\nc{\isom}{\overset {\thicksim}{\to}}
\nc{\sm}{\smallsetminus}

\begin{document}

\title{Cuspidal cohomology of stacks of shtukas}

\author{Cong Xue}
\address{Cong Xue: DPMMS, Centre for Mathematical Sciences, Wilberforce Road, Cambridge, CB3 0WB, UK}
\email{cx233@cam.ac.uk}

\maketitle

\begin{abstract}
Let $G$ be a connected split reductive group over a finite field $\Fq$ and $X$ a smooth projective geometrically connected curve over $\Fq$. The $\ell$-adic cohomology of stacks of $G$-shtukas is a generalization of the space of automorphic forms with compact support over the function field of $X$. In this paper, we construct a constant term morphism on the cohomology of stacks of shtukas which is a generalization of the constant term morphism for automorphic forms. We also define the cuspidal cohomology which generalizes the space of cuspidal automorphic forms. Then we show that the cuspidal cohomology has finite dimension and that it is equal to the (rationally) Hecke-finite cohomology defined by V. Lafforgue.
\end{abstract}


\section*{Introduction}

Let $X$ be a smooth projective geometrically connected curve over a finite field $\Fq$. We denote by $F$ its function field, by $\mb A$ the ring of adèles of $F$ and by $\mb O$ the ring of integral adèles.

Let $G$ be a connected split reductive group over $\Fq$. For simplicity, we assume in the introduction that the center of $G$ is finite.

\quad

We consider the space of automorphic forms $C_c(G(F) \backslash G(\mathbb{A}) / G(\mb O), \C)$. 
On the one hand, there is the notion of cuspidal automorphic form. 
An automorphic form is said to be cuspidal if its image under the constant term morphism along any proper parabolic subgroup of $G$ is zero. A theorem of Harder (\cite{harder} Theorem 1.2.1) says that the space of cuspidal automorphic forms has finite dimension. The proof uses the Harder-Narasimhan truncations and the contractibility of deep enough strata.

On the other hand, the space of automophic forms is equipped with an action of the Hecke algebra $C_c(G(\mb O) \backslash G(\mathbb{A}) / G(\mb O), \Q)$ by convolution on the right. An automorphic form is said to be (rationally) Hecke-finite if it belongs to a finite-dimensional subspace that is stable under the action of the Hecke algebra.

In \cite{vincent} Proposition 8.23, V. Lafforgue proved that the space of cuspidal automorphic forms and the space of Hecke-finite automorphic forms are equal. In fact, the space of cuspidal automorphic forms is stable under the action of the Hecke algebra and is finite-dimensional, thus it is included in the space of Hecke-finite automorphic forms. The converse direction follows from the following fact: any non zero image of the constant term morphism along a proper parabolic subgroup $P$ with Levi quotient $M$ is supported on the components indexed by a cone in the lattice of the cocharacters of the center of $M$. Hence it generates an infinite-dimensional vector space under the action of the Hecke algebra of $M$. Thus a non-cuspidal automorphic form can not be Hecke-finite for the Hecke algebra of $M$.

\quad

Let $\ell$ be a prime number not dividing $q$. In \cite{Dr78} and \cite{Dr87}, Drinfeld introduced the stacks classifying $GL_n$-shtukas for the representation $St \boxtimes St^*$ of $GL_n \times GL_n$, where $St$ is the standard representation of $GL_n$ and $St^*$ is its dual, and considered their $\ell$-adic cohomology. These were also used by L. Lafforgue in \cite{laurent}. 
Later in \cite{var}, Varshavsky defined the stacks classifying $G$-shtukas $\Cht_{G, I, W}$ for general $G$ and for an arbitrary representation $W$ of $\wh G^I$, where $\wh G$ is the Langlands dual group of $G$ over $\Ql$ and $I$ is a finite set (Drinfeld considered the case $G=GL_n$, $I=\{1, 2\}$ and $W=St \boxtimes St^*$). Varshavsky also defined the degree $j$ cohomology group with compact support $H_{G, I, W}^j$ of the $\ell$-adic intersection complex of $\Cht_{G, I, W}$ (this stack is smooth in the case of Drinfeld but not in general). In particular, when $I=\emptyset$ and $W = \bf 1$ is the one-dimensional trivial representation of the trivial group $\wh G^{\emptyset}$, the cohomology group $H_{G, \emptyset, \bf 1}^0$ coincides with $C_c(G(F) \backslash G(\mathbb{A}) / G(\mb O), \Ql)$.

\quad

The cohomology group $H_{G, I, W}^j$ is equipped with an action of the Hecke algebra $C_c( G(\mb O)  \backslash G(\mb A) / G(\mb O) , \Zl)$. In \cite{vincent}, V. Lafforgue defined the subspace $H_{G, I, W}^{j, \; \on{Hf}}$ of $H_{G, I, W}^j$ which consists of the cohomology classes $c$ for which $C_c( G(\mb O)  \backslash G(\mb A) / G(\mb O)  , \Zl) \cdot c$ is a finitely generated $\Zl$-submodule of $H_{G, I, W}^j$. When $I=\emptyset$ and $W=\bf 1$, the space $H_{G, \emptyset, \bf 1}^{0, \, \on{Hf}}$ coincides with the space of Hecke-finite automorphic forms, thus coincides with the space of cuspidal automorphic forms. V. Lafforgue used $H_{G, I, W}^{0, \; \on{Hf}}$ to construct the excursion operators on the space of cuspidal automorphic forms and obtained a canonical decomposition of this space indexed by the Langlands parameters.

We can also define a subspace $H_{G, I, W}^{j, \; \on{Hf-rat}}$ of $H_{G, I, W}^j$ which consists of the cohomology classes $c$ for which $C_c( G(\mb O)  \backslash G(\mb A) / G(\mb O) , \Ql) \cdot c$ is a finite-dimensional $\Ql$-vector subspace of $H_{G, I, W}^j$.
By definition, we have $H_{G, I, W}^{j, \; \on{Hf}} \subset H_{G, I, W}^{j, \; \on{Hf-rat}}$. When $I=\emptyset$ and $W=\bf 1$, it is easy to see that they are equal. 

\quad

In this paper, we are interested in the constant term morphism of the cohomology of stacks of shtukas, analogous to the case of automorphic forms. For any parabolic subgroup $P$ of $G$, let $M$ be its Levi quotient. As in \cite{var}, we can define the stack of $P$-shtukas $\Cht_{P, I, W}$ and the stack of $M$-shtukas $\Cht_{M, I, W}$. The morphisms $G \hookleftarrow P \twoheadrightarrow M$ induce a correspondence $$\Cht_{G, I, W}  \leftarrow \Cht_{P, I, W}  \rightarrow \Cht_{M, I, W} . $$
From this we construct a constant term morphism $$C_{G}^{P, \, j}: H_{G, I, W}^j \rightarrow H_{M, I, W}^j.$$ Then we define the cuspidal cohomology $H_{G, I, W}^{j, \; \on{cusp}} \subset H_{G, I, W}^j$ as the intersection of the kernels of the constant term morphisms for all proper parabolic subgroups.

This construction was suggested by V. Lafforgue. He also conjectured that 

- the cuspidal cohomology is of finite dimension;

- the following three $\Ql$-vector subspaces of  $H_{G, I, W}^{j}$ are equal:
$$H_{G, I, W}^{j, \; \on{Hf}} = H_{G, I, W}^{j, \; \on{Hf-rat}} = H_{G, I, W}^{j, \; \on{cusp}} .$$

In this paper, we prove these conjectures except for the equality with $H_{G, I, W}^{j, \; \on{Hf}}$, which we plan to treat in a future paper. The main results are:

\begin{thm}  \label{thm-cusp-dim-fini}  (Theorem \ref{thm-cusp-dim-fini-second})
The $\Ql$-vector space $H_{G, I, W}^{j, \; \on{cusp}}$ has finite dimension.
\end{thm}

\begin{prop}  \label{prop-cusp-egal-Hfrat}  (Proposition \ref{prop-cusp-egal-Hfrat-section-6})
The two $\Ql$-vector subspaces $H_{G, I, W}^{j, \; \on{cusp}}$ and $H_{G, I, W}^{j, \; \on{Hf-rat}}$ of $H_{G, I, W}^{j}$ are equal.
\end{prop}

As a consequence, $H_{G, I, W}^{j, \; \on{Hf}}$ has finite dimension.

In particular, when $I=\emptyset$ and $W = \bf 1$, the constant term morphism $C_{G}^{P, \, 0}$ coincides with the usual constant term morphism for automorphic forms. In this case, Theorem \ref{thm-cusp-dim-fini} coincides with Theorem 1.2.1 in \cite{harder}, and Proposition \ref{prop-cusp-egal-Hfrat} coincides with Proposition 8.23 in \cite{vincent} mentioned before.

Let $N \subset X$ be a finite subscheme. Theorem \ref{thm-cusp-dim-fini} and Proposition \ref{prop-cusp-egal-Hfrat} are still true for the cohomology with level structure on $N$.

\subsection*{Structure of the paper}

In Section 1 we construct the parabolic induction diagram and define Harder-Narasimhan truncations which are compatible with the parabolic induction. In Section 2 we recall the cohomology of the stacks of $G$-shtukas and define the cohomology of the stacks of $M$-shtukas. In Section 3 we construct the constant term morphism using the compatibility of the geometric Satake equivalence with the constant term functors for the Beilinson-Drinfeld affine grassmannians. 

The idea of the proofs of Theorem \ref{thm-cusp-dim-fini} and Proposition \ref{prop-cusp-egal-Hfrat} is analogous to the case of automorphic forms.
The goal of Sections 4-5 is to prove Theorem \ref{thm-cusp-dim-fini}. In Section 4 we prove the contractibility of deep enough horospheres. In Section 5 we use this result and an argument by induction on the semisimple rank to prove the finiteness of cuspidal cohomology. In Section 6 we show that the constant term morphism commutes with the action of the Hecke algebra, and we prove Proposition \ref{prop-cusp-egal-Hfrat}.

\subsection*{Notations and conventions}

\sssec{} \label{subsection-def-Z-G}
Let $G$ be a connected split reductive group over $\Fq$. Let $G^{\mr{der}}$ be the derived group of $G$ and $G^{\mr{ab}}:=G/G^{\mr{der}}$ the abelianization of $G$.
Let $Z_G$ be the center of $G$ and $G^{\mr{ad}}$ the adjoint group of $G$ (equal to $G / Z_G$).

\sssec{}  \label{subsection-def-Xi}
We fix a discrete subgroup $\Xi_G$ of $Z_G(\mb A)$ such that $\Xi_G \cap Z_G(\mb O) Z_G(F) = \{ 1 \}$, the quotient $Z_G(F) \backslash Z_G(\mb A) / Z_G(\mb O) \Xi_G $ is finite and the composition of morphisms $\Xi_G \hookrightarrow Z_G(\mb A) \hookrightarrow G(\mb A) \twoheadrightarrow G^{\mr{ab}}(\mb A)$ is injective. Note that the volume of $G(F) \backslash G(\mb A) / G(\mb O) \Xi_G $ is finite. We write $\Xi:=\Xi_G$.

\sssec{} We fix a Borel subgroup $B \subset G$. 
By a parabolic subgroup we will mean a standard parabolic subgroup (i.e. a parabolic subgroup containing $B$), unless explicitly stated otherwise.

\sssec{} \label{subsection-def-Lambda-H}
Let $H$ be a connected split reductive group over $\Fq$ with a fixed Borel subgroup.
Let $\Lambda_H$ (resp. $\wh{\Lambda}_H$) denote the weight (resp. coweight) lattice of $H$. Let $\langle \ , \ \rangle: \wh \Lambda_H \times \Lambda_H \rightarrow \Z$ denote the natural pairing between the two.

Let $\wh \Lambda_H^+ \subset \wh \Lambda_H$ denote the monoid of dominant coweights and $\wh \Lambda_H^{\mr{pos}} \subset \wh \Lambda_H$ the monoid generated by positive simple coroots. Let $\wh \Lambda_H^{\Q}:=\wh \Lambda_H \underset{\Z} \otimes \Q$. Let $\wh \Lambda_H^{\mr{pos}, \Q}$ and $\wh \Lambda_H^{+, \Q}$ denote the rational cones of $\wh \Lambda_H^{\mr{pos}}$ and $\wh \Lambda_H^{+}$. 
We use analogous notation for the weight lattice. 

We use the partial order on $\wh{\Lambda}_{H}^{\Q}$ defined by $\mu_1 \leq^H \mu_2 \Leftrightarrow \mu_2 - \mu_1 \in \wh \Lambda_H^{\mr{pos}, \Q}$
(i.e. $\mu_2 - \mu_1$ is a linear combination of simple coroots of $H$ with coefficients in $\Q_{\geq 0}$).

We will apply these notations to $H = G$, $H=G^{\mr{ad}}$ or $H=$ some Levi quotient $M$ of $G$.





\sssec{} \label{subsection-def-Gamma-G-Gamma-M}
We denote by $\Gamma_G$ the set of simple roots of $G$ and by $\wh\Gamma_G$ the set of simple coroots. The standard parabolic subgroups of $G$ are in bijection with the subsets of $\Gamma_G$ in the following way. To a parabolic subgroup $P$ with Levi quotient $M$, we associate the subset $\Gamma_M$ in $\Gamma_G$ equal to the set of simple roots of $M$. 

%


\sssec{} \label{subsection-O-N}
Let $N \subset X$ be a finite subscheme. We denote by $\mc O_N$ the ring of functions on $N$ and write $K_{G, N}:=\Ker(G(\mb O) \rightarrow G(\mc O_N) )$. 

Let $H$ be an algebraic group over $\Fq$. We denote by $H_{N}$ the Weil restriction $\on{Res}_{\mc O_{N} / \Fq} H$.


\sssec{} If not specified, all schemes are defined over $\Fq$ and all the fiber products are taken over $\Fq$.

\sssec{} For any scheme $S$ over $\Fq$ and $x$ an $S$-point of $X$, we denote by $\Gamma_x \subset X \times S$ the graph of $x$. 

\sssec{} \label{def-ta-inverse-image-of-Frob}
For any scheme $S$ over $\Fq$, we denote by $\Frob_S: S \rightarrow S$ the Frobenius morphism over $\Fq$. 
For any $G$-bundle $\mc G$ on $X \times S$, we denote by $\ta \mc G$ the $G$-bundle $(\Id_X \times_{\Fq} \Frob_S)^*\mc G$.

\sssec{} \label{def-prechat-chat-chat-alg}
We use Definition 3.1 and Definition 4.1 in \cite{chat-alg} for prestacks, stacks and algebraic stacks.

\sssec{}   \label{def-D-c-b-pour-chat-alg}
As in \cite{chat-alg} Section 18, \cite{LO08} and \cite{LO09}, for $\mathcal{X}$ an algebraic stack locally of finite type over $\Fq$, we denote by $D_c^b(\mathcal{X}, \Ql)$ the bounded derived category of constructible $\ell$-adic sheaves on $\mathcal{X}$. We have the notion of six operators and perverse sheaves. 

If $f: \mc X_1 \rightarrow \mc X_2$ is a morphism of finite type of schemes (resp. algebraic stacks) locally of finite type, we will denote by $f_!$, $f_*$, $f^*$, $f^!$ the corresponding functors between $D_c^b(\mathcal{X}_1, \Ql)$ and $D_c^b(\mathcal{X}_2, \Ql)$, always understood in the derived sense.

\sssec{}   \label{subsection-reduced-subscheme}
We will work with étale cohomology. So for any stack (resp. scheme) (for example $\Cht_{G, N, I, W}$ and $\Gr_{G, I, W}$), we consider only the reduced substack (resp. subscheme) associated to it.


\subsection*{Acknowledgments}
This paper is based on my PhD thesis \cite{these}, which I wrote at Université Paris-Sud. I thank my advisors Vincent Lafforgue and Gérard Laumon for suggestion of the project and support during the writing of my thesis. I thank Jack Thorne for his helpful comments on a draft of this paper. I thank the referees for their suggestions of improvements.

The author is supported by funding from the European Research Council (ERC) under the European Union's Horizon 2020 research and innovation programme (grant agreement No 714405).

\tableofcontents  

\section{Parabolic induction diagram of stacks of shtukas}


The goal of this section is to introduce the parabolic induction diagram of stacks of shtukas without a bound on the modifications at paws in Sections \ref{section-def-cht}-\ref{subsection-quotient-by-Xi} and to introduce the Harder-Narasimhan stratification for the parabolic induction diagram in Sections \ref{subsection-HN-truncation}-\ref{subsection-HN-Cht-G-P-M}.

In Sections \ref{section-def-cht}-\ref{subsection-quotient-by-Xi} we work in the context of prestacks (see \ref{def-prechat-chat-chat-alg}) 

\subsection{Reminder of stacks of shtukas and Beilinson-Drinfeld affine grassmannians}  \label{section-def-cht}

This subsection is based on Section 2 of \cite{var} and Sections 1 and 2 of \cite{vincent}. All the results are well-known.

\begin{defi}   \label{def-Bun-G-N}
We define $\Bun_{G, N}$ to be the prestack that associates to any affine scheme $S$ over $\Fq$ the groupoid 
$\Bun_{G, N}(S) : = \{ (\mc G, \psi), \text{ where } \mc G \text{ is a } G \text{-bundle on } X \times S, \; \psi \text{ is an isomorphism of } G \text{-bundles}: \restr{\mc G}{N \times S} \isom \restr{G}{N \times S}  \}.$
\end{defi}

\sssec{}
$\Bun_{G, N}$ is a smooth algebraic stack over $\Fq$, locally of finite type.

\begin{defi}  
We define $\on{Hecke}_{G, N, I}$ to be the prestack that associates to any affine scheme $S$ over $\Fq$ the groupoid $\on{Hecke}_{G, N, I}(S)$ that classifies the following data:

(i)  $(x_i)_{i \in I} \in (X\sm N)^I(S)$, 

(ii) $(\mc G, \psi), (\mc G', \psi') \in \Bun_{G, N}(S)$, 
 
(iii) an isomorphism of $G$-bundles $\phi:\restr{\mc G}{(X\times S)\sm(\bigcup_{i\in I}\Gamma_{x_i})} \isom \restr{\mc G'}{(X\times S)\sm(\bigcup_{i\in I}\Gamma_{x_i})}$ which preserves the $N$-level structure, i.e. $\psi'\circ \restr{\phi}{N\times S}=\psi$.
\end{defi}

\sssec{}
The prestack $\on{Hecke}_{G, N, I}$ is an inductive limit of algebraic stacks over $(X \sm N)^I$. We define the morphism of paws $\on{Hecke}_{G, N, I} \rightarrow (X \sm N)^I$ by sending $\big( (x_i)_{i\in I}, (\mc G, \psi) \xrightarrow{\phi}  (\mc G', \psi') \big)$ to $(x_i)_{i\in I}$.

\sssec{}
We denote by $\on{pr}_0$ (resp. $\on{pr}_1$) the projection $\on{Hecke}_{G, N, I} \rightarrow \Bun_{G, N}$ which sends $\big( (x_i)_{i\in I}, (\mc G, \psi) \xrightarrow{\phi}  (\mc G', \psi') \big)$ to $(\mc G, \psi)$ (resp. to $ (\mc G', \psi')$).

Let $\Frob: \Bun_{G, N} \rightarrow \Bun_{G, N}$ be the Frobenius morphism over $\Fq$. With the notation in \ref{def-ta-inverse-image-of-Frob}, for any affine scheme $S$ over $\Fq$, the morphism $\Frob: \Bun_{G, N}(S) \rightarrow \Bun_{G, N}(S)$ is given by $(\mc G, \psi) \rightarrow (\ta \mc G, \ta \psi).$

\begin{defi}   \label{def-Cht-G-N-I}
We define the prestack of shtukas $\on{Cht}_{G, N, I}$ to be the fiber product:
\begin{equation}
\xymatrixrowsep{2pc}
\xymatrixcolsep{4pc}
\xymatrix{
\on{Cht}_{G, N, I}  \ar[r]  \ar[d]
& \on{Hecke}_{G, N, I}   \ar[d]^{(\on{pr}_0, \; \on{pr}_1)} \\
\Bun_{G, N}  \ar[r]^{(\Id, \; \Frob) \quad \quad}
& \Bun_{G, N} \times_{\Fq} \Bun_{G, N}
}
\end{equation}
\end{defi}

\sssec{}
Concretely, $\on{Cht}_{G, N, I}$ is 
the prestack which associates to any affine scheme $S$ over $\Fq$ the groupoid $\on{Cht}_{G, N, I}(S)$ classifying the following data:

(i)  $(x_i)_{i \in I} \in (X\sm N)^I(S)$, 

(ii) $(\mc G, \psi) \in \Bun_{G, N}(S)$, 
 
(iii) an isomorphism of $G$-bundles $\phi:\restr{\mc G}{(X\times S)\sm(\bigcup_{i\in I}\Gamma_{x_i})} \isom \restr{\ta \mc G}{(X\times S)\sm(\bigcup_{i\in I}\Gamma_{x_i})}$ which preserves the $N$-level structure, i.e. $\ta \psi \circ \restr{\phi}{N\times S}=\psi$.

We define the morphism of paws $\mf p_G: \on{Cht}_{G, N, I} \rightarrow (X \sm N)^I$ by sending $\big( (x_i)_{i\in I}, (\mc G, \psi) \xrightarrow{\phi}  (\ta \mc G, \ta \psi) \big)$ to $(x_i)_{i\in I}$.


\sssec{}
The prestack $\on{Cht}_{G, N, I}$ is an inductive limit of algebraic stacks over $(X \sm N)^I$. 

\sssec{} We will omit the index $N$ if $N = \emptyset$.

\quad

We will need a local model of $\on{Cht}_{G, N, I}$. For this, we recall the definition of Beilinson-Drinfeld affine grassmannians.

\sssec{}    
For $(x_i)_{i \in I} \in X^I(S)$, $d \in \N$, we denote by $\Gamma_{\sum d x_i}$ the closed subscheme of $X \times S$ whose ideal is generated by $(\prod_{i \in I} t_i)^d$ locally for the Zariski topology, where $t_i$ is an equation of the graph $\Gamma_{x_i}$. We define $\Gamma_{\sum \infty x_i}:=\varinjlim _{d} \Gamma_{\sum d x_i}$ to be the formal neighborhood of $\cup_{i \in I}\Gamma_{x_i}$ in $X \times S.$

A $G$-bundle on $\Gamma_{\sum \infty x_i}$ is a projective limit of $G$-bundles on $\Gamma_{\sum d x_i}$ as $d \rightarrow \infty$.

\begin{defi} \label{def-Gr-G-I}
We define the Beilinson-Drinfeld affine grassmannian $\on{Gr}_{G, I}$ to be the ind-scheme that associates to any affine scheme $S$ over $\Fq$ the set $\on{Gr}_{G, I}(S)$ classifying the following data:

(i)  $(x_i)_{i \in I} \in X^I(S)$, 

(ii) $\mc G, \mc G'$ two $G$-bundles on $\Gamma_{\sum \infty x_i}$, 
 
(iii) an isomorphism of $G$-bundles $\phi: \restr{\mc G}{ \Gamma_{\sum \infty x_i} \sm (\bigcup_{i\in I}\Gamma_{x_i})} \isom \restr{\mc G'}{ \Gamma_{\sum \infty x_i } \sm (\bigcup_{i\in I}\Gamma_{x_i})}$ where the precise meaning is given in \cite{vincent} Notation 1.7,

(iv) a trivialization $\theta: \mc G' \isom G$ on $\Gamma_{\sum \infty x_i}$.
\end{defi}

\sssec{}   \label{subsection-usual-affine-grassmannian}
We have the morphism of paws: $\on{Gr}_{G, I} \rightarrow X^I$. The fiber over $(x_i)_{i \in I} \in X^I_{\Fqbar}$ is $\prod_{y \in \{ x_i | i \in I \}} \on{Gr}_{G, y}$, where $\on{Gr}_{G, y}$ is the usual affine grassmannian, i.e. the fpqc quotient $G_{\mc K_y} / G_{\mc O_y}$, where $\mc O_y$ is the complete local ring on $y$ and $\mc K_y$ is its field of fractions.

\begin{defi}   \label{defi-G-I-infry}
(a) For any $d \in \N$, we define $G_{I, d}$ to be the group scheme over $X^I$ that associates to any affine scheme $S$ over $\Fq$ the 
set consisting of pairs $((x_i)_{i \in I}, f)$, where
$(x_i)_{i \in I} \in X^I(S)$ and $f$ is
an automorphism of the trivial $G$-bundle on $\Gamma_{\sum d x_i}$.

(b) We define the group scheme $G_{I, \infty}:= \underset{\longleftarrow} \lim \, G_{I, d}$.
\end{defi}

\sssec{}
The fiber of $G_{I, \infty}$ over $(x_i)_{i \in I} \in X^I_{\Fqbar}$ is $\prod_{y \in \{ x_i | i \in I \}} G_{\mc O_y}$.

\sssec{}
The group scheme $G_{I, \infty}$ acts on $\on{Gr}_{G, I}$ by changing the trivialization $\theta$. 
We denote by $[G_{I, \infty} \backslash \on{Gr}_{G, I}]$ the quotient prestack. For any affine scheme $S$ over $\Fq$, $[G_{I, \infty} \backslash \on{Gr}_{G, I}](S)$ is the groupoid classifying the data (i), (ii) and (iii) in Definition \ref{def-Gr-G-I}.

\sssec{}    \label{subsection-epsilon-G-I-infty}
We have a morphism of prestacks:
\begin{equation} \label{equation-Cht-Gr-quotient}
\begin{aligned}
\epsilon_{G, N, I, \infty}: \on{Cht}_{G, N, I} & \rightarrow  [  G_{I, \infty} \backslash  \on{Gr}_{G, I}] \\
\left(  (x_i)_{i\in I}, (\mc G, \psi) \xrightarrow{\phi}  (\ta \mc G, \ta \psi) \right) & \mapsto  \left( (x_i)_{i\in I},  \restr{\mc G}{\Gamma_{\sum \infty x_i}} \xrightarrow{\phi}  \restr{\ta \mc G}{\Gamma_{\sum \infty x_i}}  \right) .
\end{aligned}
\end{equation}

\begin{rem}
The prestack $[G_{I, \infty} \backslash \on{Gr}_{G, I}]$ is not an inductive limit of algebraic stacks. But we can still use it for the construction in Sections \ref{subsection-para-induction-diagram} and \ref{subsection-quotient-by-Xi}.
We will construct a variant of morphism (\ref{equation-Cht-Gr-quotient}) for algebraic stacks in \ref{subsection-epsilon-G-N-I-d}.
\end{rem}

\quad

The following definition will be used in Section \ref{section-contractibility}.
\begin{defi} \label{def-Bun-G-N-I-d}
(a) We define $\Bun_{G, N, I, d}$ to be the prestack that associates to any affine scheme $S$ over $\Fq$ the groupoid classifying the following data:

(i) $(x_i)_{i \in I} \in (X\sm N)^I(S)$, 

(ii) $\mc G$: a $G$-bundle over $X \times S$,
 
(iii) a level structure on the divisor $(N \times S) + \Gamma_{\sum d x_i}$, i.e. an isomorphism of $G$-bundles:  
$\psi: \restr{\mc G}{ (N \times S) + \Gamma_{\sum d x_i} } \isom \restr{G}{ (N\times S) + \Gamma_{\sum d x_i}  }$

(b) We define $\Bun_{G, N, I, \infty}:= \underset{\longleftarrow}{\lim } \ \Bun_{G, N, I, d}$.
\end{defi}

\sssec{}
$\Bun_{G, N, I, d}$ is a smooth algebraic stack over $(X \sm N)^I$. 
Its fiber over a point $(x_i)_{i\in I} \in (X \sm N)^I(\Fq)$ is $\Bun_{G, N+\sum d x_i}$.

\quad

\sssec{}
The definitions and constructions in this subsection work for all affine smooth geometrically connected algebraic groups over $\Fq$ (we will use these for parabolic subgroups of $G$ and their Levi quotients).

\subsection{Parabolic induction diagrams}   \label{subsection-para-induction-diagram}

\sssec{}    \label{subsectino-def-Cht-P-Cht-M}
Let $P$ be a parabolic subgroup of $G$ and let $M$ be its Levi quotient. Applying the definitions and constructions in Section \ref{section-def-cht} to $P$ and $M$, respectively, we define $\Bun_{P, N}$, $\on{Cht}_{P, N, I}$, $\on{Gr}_{P, I}$, $P_{I, \infty}$, $\epsilon_{P, N, I, \infty}$ and $\Bun_{M, N}$, $\on{Cht}_{M, N, I}$, $\on{Gr}_{M, I}$, $M_{I, \infty}$, $\epsilon_{M, N, I, \infty}$.

\begin{rem}
When $N$ is non-empty, the prestack $\on{Cht}_{P, N, I}$ defined above is not the same as the one defined in \cite{var} 2.28. We will describe the difference in Remark \ref{rem-Cht-P-prime-is-fiber-product}.
\end{rem}

\sssec{} \label{subsection-Bun-G-Bun-P-Bun-M}
The morphisms of groups $G \hookleftarrow P \twoheadrightarrow M$ induce morphisms of prestacks over $\on{Spec} \Fq$:
\begin{equation}    \label{equation-Bun-G-P-M}
\Bun_{G, N} \xleftarrow{i^{Bun}} \Bun_{P, N} \xrightarrow{\pi^{Bun}} \Bun_{M, N}.
\end{equation}
%
%

\begin{construction}   \label{constr-Cht-G-Cht-P-Cht-M}
The morphisms of groups $G \hookleftarrow P \twoheadrightarrow M$ induce morphisms of prestacks over $(X \sm N)^I$: 
\begin{equation}   \label{diagram-Cht-G-P-M-general}
\xymatrixrowsep{1pc}
\xymatrixcolsep{1pc}
\xymatrix{
& \on{Cht}_{P, N, I} \ar[ld]_{i}  \ar[rd]^{\pi}  \ar[dd]_{\mf p_P} & \\
\on{Cht}_{G, N, I}    \ar[rd]_{\mf p_G}
& & \on{Cht}_{M, N, I}   \ar[ld]^{\mf p_M} \\
& (X \sm N)^I 
}
\end{equation}
More concretely, for any affine scheme $S$ over $\Fq$, 

$i: \on{Cht}_{P, N, I}(S) \rightarrow \on{Cht}_{G, N, I} (S)$ is given by $(\mc P \rightarrow \ta \mc P) \mapsto (\mc P \overset{P} \times G \rightarrow \ta \mc P \overset{P} \times G)$ where the level structure $\psi: \restr{\mc P}{N \times S} \isom \restr{P}{N \times S} $ is sent to $\psi \overset{P} \times G$;

$\pi: \on{Cht}_{P, N, I}(S) \rightarrow \on{Cht}_{M, N, I} (S)$ is given by $(\mc P \rightarrow \ta \mc P) \mapsto (\mc P \overset{P} \times M \rightarrow \ta \mc P \overset{P} \times M)$ where the level structure $\psi$ is sent to $\psi \overset{P} \times M$.
\end{construction}

\sssec{} 
The morphisms of groups $G \hookleftarrow P \twoheadrightarrow M$ induce morphisms of ind-schemes over $X^I$: 
\begin{equation}   \label{diagram-Gr-G-P-M}
\on{Gr}_{G, I} \xleftarrow{i^0} \on{Gr}_{P, I} \xrightarrow{\pi^0} \on{Gr}_{M, I}.
\end{equation}


\sssec{}  \label{subsection-morphism-of-quotient-stacks}
Let $\mc X$ (resp. $\mc Y$) be an (ind-)scheme over a base $S$ that is equipped with an action of a group scheme $A$ (resp. $B$) over $S$ from the right. 
Let $A \rightarrow B$ be a morphism of group schemes over $S$. Let $\mc X \rightarrow \mc Y$ be a morphism of (ind-)schemes over $S$ which is $A$-equivariant (where $A$ acts on $\mc Y$ via $A \rightarrow B$). This morphism induces a morphism of quotient prestacks
$$[A \backslash \mc X] \rightarrow [B \backslash \mc Y].$$

\sssec{}
Applying \ref{subsection-morphism-of-quotient-stacks} to $i^0: \on{Gr}_{P, I} \rightarrow \on{Gr}_{G, I}$ and $P_{I, \infty} \hookrightarrow G_{I, \infty}$, we obtain a morphism of prestacks:  $$\ov{i^0} :    [  P_{I, \infty} \backslash  \on{Gr}_{P, I}] \rightarrow [  G_{I, \infty} \backslash  \on{Gr}_{G, I}].$$

Applying \ref{subsection-morphism-of-quotient-stacks} to $\pi^0: \on{Gr}_{P, I} \rightarrow \on{Gr}_{M, I}$ and $P_{I, \infty} \twoheadrightarrow M_{I, \infty}$, we obtain a morphism of prestacks: $$ \ov{\pi^0} :  [ P_{I, \infty} \backslash  \on{Gr}_{P, I}] \rightarrow [M_{I, \infty} \backslash  \on{Gr}_{M, I}] .$$


\sssec{}
The following diagram of prestacks is commutative:
\begin{equation}   \label{diagram-TC-Cht-Gr}
\xymatrix{
\on{Cht}_{G, N, I} \ar[d]^{\epsilon_{G, N, I, \infty}} 
&\on{Cht}_{P, N, I}  \ar[l]_{i}   \ar[d]^{\epsilon_{P, N, I, \infty}}  \ar[r]^{\pi} 
&  \on{Cht}_{M, N, I}  \ar[d]^{\epsilon_{M, N, I, \infty}} \\
[G_{I, \infty} \backslash  \on{Gr}_{G, I}]       
& [ P_{I, \infty} \backslash  \on{Gr}_{P, I}]   \ar[l]_{\ov{i^0} }   \ar[r]^{\ov{\pi^0}} 
&  [  M_{I, \infty} \backslash  \on{Gr}_{M, I}]
}
\end{equation}

\subsection{Quotient by $\Xi$}   \label{subsection-quotient-by-Xi}

\sssec{}   \label{subsection-Xi-act-on-Cht-G}
Let $Z_G$ be the center of $G$ as defined in \ref{subsection-def-Z-G}. We have an action of $\Bun_{Z_G}$ on $\Bun_{G, N}$ by twisting a $G$-bundle by a $Z_G$-bundle, i.e. the action of $\mc T_Z \in \Bun_{Z_G}$ is given by $\mc G \mapsto (\mc G \times \mc T_Z) / Z_G$. Similarly, $\Bun_{Z_G}$ acts on $[  G_{I, \infty} \backslash  \on{Gr}_{G, I}] $, i.e. the action of $\mc T_Z \in \Bun_{Z_G}$ is given by 
$$\big( \mc G \xrightarrow{\phi} \mc G' \big) \mapsto \big( ( \mc G \times\restr{\mc T_G}{\Gamma_{\sum \infty x_i}} ) / Z_G \xrightarrow{\phi}   ( \mc G' \times\restr{\mc T_G}{\Gamma_{\sum \infty x_i}} ) / Z_G \big).$$

For $\mc T_Z \in \Bun_{Z_G}(\Fq)$, we have a canonical identification $\mc T_Z \simeq \ta \mc T_Z$. Thus $\Bun_{Z_G}(\Fq)$ acts on $\on{Cht}_{G, N, I}$ by twisting a $G$-bundle by a $Z_G$-bundle, i.e. the action of $\mc T_Z \in \Bun_{Z_G}(\Fq)$ is given by $(\mc G \xrightarrow{\phi}  \ta \mc G)  \mapsto \big( (\mc G \times \mc T_Z) / Z_G \xrightarrow{\phi}   \ta (\mc G \times \mc T_Z) / Z_G  \big) $.

The group $\Xi$ defined in \ref{subsection-def-Xi} acts on $\Bun_{G, N}$, $\on{Cht}_{G, N, I}$ and $[  G_{I, \infty} \backslash  \on{Gr}_{G, I}] $ via $\Xi \rightarrow Z_G(\mb A) \rightarrow  \Bun_{Z_G}(\Fq).$

\sssec{}   \label{subsection-Cht-Xi-to-G-ad}
Note that the morphism $\epsilon_{G, N, I, \infty}$ defined in (\ref{equation-Cht-Gr-quotient}) is $\Xi$-equivariant.

Now applying Definition \ref{defi-G-I-infry} to $Z_G$ (resp. $G^{\mr{ad}}$), we define a group scheme $(Z_G)_{I, \infty}$ (resp. $G^{\mr{ad}}_{I, \infty}$) over $X^I$. We have $G^{ad}_{I, \infty} = G_{I, \infty} / (Z_G)_{I, \infty}$.
The group scheme $(Z_G)_{I, \infty}$ acts trivially on $\on{Gr}_{G, I}$, so the action of $G_{I, \infty} $ on $\on{Gr}_{G, I}$ factors through $G^{\mr{ad}}_{I, \infty}$. We use this action to define the quotient prestack $[  G^{\mr{ad}}_{I, \infty} \backslash  \on{Gr}_{G, I}] $. The morphism $G_{I, \infty} \twoheadrightarrow G^{\mr{ad}}_{I, \infty}$ induces a morphism $[  G_{I, \infty} \backslash  \on{Gr}_{G, I}] \rightarrow [  G^{\mr{ad}}_{I, \infty} \backslash  \on{Gr}_{G, I}] $, which is $\Xi$-equivariant for the trivial action of $\Xi$ on $[  G^{\mr{ad}}_{I, \infty} \backslash  \on{Gr}_{G, I}] $.

Hence the composition of morphisms 
$$\Cht_{G, N, I}  \xrightarrow{\epsilon_{G, N, I, \infty}}  [  G_{I, \infty} \backslash  \Gr_{G, I}]  \rightarrow [  G_{I, \infty}^{\mr{ad}} \backslash  \Gr_{G, I}] $$ is $\Xi$-equivariant. 
Thus it factors through: 
\begin{equation}  \label{equation-epsilon-G-infty-Xi}
\epsilon_{G, N, I, \infty}^{\Xi}: \Cht_{G, N, I} / \Xi  \rightarrow  [  G_{I, \infty}^{ad} \backslash  \Gr_{G, I}] .
\end{equation}
We will construct a variant of morphism (\ref{equation-epsilon-G-infty-Xi}) for algebraic stacks in \ref{subsection-epsilon-G-N-I-d}.

\sssec{}
$Z_G$ acts on a $P$-bundle via $Z_G \hookrightarrow P$.
Just as in \ref{subsection-Xi-act-on-Cht-G}, we have an action of $\Bun_{Z_G}$ on $\Bun_{P, N}$ by twisting a $P$-bundle by a $Z_G$-bundle.
This leads to an action of $\Xi$ on $\Bun_{P, N}$, $\Cht_{P, N, I}$ and $[  P_{I, \infty} \backslash  \Gr_{P, I}] $ via $\Xi \rightarrow Z_G(\mb A) \rightarrow  \Bun_{Z_G}(\Fq)$.

Using the morphism $Z_G \hookrightarrow M$, we similarly obtain an action of $\Xi$ on $\Bun_{M, N}$, $\Cht_{M, N, I}$ and $[  M_{I, \infty} \backslash  \Gr_{M, I}]$. 

\sssec{}   \label{subsection-ov-P-ov-M}
Applying Definition \ref{defi-G-I-infry} to $\ov P:=P / Z_G$ (resp. $\ov M:=M / Z_G$), we define a group scheme $\ov P_{I, \infty} $ (resp. $\ov M_{I, \infty} $) over $X^I$. We have $\ov P_{I, \infty}= P_{I, \infty} / (Z_G)_{I, \infty}$ and $\ov M_{I, \infty}= M_{I, \infty} / (Z_G)_{I, \infty}$.

The morphism $\epsilon_{P, N, I, \infty}$ defined in \ref{subsectino-def-Cht-P-Cht-M} is $\Xi$-equivariant.
Since the group scheme $(Z_G)_{I, \infty}$ acts trivially on $\Gr_{P, I}$, the action of $P_{I, \infty}$ on $\Gr_{P, I}$ factors through $\ov P_{I, \infty}$. We denote by $[  \ov P_{I, \infty} \backslash  \Gr_{P, I}] $ the resulting quotient prestack.
The morphism $P_{I, \infty} \twoheadrightarrow \ov P_{I, \infty} $ induces a morphism $[  P_{I, \infty} \backslash  \Gr_{P, I}] \rightarrow [  \ov P_{I, \infty} \backslash  \Gr_{P, I}]$, which is $\Xi$-equivariant for the trivial action of $\Xi$ on $[  \ov P_{I, \infty} \backslash  \Gr_{P, I}] $.
Hence the composition of morphisms 
$\Cht_{P, N, I}  \xrightarrow{\epsilon_{P, N, I, \infty}}  [  P_{I, \infty} \backslash  \Gr_{P, I}]  \rightarrow [  \ov P_{I, \infty} \backslash  \Gr_{P, I}] $ is $\Xi$-equivariant. Thus it factors through:   
\begin{equation}
\epsilon_{P, N, I, \infty}^{\Xi}: \Cht_{P, N, I} / \Xi \rightarrow [\ov P_{I, \infty} \backslash \Gr_{P, I}].
\end{equation}

Similarly, the composition of morphisms $\Cht_{M, N, I}  \xrightarrow{\epsilon_{M, N, I, \infty}}  [  M_{I, \infty} \backslash  \Gr_{M, I}]  \rightarrow [  \ov M_{I, \infty} \backslash  \Gr_{M, I}] $ is $\Xi$-equivariant for the trivial action of $\Xi$ on $[  \ov M_{I, \infty} \backslash  \Gr_{M, I}] $. Thus it factors through: 
\begin{equation}
\epsilon_{M, N, I, \infty}^{\Xi}: \Cht_{M, N, I} / \Xi \rightarrow [\ov M_{I, \infty} \backslash \Gr_{M, I}].
\end{equation}

\sssec{}
The morphisms $i$ and $\pi$ in (\ref{diagram-TC-Cht-Gr}) are $\Xi$-equivariant. Diagram (\ref{diagram-TC-Cht-Gr}) induces a commutative diagram of prestacks:
\begin{equation}   \label{diagram-TC-Cht-Gr-Xi}
\xymatrix{
\Cht_{G, N, I}/ \Xi  \ar[d]^{\epsilon_{G, N, I, \infty}^{\Xi}} 
&\Cht_{P, N, I}   / \Xi  \ar[l]_i   \ar[d]^{\epsilon_{P, N, I, \infty}^{\Xi}}  \ar[r]^{\pi} 
&  \Cht_{M, N, I}  / \Xi  \ar[d]^{\epsilon_{M, N, I, \infty}^{\Xi}} \\
[G_{I, \infty}^{ad} \backslash  \Gr_{G, I}]       
& [ \ov P_{I, \infty} \backslash  \Gr_{P, I}]   \ar[l]_{\ov{i^0} }   \ar[r]^{\ov{\pi^0}} 
&  [  \ov M_{I, \infty} \backslash  \Gr_{M, I}]
}
\end{equation} 

\quad

\quad

{\it{ In the remaining part of Section 1, we introduce the Harder-Narasimhan stratification (compatible with the action of $\Xi$) for the parabolic induction diagram (\ref{diagram-Cht-G-P-M-general}). In order to do so, we use the Harder-Narasimhan stratification for the parabolic induction diagram (\ref{equation-Bun-G-P-M}).
From now on we work in the context of agebraic (ind-)stacks.

In Section \ref{subsection-HN-truncation}, we recall the usual Harder-Narasimhan stratification $\Bun_G^{\leq^G \mu} \subset \Bun_G$ and a variant $\Bun_G^{\leq^{G^{\mr{ad}}} \mu} \subset \Bun_G$ which is compatible with the action by $\Xi$.

In Section \ref{subsection-HN-Bun-M}, we introduce the Harder-Narasimhan stratification 
$\Bun_M^{\leq^{G^{\mr{ad}}} \mu} \subset \Bun_M$,
which allows us to construct in Section \ref{subsection-HN-Bun-G-P-M} the truncated parabolic induction diagrams (\ref{diagram-Bun-G-Bun-P-Bun-M-leq-ad-mu-Xi}):
\begin{equation*}  
\Bun_{G}^{\leq^{G^{\mr{ad}}} \mu} / \Xi \leftarrow \Bun_{P}^{\leq^{G^{\mr{ad}}} \mu} / \Xi \rightarrow \Bun_{M}^{\leq^{G^{\mr{ad}}} \mu} / \Xi .
\end{equation*}

In Section \ref{subsection-HN-Cht-G-P-M}, we define the Harder-Narasimhan stratification on the stacks of shtukas using Sections \ref{subsection-HN-truncation}-\ref{subsection-HN-Bun-G-P-M}.
 }}


\subsection{Harder-Narasimhan stratification of $\Bun_G$}   \label{subsection-HN-truncation}

In \ref{subsection-def-Lambda-G}-\ref{subsection-Bun-G-leq-mu-union-de-strate}, we recall the Harder-Narasimhan stratification of $\Bun_G$ defined in \cite{schieder} and \cite{DG15} Section 7. (In these papers, the group is reductive over an algebraically closed field. Since our group $G$ is split over $\Fq$, we use Galois descent to obtain the stratification over $\Fq$.)

In \ref{subsection-Lambda-G-ad}-\ref{rem-proof-Var-type-fini}, we recall a variant of the Harder-Narasimhan stratification of $\Bun_G$ which is compatible with the quotient by $\Xi$, as in \cite{var} Section 2 and \cite{vincent} Section 1.

\sssec{}   \label{subsection-def-Lambda-G}
Applying \ref{subsection-def-Lambda-H} to group $G$, we define $\wh \Lambda_G$, $\wh \Lambda_G^+$, $\wh \Lambda_G^{\mr{pos}}$, $\wh \Lambda_G^{\Q}$, $\wh \Lambda_G^{+, \Q}$, $\wh \Lambda_G^{\mr{pos}, \Q}$ and the partial order "$\leq^G$" on $\wh{\Lambda}_{G}^{\Q}$.

\sssec{} (\cite{schieder} 2.1.2)  \label{subsection-def-Lambda-G-P}
Let $P$ be a parabolic subgroup of $G$ and $M$ its Levi quotient. 
Consider the sublattice $\wh\Lambda_{[M, M]_{\mr{sc}}} \subset \wh \Lambda_G$ spanned by the simple coroots of $M$.
We define
\begin{equation}
\wh \Lambda_{G, P}:= \wh \Lambda_G / \wh\Lambda_{[M, M]_{\mr{sc}}}.
\end{equation}
Let $\wh \Lambda_{G, P}^{\Q}:=\wh \Lambda_{G, P} \otimes_{\Z} \Q$. We denote by $\wh \Lambda_{G, P}^{\mr{pos}}$ the image of $\wh \Lambda_G^{\mr{pos}}$ in $\wh \Lambda_{G, P}$, and by $\wh \Lambda_{G, P}^{\mr{pos}, \Q}$ the image of $\wh \Lambda_G^{\mr{pos}, \Q}$ in $\wh \Lambda_{G, P}^{\Q}$.
We introduce the partial order on $\wh{\Lambda}_{G, P}$ by $$\mu_1 \leq^G \mu_2 \Leftrightarrow \mu_2 - \mu_1 \in \wh \Lambda_{G, P}^{\mr{pos}}.$$

\sssec{} (\cite{schieder} 2.1.3, \cite{DG15} 7.1.3, 7.1.5)
Let $Z_M$ be the center of $M$. Let $\wh \Lambda_{Z_M}$ be the coweight lattice of $Z_M$, i.e. $\on{Hom}(\mb G_m, Z_M)$. Note that it equals to $\wh \Lambda_{Z_M^0} = \on{Hom}(\mb G_m, Z_M^0)$, where $Z_M^0$ is the neutral connected component of $Z_M$.

We have a natural inclusion $\wh \Lambda_{Z_M} \subset \wh \Lambda_G$ (because $Z_M$ is included in the image of $B \hookrightarrow P \twoheadrightarrow M$). The composition $\wh \Lambda_{Z_M}^{\Q} \hookrightarrow \wh \Lambda_G^{\Q} \twoheadrightarrow \wh \Lambda_{G, P}^{\Q}$ is an isomorphism: 
\begin{equation}
\wh \Lambda_{Z_M}^{\Q}  \isom \wh \Lambda_{G, P}^{\Q} .
\end{equation}

We define the slope map to be the composition
\begin{equation}    \label{equation-def-phi-p}
\phi_P: \wh \Lambda_{G, P} \rightarrow \wh \Lambda_{G, P}^{\Q} \cong \wh \Lambda_{Z_M}^{\Q} \hookrightarrow \wh \Lambda_G^{\Q} .
\end{equation} 

We define $\on{pr}_P$ to be the composition
\begin{equation}   \label{equation-def-pr-p}
\on{pr}_P: \wh \Lambda_G^{\Q} \twoheadrightarrow \wh \Lambda_{G, P}^{\Q} \simeq \wh \Lambda_{Z_M}^{\Q}.
\end{equation}

By definition, we have $\wh \Lambda_{G, G}^{\Q} = \wh \Lambda_{Z_G}^{\Q}$, $\wh \Lambda_{G, P} = \wh \Lambda_{M, M}$ and $\wh \Lambda_{G, B}= \wh \Lambda_{G}$. So $\phi_B$ is just the inclusion $\wh \Lambda_G \hookrightarrow \wh \Lambda_G^{\Q}$.

\begin{lem} (\cite{schieder} Proposition 3.1)  \label{lem-slope-map-preserves-order}
The slope map $\phi_P$ preserves the partial orders $"\leq^G"$ on $\wh{\Lambda}_{G, P}$ and $\wh{\Lambda}_{G}^{\Q}$ in the sense that it maps $\wh \Lambda_{G, P}^{\mr{pos}}$ to $\wh \Lambda_G^{\mr{pos}, \Q}$.
\cqfd
\end{lem}

\sssec{} (\cite{var} Lemma 2.2, \cite{schieder} 2.2.1, 2.2.2, \cite{DG15} 7.2.3)  \label{subsection-deg-M}
The map $\Bun_P \rightarrow \Bun_M$ in \ref{subsection-Bun-G-Bun-P-Bun-M} induces a bijection on the set of connected components of $\Bun_P $ and $\Bun_M$. We have $\pi_0(\Bun_P) \cong \pi_0(\Bun_M) \cong \wh \Lambda_{G, P}$. 
Let $\deg_M: \Bun_M \rightarrow \pi_0(\Bun_M)  \cong \wh \Lambda_{G, P}$ and $\deg_P: \Bun_P \rightarrow \Bun_M  \rightarrow \wh \Lambda_{G, P}$.

\begin{defi} (\cite{DG15} 7.3.3, 7.3.4)  \label{def-Bun-G-leq-mu}
For any $\mu\in \wh{\Lambda}_{G}^{+, \Q}$, we define $\Bun_G^{\leq^G \mu}$ to be the stack that associates to any affine scheme $S$ over $\Fq$ the groupoid 
\begin{gather*} 
   \Bun_{G}^{\leq^G  \mu}(S) := \{\mc G\in \Bun_{G}(S) |
    \text{ for each geometric point } s\in S, \text{ each parabolic }     \\ \nonumber
 \text{ subgroup }  P  \text{ and each } P \text{-structure } \mc P \text{ of } \mc G_s, 
   \text{ we have } \phi_P \circ \on{deg}_P(\mc P) \leq^G \mu \},
\end{gather*}
where a $P$-structure of $\mc G_s$ is a $P$-bundle $\mc P$ on $X_s$ such that $\mc P \overset{P}{\times} G \simeq \mc G_s$.
\end{defi}

\begin{rem}  \label{rem-Bun-G-leq-mu-equivalent-P-B}
(a) By \cite{schieder} Lemma 3.3, 
the above Definition \ref{def-Bun-G-leq-mu} is equivalent to 
\begin{gather*} 
   \Bun_{G}^{\leq^G  \mu}(S) := \{\mc G\in \Bun_{G}(S) |
    \text{ for each geometric point } s\in S,    \\ \nonumber
     \text{ and each } B \text{-structure } \mc B \text{ of } \mc G_s, 
   \text{ we have } \on{deg}_B(\mc B) \leq^G \mu \}. 
\end{gather*}
(the argument repeats the proof in $loc.cit.$ by replacing $\phi_G(\check{\lambda}_G)$ by $\mu$).

(b) 
By \cite{schieder} Proposition 3.2 and Remark 3.2.4, 
the definition of $\Bun_G^{\leq^G \mu}$ in (a) is equivalent to the Tannakian description:
\begin{gather*} 
   \Bun_{G}^{\leq^G  \mu}(S) := \{\mc G\in \Bun_{G}(S) |
    \text{ for each geometric point } s\in S,   \\ \nonumber
     \text{ each } B \text{-structure } \mc B \text{ of } \mc G_s  \text{ and each } \lambda \in \Lambda_{G},
   \text{ we have } \on{deg} \mc B_{\lambda} \leq \langle \mu, \lambda \rangle \}. 
\end{gather*}
where $\mc B_{\lambda}$ is the line bundle associated to $\mc B$ and $B  \rightarrow T \xrightarrow{\lambda} \mb G_m$.

(c) The reason why we use Definition \ref{def-Bun-G-leq-mu} (rather than its equivalent forms) is that it will be useful for non split groups in future works.
\end{rem}

\begin{lem}   \label{lem-Bun-G-leq-mu-open}
(\cite{DG15} 7.3.4, Proposition 7.3.5)

(a) For any $\mu\in \wh{\Lambda}_{G}^{+, \Q}$, the stack $\Bun_G^{\leq^{G}  \mu}$ is an open substack of $\Bun_{G}$. 

(b) For any $\mu_1 \leq^G \mu_2$, we have an open immersion $\Bun_G^{\leq^{G}  \mu_1} \hookrightarrow \Bun_G^{\leq^{G}  \mu_2}.$ 

(c) $\Bun_{G} = \bigcup_{\mu \in \wh{\Lambda}_{G}^{+, \Q}} \Bun_G^{\leq^{G}  \mu}.$

(d) The open substack $\Bun_G^{\leq^{G}  \mu}$ is of finite type.
\cqfd
\end{lem}

\begin{defi}  \label{def-Bun-G-equal-mu}
For any $\lambda \in \wh \Lambda_G^{+, \Q}$, let $\Bun_G^{(\lambda)} \subset \Bun_G$ be the quasi-compact locally closed reduced substack defined in \cite{schieder} Theorem 2.1 and \cite{DG15} Theorem 7.4.3. It is called a Harder-Narasimhan stratum of $\Bun_G$. 
\end{defi}

\sssec{}   \label{subsection-Bun-G-leq-mu-union-de-strate}
(\cite{DG15} Corollary 7.4.5) We have
$$\Bun_G^{(\lambda)} \neq \emptyset   \Rightarrow   \lambda \in \bigcup_{P \subsetneq G} \iota \circ \on{pr}_P( \wh \Lambda_G ),$$
where $\on{pr}_P$ is defined in (\ref{equation-def-pr-p}) and $\iota: \wh \Lambda_{Z_M}^{\Q} \hookrightarrow \wh \Lambda_G^{\Q}$ is the inclusion.
For any $\mu \in \wh \Lambda_G^{+, \Q}$, we have 
$$\Bun_G^{\leq^G \mu} = \bigcup_{\lambda \in \wh \Lambda_G^{+, \Q}, \; \lambda \leq^G \mu} \Bun_G^{(\lambda)} .$$
The set $\{ \lambda \in \wh \Lambda_G^{+, \Q} \, | \, \lambda \leq^G \mu \text{ and } \Bun_G^{(\lambda)} \neq \emptyset \}$ is finite. This gives another proof of Lemma \ref{lem-Bun-G-leq-mu-open} (d).

\quad

The above open substack $\Bun_G^{\leq^G \mu}$ is not preserved by the action of $\Xi$ on $\Bun_G$. Now we introduce open substacks which are preserved by the action of $\Xi$.

\sssec{}    \label{subsection-Lambda-G-ad}
Applying \ref{subsection-def-Lambda-H} to group $G^{\mr{ad}}$, we define $\wh \Lambda_{G^{\mr{ad}}}$, $\wh \Lambda_{G^{\mr{ad}}}^+$, $\wh \Lambda_{G^{\mr{ad}}}^{\mr{pos}}$, $\wh \Lambda_{G^{\mr{ad}}}^{\Q}$, $\wh \Lambda_{G^{\mr{ad}}}^{+, \Q}$, $\wh \Lambda_{G^{\mr{ad}}}^{\mr{pos}, \Q}$ and the partial order "$ \leq^{G^{\mr{ad}} } $" on $\wh \Lambda_{G^{\mr{ad}}}$.


The morphism $G \twoheadrightarrow G/ Z_G = G^{\mr{ad}}$ induces a morphism
\begin{equation}   \label{equation-Upsilon-G}
\Upsilon_G: \wh{\Lambda}_G^{\Q} \rightarrow \wh{\Lambda}_{G^{\mr{ad}}}^{\Q}.
\end{equation}
It maps $\wh \Lambda_G^{\mr{pos}, \Q}$ to $\wh \Lambda_{G^{\mr{ad}}}^{\mr{pos}, \Q}$.

\begin{defi} \label{def-Bun-G-leq-ad-mu}
For any $\mu\in \wh{\Lambda}_{G^{\mr{ad}}}^{+, \Q}$, we define $\Bun_G^{\leq^{G^{\mr{ad}} } \mu}$ to be the stack that associates to any affine scheme $S$ over $\Fq$ the groupoid 
\begin{gather*} 
   \Bun_{G}^{\leq^{G^{\mr{ad}} }  \mu}(S) := \{\mc G\in \Bun_{G}(S) |
    \text{ for each geometric point } s\in S, \text{ each parabolic }     \\ \nonumber
   \text{ subgroup }  P  \text{ and each } P \text{-structure } \mc P \text{ on } \mc G_s, 
   \text{ we have } \Upsilon_G \circ \phi_P \circ \on{deg}_P(\mc P) \leq^{G^{\mr{ad}} }  \mu \}. 
\end{gather*}
\end{defi}

\begin{rem}
For the same reason as in Remark \ref{rem-Bun-G-leq-mu-equivalent-P-B}, Definition \ref{def-Bun-G-leq-ad-mu} is equivalent to \cite{var} Notation 2.1 b) and \cite{vincent} (1.3).
\end{rem}

\sssec{}    \label{subsection-Bun-G-leq-ad-mu-union-strate}
Just as in \ref{subsection-Bun-G-leq-mu-union-de-strate}, for $\mu \in \wh \Lambda_{G^{\mr{ad}}}^{+, \Q}$, we have 
$$\Bun_G^{\leq^{G^{\mr{ad}}} \mu} = \bigcup_{\lambda \in \wh \Lambda_G^{+, \Q}, \; \Upsilon_G(\lambda) \leq^{G^{\mr{ad}}} \mu} \Bun_G^{(\lambda)} .$$
The set $\{ \lambda \in \wh \Lambda_G^{+, \Q} \, | \, \Upsilon_G(\lambda) \leq^{G^{\mr{ad}}}  \mu \text{ and } \Bun_G^{(\lambda)} \neq \emptyset \}$ is finite modulo $\wh \Lambda_{Z_G}$. 

\sssec{}   \label{subsection-Bun-G-leq-ad-mu-Xi}
The action of $\Xi$ on $\Bun_G$ preserves $\Bun_{G}^{\leq^{G^{\mr{ad}} } \mu}$. We define the quotient $\Bun_G^{ \leq^{G^{\mr{ad}} } \mu} / \Xi$.

\begin{lem}   \label{lem-Bun-G-leq-ad-mu-open} 
(a) For any $\mu\in \wh{\Lambda}_{G^{\mr{ad}}}^{+, \Q}$, the stack $\Bun_G^{\leq^{G^{\mr{ad}} } \mu}$ is an open substack of $\Bun_{G}$. 

(b) For any $\mu_1 \leq^{G^{\mr{ad}} } \mu_2$, we have an open immersion $\Bun_G^{\leq^{G^{\mr{ad}} } \mu_1} \hookrightarrow \Bun_G^{\leq^{G^{\mr{ad}} } \mu_2}.$ 

(c) $\Bun_{G}$ is the inductive limit of these open substacks: $\Bun_{G} = \bigcup_{\mu \in \wh{\Lambda}_{G^{\mr{ad}}}^{+, \Q}} \Bun_G^{\leq^{G^{\mr{ad}} } \mu}.$

(d) The stacks $\Bun_{G}^{\leq^{G^{\mr{ad}} }  \mu}  / \Xi$ is of finite type.
\end{lem}
%
%
\dem
(a), (b) and (c) are induced by Lemma \ref{lem-Bun-G-leq-mu-open} (see also \cite{var} Lemme A.3). (d) follows from \ref{subsection-Bun-G-leq-ad-mu-union-strate}. 
%
%
%
\cqfd

\begin{rem}   \label{rem-proof-Var-type-fini}
See \cite{var} Lemma 3.1 and 3.7 for another proof of Lemma \ref{lem-Bun-G-leq-mu-open} (d) and Lemma \ref{lem-Bun-G-leq-ad-mu-open} (d). 
\end{rem}


\subsection{Harder-Narasimhan stratification of $\Bun_M$}     \label{subsection-HN-Bun-M}

Let $P$ be a proper parabolic subgroup of $G$ and $M$ its Levi quotient. 



\sssec{}
Applying \ref{subsection-def-Lambda-H} to group $M$, we define $\wh \Lambda_M$, $\wh \Lambda_M^+$, $\wh \Lambda_M^{\mr{pos}}$, $\wh \Lambda_M^{\Q}$, $\wh \Lambda_M^{+, \Q}$, $\wh \Lambda_M^{\mr{pos}, \Q}$
and the partial order "$\leq^M$" on $\wh{\Lambda}_{M}^{\Q}$.

\sssec{}
The Sections \ref{subsection-def-Lambda-G-P}-\ref{subsection-Bun-G-leq-mu-union-de-strate} work also for $M$. 
In particular, 
let $P'$ be a parabolic subgroup of $M$, we have the slope map $\phi_{P'} : \wh \Lambda_{M, P'} \rightarrow \wh \Lambda_M^{\Q}$ and $\deg_{P'}: \Bun_{P'} \rightarrow \wh \Lambda_{M, P'}$.
%
%
%

\begin{defi}   \label{def-Bun-M-equal-lambda}
Applying Definition \ref{def-Bun-G-equal-mu} to $M$, for any $\lambda \in \wh \Lambda_M^{+, \Q}$, we define a quasi-compact locally closed substack $\Bun_M^{(\lambda)} \subset \Bun_M$, called a Harder-Narasimhan stratum of $\Bun_M$. 
\end{defi}

Now we introduce $\Bun_M^{\leq^{G^{\mr{ad}}} \mu} \subset \Bun_M$ which will be used to construct diagram (\ref{diagram-Bun-G-Bun-P-Bun-M-leq-ad-mu-Xi}). 

\begin{defi} \label{def-troncature-Bun-M-leq-ad-mu}
For any $\mu\in \wh{\Lambda}_{G^{\mr{ad}}}^{+, \Q}$, we define $\Bun_M^{\leq^{G^{\mr{ad}} } \mu}$ to be the stack that associates to any affine scheme $S$ over $\Fq$ the groupoid $\Bun_{M}^{\leq^{G^{\mr{ad}}}  \mu}(S) :=$
\begin{gather*}   
   \{\mc M \in \Bun_{M}(S) |
    \text{ for each geometric point } s\in S, \text{ each parabolic subgroup } P'  \\ \nonumber
   \text{ of } M   \text{ and each } P' \text{-structure } \mc P' \text{ of } \mc M_s,  
   \text{ we have } \Upsilon_G \circ \phi_{P'} \circ \on{deg}_{P'}(\mc P') \leq^{G^{\mr{ad}}} \mu \}. 
\end{gather*}
where $\Upsilon_G: \wh{\Lambda}_{M}^{\Q} = \wh{\Lambda}_G^{\Q} \rightarrow \wh{\Lambda}_{G^{\mr{ad}}}^{\Q}$ is defined in (\ref{equation-Upsilon-G}).
\end{defi}

Similarly to Lemma \ref{lem-Bun-G-leq-ad-mu-open}, we have 
\begin{lem}   \label{lem-Bun-M-leq-ad-mu-open}

(a) For any $\mu\in \wh{\Lambda}_{G^{\mr{ad}}}^{+, \Q}$, the stack $\Bun_{M}^{\leq^{G^{\mr{ad}}} \mu}$ is an open substack of $\Bun_{M}$. 

(b) For any $\mu_1 \leq^{G^{\mr{ad}}} \mu_2$, we have an open immersion $\Bun_{M}^{\leq^{G^{\mr{ad}}} \mu_1} \hookrightarrow \Bun_{M}^{\leq^{G^{\mr{ad}}} \mu_2}.$ 

(c) $\Bun_{M}$ is the inductive limit of these open substacks: $\Bun_{M} = \bigcup_{\mu \in \wh{\Lambda}_{G^{\mr{ad}}}^{+, \Q}} \Bun_{M}^{\leq^{G^{\mr{ad}}} \mu}.$
\cqfd
\end{lem}

\sssec{}    \label{subsection-Xi-not-lattice-in-Z-M-for-Bun-M}
The action of $\Xi$ on $\Bun_M$ preserve $\Bun_M^{ \leq^{G^{\mr{ad}}} \mu}$. We define the quotient $\Bun_M^{ \leq^{G^{\mr{ad}}} \mu} / \Xi$. Note that $\Xi$ is a lattice in $Z_G(F) \backslash Z_G(\mb A)$. However, $\Xi$ is only a discrete subgroup but not a lattice in $Z_M(F) \backslash Z_M(\mb A)$ (since $P \subsetneq G$). 
We will see that $\Bun_M^{ \leq^{G^{\mr{ad}}} \mu} / \Xi$ is locally of finite type but not necessarily of finite type.

\sssec{}    \label{subsection-def-Bun-M-nu}
Note that $\wh \Lambda_{G, P} = \wh \Lambda_{M, M}$. Consider the composition of morphisms 
\begin{equation}    \label{equation-Bun-M-deg-to-Lambda-Z-M-Z-G}
\Bun_M \xrightarrow{\deg_M} \wh \Lambda_{M, M} \rightarrow \wh \Lambda_{M, M}^{\Q} \simeq \wh \Lambda_{Z_M}^{\Q} \twoheadrightarrow \wh \Lambda_{Z_M / Z_G}^{\Q} ,
\end{equation}
where $\deg_M$ is defined in \ref{subsection-deg-M}.
We denote by $A_M$ the image of $\wh \Lambda_{M, M}$ in $\wh \Lambda_{Z_M / Z_G}^{\Q}$.
For any $\nu \in \wh \Lambda_{Z_M / Z_G}^{\Q}$, we denote by $\Bun_M^{\nu}$ its inverse image in $\Bun_M$. It is non-empty if and only if $\nu \in A_M$.

\begin{defi}
We define $\Bun_M^{\leq^{G^{\mr{ad}}}  \mu, \, \nu}$ to be the intersection of $\Bun_M^{\leq^{G^{\mr{ad}} }  \mu}$ and $\Bun_M^{\nu}$. 
\end{defi}

\sssec{}    \label{subsection-def-Bun-M-leq-mu-nu}
$\Bun_M^{\leq^{G^{\mr{ad}}}  \mu, \, \nu}$ is open and closed in $\Bun_M^{\leq^{G^{\mr{ad}}} \mu}$ and is open in $\Bun_M^{\nu}$. 
We have a decomposition
\begin{equation}   \label{equation-decom-Bun-M-leq-G-ad-mu}
\Bun_M^{ \leq^{G^{\mr{ad}}}  \mu} = \underset{ \nu \in \wh{\Lambda}_{Z_M / Z_G}^{\Q} } \bigsqcup \Bun_M^{ \leq^{G^{\mr{ad}}} \mu, \, \nu}. 
\end{equation}

\sssec{}   \label{subsection-Bun-M-leq-mu-union-of-Bun-M-lambda}
Just as in \ref{subsection-Bun-G-leq-ad-mu-union-strate}, we have $$\Bun_M^{\leq^{G^{\mr{ad}}}  \mu} =  \underset{ \lambda \in \wh{\Lambda}_{M}^{+, \Q}, \; \Upsilon_G(\lambda) \leq^{G^{\mr{ad}}} \mu }  \bigcup \Bun_M^{(\lambda) }. $$

\sssec{}   \label{subsection-Bun-M-lambda-in-which-component}
Similarly to (\ref{equation-def-pr-p}), we define
\begin{equation}     \label{equation-pr-p-ad}
\on{pr}_P^{\mr{ad}}: \wh \Lambda_{G^{\mr{ad}}}^{\Q} \rightarrow\wh \Lambda_{Z_M / Z_G}^{\Q} .
\end{equation}

Taking into account that $\wh \Lambda_G = \wh \Lambda_M$ and $\wh \Lambda_{G, P} = \wh \Lambda_{M, M}$, for any $\lambda \in \wh{\Lambda}_{M}^{+, \Q},$ we deduce that $\Bun_M^{(\lambda)} \subset \Bun_M^{\nu}$ if and only if $\nu = \on{pr}_P^{\mr{ad}} \circ \Upsilon_G(\lambda)$. 

\sssec{}    \label{subsection-Bun-M-leq-ad-mu-nu-union-of-Bun-M-lambda}
We deduce from \ref{subsection-Bun-M-leq-mu-union-of-Bun-M-lambda} and \ref{subsection-Bun-M-lambda-in-which-component} that
\begin{equation}      
\Bun_M^{\leq^{G^{\mr{ad}}} \mu, \, \nu} =  \underset{ \lambda \in \wh{\Lambda}_{M}^{+, \Q}, \; \Upsilon_G(\lambda) \leq^{G^{\mr{ad}}} \mu, \; \on{pr}_P^{\mr{ad}} \circ \Upsilon_G(\lambda) = \nu }  \bigcup \Bun_M^{(\lambda) }. 
\end{equation}

\sssec{}      \label{subsection-leq-in-Lambda-Z-M-Z-G}
We denote by $\wh{\Lambda}_{Z_M / Z_G}^{\on{pos}, \Q} := \on{pr}_P^{\mr{ad}} ( \wh \Lambda_{G^{\mr{ad}} }^{\mr{pos}, \Q} )$. 
We introduce the partial order on $\wh{\Lambda}_{Z_M / Z_G}^{\Q}$ by $$\mu_1 \leq^{G^{\mr{ad}} } \mu_2 \Leftrightarrow \mu_2 - \mu_1 \in \wh{\Lambda}_{Z_M / Z_G}^{\on{pos}, \Q}.$$
By definition, for $\check{\gamma} \in \wh \Gamma_M$, we have $\on{pr}_P^{ad} \circ \Upsilon_G (\check{\gamma}) =0.$ 
By \cite{schieder} Proposition 3.1, for $\check{\gamma} \in \wh \Gamma_G - \wh \Gamma_M$ we have $\on{pr}_P^{ad} \circ \Upsilon_G (\check{\gamma}) > 0$ and these $\on{pr}_P^{ad} \circ \Upsilon_G (\check{\gamma}) $ are linearly independent. 
Thus for $\lambda_1, \lambda_2 \in \wh{\Lambda}_{ G^{\mr{ad}}  }^{\Q}$ and $\lambda_1 \leq^{G^{\mr{ad}} } \lambda_2 $, we have $\on{pr}_P^{\mr{ad}} (\lambda_1) \leq^{G^{\mr{ad}} } \on{pr}_P^{\mr{ad}} (\lambda_2) $.
And the inclusion $\wh{\Lambda}_{Z_M / Z_G}^{\Q} \subset \wh{\Lambda}_{G^{\mr{ad}}}^{\Q}$ maps $\wh{\Lambda}_{Z_M / Z_G}^{\on{pos}, \Q}$ to $\wh{\Lambda}_{G^{\mr{ad}}}^{\on{pos}, \Q}$.

\begin{lem}   \label{lem-Bun-M-leq-ad-mu-nu-non-vide}
Let $\mu \in \wh{\Lambda}_{ G^{\mr{ad}}  }^{+, \Q}$. Then the stack $\Bun_M^{ \leq^{G^{\mr{ad}}} \mu , \; \nu}$ is empty unless $\nu \in A_M$ defined in \ref{subsection-def-Bun-M-nu} and $\nu \leq^{G^{\mr{ad}} }   \on{pr}_P^{\mr{ad}}(\mu)$.
\end{lem}
\dem
The first condition follows from \ref{subsection-def-Bun-M-nu}.
To prove the second condition, note that for the set $\{  \lambda \in \wh{\Lambda}_{M}^{+, \Q} \; | \; \Upsilon_G(\lambda) \leq^{G^{\mr{ad}}} \mu, \; \on{pr}_P^{\mr{ad}} \circ \Upsilon_G(\lambda) = \nu   \}$ to be non-empty, by \ref{subsection-leq-in-Lambda-Z-M-Z-G} we must have $\nu  \leq^{G^{\mr{ad}} } \on{pr}_P^{\mr{ad}} (\mu)$.
\cqfd

\sssec{}
Let $\ov M= M / Z_G$ as in \ref{subsection-ov-P-ov-M}. For $\lambda, \mu \in \wh{\Lambda}_{G^{\mr{ad}}}^{\Q}$, we define $\lambda \leq^{\ov M} \mu$ if and only if $\mu - \lambda$ is a linear combination of simple coroots of $M$ with coefficients in $\Q_{\geq 0}$ modulo $\wh{\Lambda}_{Z_G}^{\Q}$.

\sssec{}   \label{subsection-pr-p-ad-gamma-geq-0}


Let $\lambda, \mu \in \wh{\Lambda}_{G^{\mr{ad}}}^{\Q}$ and $\lambda \leq^{G^{\mr{ad}}} \mu$. We write $\lambda = \mu - \sum_{\check{\gamma} \in \wh \Gamma_G} c_{\check{\gamma}} \Upsilon_G( \check{\gamma})$ for some $c_{\check{\gamma}} \in \Q_{\geq 0} $. We deduce from \ref{subsection-leq-in-Lambda-Z-M-Z-G} that $\on{pr}_P^{ad}(\lambda) = \on{pr}_P^{ad}(\mu)$ if and only if $c_{\check{\gamma}} =0$ for all $\check{\gamma} \in \wh \Gamma_G - \wh \Gamma_M$.
Hence 
\begin{equation}    \label{equation-leq-M-bar-equal-leq-mu-and-pr-same}
\lambda \leq^{G^{\mr{ad}}} \mu  \text{ and } \on{pr}_P^{ad}(\lambda) = \on{pr}_P^{ad}(\mu) \; \Leftrightarrow \; \lambda \leq^{\ov M} \mu. 
\end{equation}

\sssec{}   \label{subsection-leq-mu-fix-nu-equal-leq-M-mu-nu}
Let $\mu \in \wh{\Lambda}_{ G^{\mr{ad}}  }^{+, \Q}$ and $\nu \leq^{G^{\mr{ad}} }   \on{pr}_P^{\mr{ad}}(\mu)$. For every $\check{\gamma} \in \wh \Gamma_G - \wh \Gamma_M$, let $c_{\gamma} \in \Q_{\geq 0}$ be the unique coefficient such that $$\on{pr}_P^{\mr{ad}} (\mu ) - \sum_{\check{\gamma} \in \wh \Gamma_G - \wh \Gamma_M} c_\gamma \on{pr}_P^{\mr{ad}} \circ \Upsilon_G ( \check{\gamma} ) = \nu.$$
We define $\mu_{\nu} := \mu - \sum_{\check{\gamma} \in \wh \Gamma_G - \wh \Gamma_M} c_\gamma \Upsilon_G(\check{\gamma} )$. 
As in \ref{subsection-pr-p-ad-gamma-geq-0}, we deduce that 
\begin{equation}      \label{equation-leq-G-ad-fix-nu-equal-leq-M-mu-nu}
\lambda \leq^{G^{\mr{ad}}} \mu  \text{ and } \on{pr}_P^{ad}(\lambda) = \nu \; \Leftrightarrow  \; \lambda \leq^{\ov M} \mu_{\nu}. 
\end{equation}


\sssec{}
The action of $\Xi$ on $\Bun_M$ preserves $\Bun_M^{ \leq^{G^{\mr{ad}}} \mu, \; \nu} $. We define the quotient $\Bun_M^{ \leq^{G^{\mr{ad}}} \mu, \; \nu} / \Xi$.

\begin{lem}    \label{lem-Bun-M-leq-mu-nu-Xi-finite-type}
$\Bun_M^{ \leq^{G^{\mr{ad}}} \mu, \; \nu} / \Xi$ is of finite type.
\end{lem}
\dem
By (\ref{equation-leq-G-ad-fix-nu-equal-leq-M-mu-nu}), we have $$\{  \lambda \in \wh{\Lambda}_{M}^{+, \Q} \; | \; \Upsilon_G(\lambda) \leq^{G^{\mr{ad}}} \mu, \; \on{pr}_P^{\mr{ad}} \circ \Upsilon_G(\lambda) = \nu  \} = \{  \lambda \in \wh{\Lambda}_{M}^{+, \Q} \; | \;    \Upsilon_G(\lambda) \leq^{\ov{M}} \mu_{\nu}     \} .$$
We deduce from \ref{subsection-Bun-G-leq-mu-union-de-strate} (applied to $M$) that the set $\{  \lambda \in \wh{\Lambda}_{M}^{+, \Q} \; | \; \Upsilon_G(\lambda) \leq^{G^{\mr{ad}}} \mu, \; \on{pr}_P^{\mr{ad}} \circ \Upsilon_G(\lambda) = \nu, \;  \Bun_M^{( \lambda)} \neq \emptyset  \}$ is finite modulo $\wh \Lambda_{Z_G}$.
By Definition \ref{def-Bun-M-equal-lambda}, $\Bun_M^{(\lambda)}$ is of finite type. 
From \ref{subsection-Bun-M-leq-ad-mu-nu-union-of-Bun-M-lambda} we deduce the lemma. 
\cqfd

\sssec{}    \label{subsection-def-translated-cone-ad}
By Lemma \ref{lem-Bun-M-leq-ad-mu-nu-non-vide}, the decomposition (\ref{equation-decom-Bun-M-leq-G-ad-mu}) is in fact indexed by a translated cone in $\wh{\Lambda}_{Z_M / Z_G}^{\Q}$:
\begin{equation}  
\wh{\Lambda}_{Z_M / Z_G}^{\mu} := \{  \nu \in \wh{\Lambda}_{Z_M / Z_G}^{\Q}, \; \nu \leq^{G^{\mr{ad}} }  \on{pr}_P^{\mr{ad}}(\mu) \} .
\end{equation}
We deduce that
\begin{equation}   \label{formule-decom-Bun-M-leq-ad-mu-cone}
\Bun_M^{ \leq^{G^{\mr{ad}}} \mu}= \underset{ \nu \in \wh{\Lambda}_{Z_M / Z_G}^{\mu} } \bigsqcup \Bun_M^{ \leq^{G^{\mr{ad}}} \mu, \; \nu}.
\end{equation}
and
\begin{equation}   \label{formule-decom-Bun-M-leq-ad-mu-cone-Xi}
\Bun_M^{ \leq^{G^{\mr{ad}}} \mu} / \Xi = \underset{ \nu \in \wh{\Lambda}_{Z_M / Z_G}^{\mu} } \bigsqcup \Bun_M^{ \leq^{G^{\mr{ad}}} \mu, \; \nu} / \Xi.
\end{equation}

\subsection{Harder-Narasimhan stratification of parabolic induction}     \label{subsection-HN-Bun-G-P-M}

Recall that we have morphisms (\ref{equation-Bun-G-P-M}): $\Bun_{G} \xleftarrow{i^{Bun}} \Bun_{P} \xrightarrow{\pi^{Bun}} \Bun_{M}.$

\begin{defi}   \label{def-Bun-P-leq-ad-mu}
Let $\mu \in \wh{\Lambda}_{G^{\mr{ad}}}^{+,\Q}$. We define $\Bun_{P}^{\leq^{G^{\mr{ad}}} \mu}$ to be the inverse image of $\Bun_G^{\leq^{G^{\mr{ad}}} \mu}$ in $\Bun_P$.
\end{defi}

\begin{lem}     \label{lem-Bun-P-leq-mu-dans-Bun-M-leq-mu}
The image of $\Bun_{P}^{\leq^{G^{\mr{ad}}} \mu}$ in $\Bun_{M}$ is included in $\Bun_{M}^{\leq^{G^{\mr{ad}}} \mu}$.
\end{lem}
\dem
Let $\mc P \in \Bun_{P}^{\leq^{G^{\mr{ad}}} \mu}$ and let $\mc M$ be its image in $\Bun_M$. 
We will check that $\mc M \in \Bun_{M}^{\leq^{G^{\mr{ad}}}  \mu}$. 
For any parabolic subgroup $P'$ of $M$, let $M'$ be its Levi quotient. Let $\mc P'$ be a $P'$-structure of $\mc M$ and $\mc M':= \mc P' \overset{P'} \times M'$. By Definition \ref{def-troncature-Bun-M-leq-ad-mu}, we need to prove that $\Upsilon_G \circ \phi_{P'} \circ \on{deg}_{P'} (\mc P') \leq^{G^{\mr{ad}}} \mu$.

Let $P'':=P \underset{M}\times P'$. It is a parabolic subgroup of $G$ with Levi quotient $M'$. We have
$$
\xymatrixrowsep{0.5pc}
\xymatrixcolsep{0.5pc}
\xymatrix{
 & & P''  \ar[rd]  \ar[ld]  \\
 & P   \ar[rd]  \ar[ld]    & & P'  \ar[rd]  \ar[ld]    \\
G & & M & & M'
}$$
By \cite{DG16} Lemma 2.5.8, we can define a $P''$-bundle $\mc P'':= \mc P \underset{\mc M} \times \mc P'$.
We have $\deg_{P'} \mc P' = \deg_{M'} \mc M' = \deg_{P''} \mc P''$. Taking into account that $\wh \Lambda_G^{\Q}=\wh \Lambda_M^{\Q}$, we deduce that $\Upsilon_G \circ \phi_{P'} \circ \on{deg}_{P'} (\mc P') = \Upsilon_G \circ \phi_{P''} \circ \on{deg}_{P''} (\mc P'') \leq^{G^{\mr{ad}}}  \mu$, where the last inequality follows from the definition of $\Bun_{P}^{\leq^{G^{\mr{ad}}} \mu}$.
\cqfd

\sssec{}   \label{subsection-Bun-P-leq-ad-mu-dans-Bun-M-leq-ad-mu}
By Lemma \ref{lem-Bun-P-leq-mu-dans-Bun-M-leq-mu}, morphisms (\ref{equation-Bun-G-P-M}) induce morphisms:
\begin{equation}   \label{diagram-Bun-G-Bun-P-Bun-M-leq-ad-mu}
\Bun_G^{\leq^{G^{\mr{ad}}} \mu} \leftarrow \Bun_P^{\leq^{G^{\mr{ad}}} \mu} \rightarrow \Bun_M^{\leq^{G^{\mr{ad}}} \mu}.
\end{equation}

The group $\Xi$ acts on all these stacks. All the morphisms are $\Xi$-equivariant. Thus morphisms (\ref{diagram-Bun-G-Bun-P-Bun-M-leq-ad-mu}) induce morphisms:
\begin{equation}   \label{diagram-Bun-G-Bun-P-Bun-M-leq-ad-mu-Xi}
\Bun_{G}^{\leq^{G^{\mr{ad}}} \mu} / \Xi \leftarrow \Bun_{P}^{\leq^{G^{\mr{ad}}} \mu} / \Xi \rightarrow \Bun_{M}^{\leq^{G^{\mr{ad}}} \mu} / \Xi .
\end{equation}

\sssec{}    \label{subsection-Bun-P-leq-ad-mu-Bun-M-leq-ad-mu-nu}
For any $\nu \in \wh \Lambda_{Z_M / Z_G}^{\Q}$, we define $\Bun_{P}^{\nu}$ to be the inverse image of $\Bun_M^{\nu}$ in $\Bun_P$. We define $\Bun_{P}^{\leq^{G^{\mr{ad}}} \mu, \, \nu}:= \Bun_{P}^{\leq^{G^{\mr{ad}}} \mu} \cap \Bun_P^{\nu}$. Morphisms (\ref{diagram-Bun-G-Bun-P-Bun-M-leq-ad-mu-Xi}) induce morphisms:
\begin{equation}   \label{diagram-Bun-G-Bun-P-Bun-M-leq-ad-mu-Xi-nu}
\Bun_{G}^{\leq^{G^{\mr{ad}}} \mu} / \Xi \leftarrow \Bun_{P}^{\leq^{G^{\mr{ad}}} \mu, \, \nu} / \Xi \rightarrow \Bun_{M}^{\leq^{G^{\mr{ad}}} \mu, \, \nu} / \Xi .
\end{equation}

\subsection{Harder-Narasimhan stratification of stack of shtukas}    \label{subsection-HN-Cht-G-P-M}

\begin{nota}    \label{notation-leq-means-leq-G-ad}
In the remaining part of the paper, we will only use the truncations indexed by "$\leq^{G^{\mr{ad}}}$" (rather than "$\leq^G$"). 
To simplify the notation, from now on, "$\leq$" means "$\leq^{G^{\mr{ad}}}$".
\end{nota}

\begin{defi}   \label{def-Cht-G-P-M-mu}
Let $\mu \in \wh{\Lambda}_{G^{\mr{ad}}}^{+,\Q}$ (resp. $\lambda \in \wh{\Lambda}_{G}^{+,\Q}$). 
We define $\Cht_{G, N, I}^{\leq \mu}$ 
(resp. $\Cht_{G, N, I}^{(\lambda)}$) 
to be the inverse image of $\Bun_{G}^{\leq \, \mu}$ 
(resp. $\Bun_{G}^{(\lambda)}$) 
by the morphism $$\Cht_{G, N, I} \rightarrow \Bun_{G}, \quad \left( (x_i)_{i \in I}, (\mc G, \psi) \xrightarrow{\phi}  (\ta \mc G, \ta \psi) \right) \mapsto \mc G.$$ 
Similarly, we define $\Cht_{M, N, I}^{\leq \mu}$ (resp. $\Cht_{M, N, I}^{\leq \mu, \, \nu}$, $\Cht_{M, N, I}^{(\lambda)}$) using the morphism $\Cht_{M, N, I} \rightarrow \Bun_{M}$
and $\Cht_{P, N, I}^{\leq \mu}$ (resp. $\Cht_{P, N, I}^{\leq \mu, \, \nu}$) using the morphism $\Cht_{P, N, I} \rightarrow \Bun_{P}$.
%
%
%
\end{defi}

\sssec{}
The following diagram is commutative:
\begin{equation}      \label{equation-Cht-G-P-M-Bun-G-P-M}
\xymatrixrowsep{1pc}
\xymatrixcolsep{2pc}
\xymatrix{
\Cht_{G, N, I} \ar[d]
&\Cht_{P, N, I}  \ar[l]_{i}   \ar[d] \ar[r]^{\pi} 
&  \Cht_{M, N, I}  \ar[d] \\
\Bun_G     
& \Bun_P   \ar[l]_{i^{Bun} }   \ar[r]^{\pi^{Bun}} 
&  \Bun_M
}
\end{equation}
where the first line is defined in (\ref{diagram-Cht-G-P-M-general}). We deduce that $\Cht_{P, N, I}^{\leq \, \mu}$ is the inverse image of $\Cht_{G, N, I}^{\leq \, \mu}$ in $\Cht_{P, N, I}$.

\begin{lem}    \label{lem-Cht-P-leq-mu-to-Cht-M-leq-mu}
The image of $\Cht_{P, N, I}^{\leq \, \mu}$ in $\Cht_{M, N, I}$ is included in $\Cht_{M, N, I}^{\leq \, \mu}$.
\end{lem}
\dem
This follows from Lemma \ref{lem-Bun-P-leq-mu-dans-Bun-M-leq-mu} and the commutativity of (\ref{equation-Cht-G-P-M-Bun-G-P-M}).
\cqfd

\sssec{}
Just as in \ref{subsection-Bun-P-leq-ad-mu-dans-Bun-M-leq-ad-mu} and \ref{subsection-Bun-P-leq-ad-mu-Bun-M-leq-ad-mu-nu}, morphisms (\ref{diagram-Cht-G-P-M-general}) induce morphisms:
\begin{equation}     \label{equation-Cht-G-P-M-leq-mu-Xi}
\Cht_{G, N, I}^{\leq \mu} / \Xi \leftarrow \Cht_{P, N, I}^{\leq \mu} / \Xi \rightarrow \Cht_{M, N, I}^{\leq \mu} / \Xi
\end{equation}
\begin{equation}
\Cht_{G, N, I}^{\leq \mu} / \Xi \leftarrow \Cht_{P, N, I}^{\leq \mu, \, \nu} / \Xi \rightarrow \Cht_{M, N, I}^{\leq \mu, \, \nu} / \Xi .
\end{equation}
We deduce from (\ref{formule-decom-Bun-M-leq-ad-mu-cone-Xi}) a decomposition:
\begin{equation}    \label{equation-Cht-M-leq-mu-decomp-en-nu}
\Cht_{M, N, I}^{\leq \mu} / \Xi = \underset{ \nu \in \wh{\Lambda}_{Z_M / Z_G}^{\mu}  } \bigsqcup \Cht_{M, N, I}^{\leq \mu,  \; \nu}  / \Xi
\end{equation}

\quad

\section{Cohomology of stacks of shtukas}

In Sections \ref{subsection-geo-Satake}-\ref{subsection-coho-Cht-G} we recall the definition of the cohomology of stacks of $G$-shtukas with values in perverse sheaves coming from $[  G_{I, \infty} \backslash  \on{Gr}_{G, I}] $ via $\epsilon_{G, N, I, \infty}$, i.e. coming from $G_{I, \infty} $-equivariant perverse sheaves over $\on{Gr}_{G, I}$. 
These sections are based on \cite{vincent} Sections 1, 2 and 4.

In Section \ref{subsection-coho-Cht-M} we define the cohomology of stacks of $M$-shtukas.


\begin{nota}
Our results are of geometric nature, i.e. we will not consider the action of $Gal(\Fqbar / \Fq)$. 
From now on, 
we pass to the base change over $\Fqbar$. We keep the same notations $X$, $\Bun_{G, N}$, $\Cht_{G, N, I}$, $\Gr_{G, I}$, etc... but now everything is over $\Fqbar$ and the fiber products are taken over $\Fqbar$.
\end{nota}



\subsection{Reminder of a generalization of the geometric Satake equivalence}    \label{subsection-geo-Satake}


\sssec{}
The geometric Satake equivalence for the affine grassmannian is established in \cite{mv} over the ground field $\C$. By \cite{mv} Section 14, \cite{ga-de-jong} Section 1.6 and \cite{zhu}, the constructions in \cite{mv} carries over to the case of an arbitrary algebraically closed ground field of characteristic prime to $\ell$.

\sssec{}   \label{subsection-def-G-dual}
Let $\wh G$ be the Langlands dual group of $G$ over $\Ql$ defined by the geometric Satake equivalence for the affine grassmannian, as in \cite{mv} Theorem 7.3 and \cite{ga-de-jong} Theorem 2.2.

\sssec{}
(\cite{mv} Section 2, \cite{ga-iwahori} 1.1.1 and Section 6)
The Beilinson-Drinfeld affine grassmannian $\Gr_{G, I}$ is an ind-scheme. Every finite-dimensional closed subscheme of $\Gr_{G, I}$ is contained in some finite-dimensional closed subscheme of $\Gr_{G, I}$ stable under the action of $G_{I, \infty}$. 

We denote by $\Perv_{G_{I, \infty}}(\Gr_{G, I}, \Ql)$ the category of $G_{I, \infty}$-equivariant perverse sheaves with $\Ql$-coefficients on $\Gr_{G, I}$ (for the perverse normalization relative to $X^I$).

\sssec{} 
As in \cite{ga-de-jong} 2.5, we denote by $P^{\wh G, I}$ the category of perverse sheaves with $\Ql$-coefficients on $X^I$ (for the perverse normalization relative to $X^I$) endowed with an extra structure given in {\it loc.cit}. 

\begin{thm}  (\cite{ga-de-jong} Theorem 2.6)  \label{thm-equivalence-Perv-Gr-G-to-P-G-I}
There is a canonical equivalence of categories $\Perv_{G_{I, \infty}}(\Gr_{G, I}, \Ql) \isom P^{\wh G, I}$, compatible with the tensor structures defined in {\it loc.cit}.
\cqfd
\end{thm} 

\sssec{}
We denote by $\on{Rep}_{\Ql}(\wh G^I)$ the category of finite dimensional $\Ql$-linear representations of $\wh G^I$. We have a fully faithful functor
$\on{Rep}_{\Ql}(\wh G^I) \rightarrow P^{\wh G, I}:  W \mapsto W \otimes {\Ql}_{ X^I }.$ The composition of this functor and the inverse functor $ P^{\wh G, I} \isom \Perv_{G_{I, \infty}}(\Gr_{G, I}, \Ql) $ in Theorem \ref{thm-equivalence-Perv-Gr-G-to-P-G-I} gives:

\begin{cor}  \label{cor-Satake-functor-I}
We have a canonical natural fully faithful $\Ql$-linear fiber functor:
$$Sat_{G, I}: \on{Rep}_{\Ql}(\wh G^I) \rightarrow \Perv_{G_{I, \infty}}(\Gr_{G, I}, \Ql).$$ 
\end{cor}


\begin{defi}  \label{def-Gr-G-I-W}
For any $W \in \on{Rep}_{\Ql}(\wh G^I)$, we define $\mc S_{G, I, W} := Sat_{G, I}(W)$.
We define $\Gr_{G, I, W}$ to be the support of $\mc S_{G, I, W}$. 
\end{defi}

\sssec{}  \label{subsection-S-G-I-W-1-plus-W-2}
When $W = W_1 \oplus W_2$, by the functoriality of $Sat_{G, I}$, we have $\mc S_{G, I, W} = \mc S_{G, I, W_1} \oplus \mc S_{G, I, W_2}$. Then $\Gr_{G, I, W} = \Gr_{G, I, W_1} \cup \Gr_{G, I, W_2}$.

\sssec{}   \label{subsection-grass-aff-globalization}
By \cite{vincent} Théorème 1.17, the above definition of $\Gr_{G, I, W}$ is equivalent to $loc.cit.$ Définition 1.12 and the definition after (1.14) (which describes $\Gr_{G, I, W}$ as a generalization of the Zariski closure of the Schubert cell in affine grassmannian). 
It is well-known that $\Gr_{G, I, W}$ is a closed subscheme of $\Gr_{G, I}$ 
and it is projective (see \cite{mv} Sections. 2-3, \cite{zhu} Proposition 2.1.5). $\Gr_{G, I}$ is an inductive limit of $\Gr_{G, I, W}$.

\begin{rem}    \label{rem-S-G-I-W-isom-IC}
By \cite{vincent} Théorème 1.17, when $W$ is irreducible, the perverse sheaf $\mc S_{G, I, W}$ is (not canonically) isomorphic to the intersection complex (with coefficient in $\Ql$ and the perverse normalization relative to $X^I$) of $\Gr_{G, I, W}$. 
\end{rem}


\subsection{Satake perverse sheaves on quotient stacks}

The stacks $[G_{I, \infty} \backslash \Gr_{G, I}]$ or $[G_{I, \infty} \backslash \Gr_{G, I, W}]$ are not algebraic because the group scheme $G_{I, \infty}$ is of infinite dimension. For technical reasons, we will need algebraic stacks.

%
%
%
%

\begin{prop}   \label{prop-d-assez-grand}  (\cite{ga-iwahori} 1.1.1)
For $d \in \Z_{\geq 0}$ large enough depending on $W$, the action of $\Ker(G_{I, \infty} \rightarrow G_{I, d})$ on $\Gr_{G, I, W}$ is trivial. Thus the action of $G_{I, \infty}$ on $\Gr_{G, I, W}$ factors through $G_{I, d}$.
\cqfd
\end{prop}


\sssec{}
For $d$ as in Proposition \ref{prop-d-assez-grand}, we define the quotient stack $[G_{I, d} \backslash \Gr_{G, I, W}]$. Since the group scheme $G_{I, d}$ is of finite dimension, the stack $[G_{I, d} \backslash \Gr_{G, I, W}]$ is algebraic.




\sssec{}    \label{subsection-def-S-G-I-W-d}
Let $\mc S_{G, I, W}$ be the $G_{I, \infty}$-equivariant perverse sheaf on $\Gr_{G, I, W}$ defined in Definition \ref{def-Gr-G-I-W}.
By Proposition \ref{prop-d-assez-grand}, the action of $G_{I, \infty}$ on $\mc S_{G, I, W}$ factors through $G_{I, d}$. Since the kernel of $G_{I, \infty} \twoheadrightarrow G_{I, d}$ is connected, by \cite{bbd} Proposition 4.2.5, we deduce that $\mc S_{G, I, W}$ is also $G_{I, d}$-equivariant. 


Let $\xi_{G, I, d}: \Gr_{G, I, W} \rightarrow [G_{I, d} \backslash \Gr_{G, I, W}]$ be the canonical morphism. It is smooth of dimension $\dim G_{I, d}$. By $loc.cit.$ Corollaire 4.2.6.2 and the discussion after it, there exists a perverse sheaf (up to shift $[\dim G_{I, d}]$) (for the perverse normalization relative to $X^I$) $\mc S_{G, I, W}^d$ on $[G_{I, d} \backslash \Gr_{G, I, W}]$ such that $\mc S_{G, I, W} = \xi_{G, I, d}^* \mc S_{G, I, W}^d$.

\sssec{}   \label{subsection-proj-sys-quotient}
Let $d \leq d' $ be two integers large enough as in Proposition \ref{prop-d-assez-grand}. Then the morphisms $G_{I, \infty} \twoheadrightarrow G_{I, d'} \twoheadrightarrow G_{I, d}$
induce a commutative diagram:
\begin{equation}  
\xymatrix{
 &   \Gr_{G, I, W}   \ar[ld]  \ar[d]_{\xi_{G, I, d'}}   \ar[rd]^{\xi_{G, I, d}} & \\
[G_{I, \infty} \backslash  \Gr_{G, I, W}]  \ar[r]  & [G_{I, d'} \backslash  \Gr_{G, I, W}]  \ar[r]^{\on{pr}_{d'}^{d}}  & [G_{I, d} \backslash  \Gr_{G, I, W}] 
}
\end{equation}


We have $(\xi_{G, I, d'})^* \mc S_{G, I, W}^{d'} = \mc S_{G, I, W} =  (\xi_{G, I, d})^* \mc S_{G, I, W}^{d} = (\xi_{G, I, d'})^* (\on{pr}_{d'}^{d})^* \mc S_{G, I, W}^d$. By \cite{bbd} Proposition 4.2.5, the functor $(\xi_{G, I, d'})^*$ (up to shift) is fully faithful. We deduce that $\mc S_{G, I, W}^{d'} =  (\on{pr}_{d'}^{d})^* \mc S_{G, I, W}^d$.


\sssec{}
By Proposition \ref{prop-d-assez-grand}, the action of $G_{I, \infty}^{\mr{ad}}$ on $\Gr_{G, I, W}$ factors through $G_{I, d}^{\mr{ad}}$. We define the quotient stack $[G_{I, d}^{\mr{ad}} \backslash \Gr_{G, I, W}]$.

As in the discussion after Définition 2.14 in \cite{vincent}, since $(Z_G)_{I, \infty}$ acts trivially on $\Gr_{G, I, W}$, the $G_{I, \infty}$-equivariant perverse sheaf $\mc S_{G, I, W}$ on $\Gr_{G, I, W}$ is also $G^{\mr{ad}}_{I, \infty}$-equivariant and $G^{\mr{ad}}_{I, d}$-equivariant. Indeed, by \ref{subsection-S-G-I-W-1-plus-W-2} it is enough to prove this for $W$ irreducible. By Remark \ref{rem-S-G-I-W-isom-IC},
in this case $\mc S_{G, I, W}$ is isomorphic to the intersection complex 
of $\Gr_{G, I, W}$, hence is $G^{\mr{ad}}_{I, \infty}$-equivariant.

Just as in \ref{subsection-def-S-G-I-W-d}, let $\xi_{G, I, d}^{\mr{ad}}: \Gr_{G, I, W} \rightarrow [G_{I, d}^{\mr{ad}} \backslash \Gr_{G, I, W}]$ be the canonical morphism. There exists a perverse sheaf (up to shift $[\dim G_{I, d}^{\mr{ad}} ]$) (for the perverse normalization relative to $X^I$) $\mc S_{G, I, W}^{\mr{ad}, \, d}$ on $[G_{I, d}^{\mr{ad}}  \backslash \Gr_{G, I, W}]$ such that $\mc S_{G, I, W} = (\xi_{G, I, d}^{\mr{ad}})^* \mc S_{G, I, W}^{\mr{ad}, \, d}$.

\subsection{Representability of stacks of shtukas}   

\begin{defi}
We define $\Cht_{G, N, I, W}$ to be the inverse image of $[G_{I, \infty} \backslash \Gr_{G, I, W}]$ in $\Cht_{G, N, I}$ by $\epsilon_{G, N, I, \infty}$.
\end{defi}

\sssec{}
$\Cht_{G, N, I}$ is an inductive limit of closed subtacks $\Cht_{G, N, I, W}$. 

\sssec{}
Let $\mu \in \wh{\Lambda}_{G^{\mr{ad}}}^{+, \Q}$. We define $\Cht_{G, N, I, W}^{\leq \mu} := \Cht_{G, N, I, W} \cap \Cht_{G, N, I}^{\leq \mu} $, where $\Cht_{G, N, I}^{\leq \mu} $ is defined in Definition \ref{def-Cht-G-P-M-mu}.
We define the quotient $\Cht_{G, N, I, W} / \Xi$ and $\Cht_{G, N, I, W}^{\leq \mu} / \Xi$.

\begin{prop}  \label{prop-Cht-G-representable}  (\cite{var}   Proposition 2.16  )
$\Cht_{G, N, I, W}$ is a Deligne-Mumford stack locally of finite type. 
$\Cht_{G, N, I, W}^{\leq \mu} / \Xi$ is a Deligne-Mumford stack of finite type. 
\cqfd
\end{prop}

\sssec{}
$\Cht_{G, N, I, W} / \Xi = \varinjlim _{\mu \in \wh{\Lambda}_{G^{\mr{ad}}}^{+, \Q} } \Cht_{G, N, I, W}^{\leq \mu} / \Xi$ is locally of finite type.

\subsection{Satake perverse sheaf on stacks of shtukas}   \label{section-perv-sheaf-on-cht}


\sssec{}    \label{subsection-epsilon-G-N-I-d}
For any $d \in \Z_{\geq 0}$ large enough as in Proposition \ref{prop-d-assez-grand}, 
we define $\epsilon_{G, N, I, d}$ to be the composition of morphisms
\begin{equation}   \label{equation-epsilon-G-d}
\epsilon_{G, N, I, d}: \Cht_{G, N, I, W} \xrightarrow{ \epsilon_{G, N, I, \infty} } [G_{I, \infty} \backslash \Gr_{G, I, W}] \rightarrow [G_{I, d} \backslash \Gr_{G, I, W}].
\end{equation}
This is morphism (2.3) in \cite{vincent}. 

Just as in \ref{subsection-Cht-Xi-to-G-ad}, we define a morphism
\begin{equation}   \label{equation-epsilon-G-d-Xi}
\epsilon_{G, N, I, d}^{\Xi}: \Cht_{G, N, I, W} / \Xi \rightarrow [G_{I, d}^{\mr{ad}} \backslash \Gr_{G, I, W}] .
\end{equation}
This is morphism (2.10) in {\it loc.cit.}.

\sssec{}
We denote by $\dim_{X^I} G_{I, d}$ the relative dimension of $G_{I, d}$ over $X^I$ and by $|I|$ the cardinal of $I$. We have $\dim_{X^I} G_{I, d} = d \cdot |I| \cdot \on{dim} G$.

\begin{prop} \label{prop-epsilon-d-smooth} (\cite{vincent} Proposition 2.8)
The morphisms $\epsilon_{G, N, I, d}$ (resp. $\epsilon_{G, N, I, d}^{\Xi}$) is smooth of dimension $\dim_{X^I} G_{I, d}$ (resp. $\dim_{X^I} G_{I, d}^{ad}$). 
\cqfd
\end{prop}

\sssec{}
For all $d \in \Z_{\geq 0}$ large enough as in Proposition \ref{prop-d-assez-grand}, we have morphisms over $(X \sm N)^I$:
$$
\xymatrixrowsep{1pc}
\xymatrixcolsep{1pc}
\xymatrix{
\Cht_{G, N, I, W}   \ar[rd]^{\epsilon_{G, N, I, d} }  
& & \Gr_{G, I, W}   \ar[ld]_{\xi_{G, I, d} }      \\
& [G_{I, d} \backslash  \Gr_{G, I, W}]  
}
$$

We deduce from Proposition \ref{prop-epsilon-d-smooth} that $\dim \Cht_{G, N, I, W} = \dim \Gr_{G, I, W}$. We refer to \cite{vincent} Proposition 2.11 for the fact that $\Cht_{G, N, I, W}$ is locally isomorphic to $\Gr_{G, I, W}$ for the étale topology. We will not use this result in this paper.

\begin{defi}  
Let $d \in \Z_{\geq 0}$ large enough as in Proposition \ref{prop-d-assez-grand}.
We define $\mc{F}_{G, N, I, W} := (\epsilon_{G, N, I, d})^* \mc S_{G, I, W}^d $. 
\end{defi}

\begin{rem}   \label{rem-F-G-I-W-independent-of-d}
As in \ref{subsection-proj-sys-quotient}, let $d, d' \in \Z_{\geq 0}$ both large enough with $d \leq d'$. Then we have $\epsilon_{G, N, I, d} =  \on{pr}_{d'}^{d} \circ \epsilon_{G, N, I, d'}$. Thus $(\epsilon_{G, N, I, d})^* \mc S_{G, I, W}^d = (\epsilon_{G, N, I, d'})^* (\on{pr}_{d'}^{d} )^* \mc S_{G, I, W}^{d} = (\epsilon_{G, N, I, d'})^* \mc S_{G, I, W}^{d'}$. Hence $\mc{F}_{G, N, I, W}$ is independent of $d$.  
\end{rem}

\begin{defi}    
We define $\mc{F}_{G, N, I, W}^\Xi := (\epsilon_{G, N, I, d}^{\Xi})^* \mc S_{G, I, W}^{\mr{ad}, \, d}$. 
\end{defi}

Just as in Remark \ref{rem-F-G-I-W-independent-of-d}, $\mc{F}_{G, N, I, W}^\Xi$ is independent of $d$.

\begin{lem}
$\mc{F}_{G, N, I, W}$ (resp. $\mc{F}_{G, N, I, W}^\Xi $) is a perverse sheaf (for the perverse normalization relative to $(X \sm N)^I$) on $\Cht_{G, N, I}$ (resp. $\Cht_{G, N, I} / \Xi$) supported on $\Cht_{G, N, I, W} $ (resp. $\Cht_{G, N, I, W} / \Xi$) (in the context of \ref{def-D-c-b-pour-chat-alg}). When $W$ is irreducible, $\mc{F}_{G, N, I, W}$ (resp. $\mc{F}_{G, N, I, W}^\Xi $) is (not canonically) isomorphic to the intersection complex (with coefficient in $\Ql$ and the perverse normalization relative to $(X \sm N)^I$) of $\Cht_{G, N, I, W}$ (resp. $\Cht_{G, N, I, W} / \Xi$). 
\end{lem}
\dem
The lemma follows from Corollary \ref{cor-Satake-functor-I}, Remark \ref{rem-S-G-I-W-isom-IC} and Proposition \ref{prop-epsilon-d-smooth}.
\cqfd


\subsection{Cohomology of stacks of $G$-shtukas}     \label{subsection-coho-Cht-G}




Recall that we have the morphism of paws $\mf{p}_G: \Cht_{G, I, N} / \Xi \rightarrow (X \sm N)^I.$ 
\begin{defi} (\cite{vincent} Definition 4.1, 4.7)   \label{def-mc-H-G-leq-mu}
For any $\mu \in \wh{\Lambda}_{G^{\mr{ad}}}^{+, \Q}$, we define 
$$
\mc H_{G, N, I, W}^{\leq \mu} := R(\mf{p}_G)_! ( \restr{\mc{F}_{G, N, I, W}^\Xi } {\Cht_{G, N, I, W}^{\leq \mu} / \Xi } ) \in D_c^b( (X \sm N)^I , \Ql ).
$$
For any $j \in \Z$, we define degree $j$ cohomology sheaf (for the ordinary $t$-structure):
$$
\mc H _{G, N, I, W}^{j, \; \leq\mu}:=R^j(\mf{p}_G)_! ( \restr{\mc{F}_{G, N, I, W}^\Xi} {\Cht_{G, N, I, W}^{\leq \mu} / \Xi } ).
$$
This is a $\Ql$-constructible sheaf on $ (X \sm N)^I$. 
\end{defi}

$\mc H_{G, N, I, W}^{\leq \mu} $ and $\mc H _{G, N, I, W}^{j, \; \leq\mu}$ depend on $\Xi$. We do not write $\Xi$ in the index to simplify the notations.

\sssec{}   \label{subsection-H-G-leq-mu-1-to-leq-mu-2}
Let $\mu_1, \mu_2 \in \wh{\Lambda}_{G^{\mr{ad}}}^{+, \Q}$ and $\mu_1 \leq \mu_2$. We have an open immersion:
\begin{equation}   \label{equation-Cht-G-leq-mu-1-to-leq-mu-2}
\Cht_{G, N, I, W}^{\leq \mu_1} / \Xi  \hookrightarrow \Cht_{G, N, I, W}^{\leq \mu_2} / \Xi.
\end{equation}
For any $j$, morphism (\ref{equation-Cht-G-leq-mu-1-to-leq-mu-2}) induces a morphism of sheaves:
$$\mc H _{G, N, I, W}^{j, \, \leq\mu_1} \rightarrow \mc H _{G, N, I, W}^{j, \, \leq \mu_2}.$$

\begin{defi} We define
$$\mc H _{G, N, I, W}^j: = \varinjlim _{\mu}  \mc H _{G, N, I, W}^{j, \; \leq\mu} $$
as an inductive limit in the category of constructible sheaves on $ (X \sm N)^I$. 
\end{defi}

\sssec{} \label{subsection-eta-I-bar}
Let $\ov{\eta^I}$ be a geometric point over the generic point $\eta^I$ of $X^I$.

\begin{defi}    \label{def-H-G-j}
We define
\begin{equation}   \label{equation-H-G-est-mc-H-G-restriction}
H _{G, N, I, W}^{j, \, \leq \mu}:=\restr{ \mc H _{G, N, I, W}^{j, \, \leq \mu} }{ \ov{\eta^I}  }  \; ,  \quad     H _{G, N, I, W}^j:=\restr{ \mc H _{G, N, I, W}^j }{ \ov{\eta^I}  }.
\end{equation}
$H _{G, N, I, W}^{j, \, \leq \mu}$ is a $\Ql$-vector space of finite dimension. $H _{G, N, I, W}^j = \varinjlim _{\mu}  H _{G, N, I, W}^{j, \; \leq\mu} $.
\end{defi}

\subsection{Cohomology of stacks of $M$-shtukas}      \label{subsection-coho-Cht-M}

Let $P$ be a proper parabolic subgroup of $G$ and let $M$ be its Levi quotient. 

\sssec{}
Let $\wh M$ be the Langlands dual group of $M$ over $\Ql$ defined by the geometric Satake equivalence.
The compatibility between the geometric Satake equivalence and the constant term functor along $P$ (that we will recall in Theorem \ref{thm-geo-satake-CT-I-paws} below) induces a canonical inclusion $\wh M \hookrightarrow \wh G$ (compatible with pinning). 

\sssec{}    \label{subsection-epsilon-M-N-I-d}
We view $W \in \on{Rep}_{\Ql}(\wh G^I)$ as a representation of $\wh M^I$ via $\wh M^I \hookrightarrow \wh G^I$. 
As in Sections \ref{subsection-geo-Satake}-\ref{section-perv-sheaf-on-cht}, we define $\Gr_{M, I, W}$ and $\Cht_{M, N, I, W}$. 
For $d \in \Z_{\geq 0}$ large enough such that the action of $M_{I, \infty}$ on $\Gr_{M, I, W}$ factors through $M_{I, d}$, we define
$$\epsilon_{M, N, I, d}: \Cht_{M, N, I, W} \rightarrow [M_{I, d} \backslash \Gr_{M, I, W}],$$
$$\epsilon_{M, N, I, d}^{\Xi}: \Cht_{M, N, I, W} / \Xi \rightarrow [\ov{M}_{I, d} \backslash \Gr_{M, I, W}].$$
We define perverse sheaf $\mc S_{M, I, W}$ on $\Gr_{M, I, W}$, perverse sheaves (up to shift) $\mc S_{M, I, W}^{d}$ on $[M_{I, d} \backslash \Gr_{M, I, W}]$ and $\mc S_{M, I, W}^{\mr{ad}, d}$ on $[\ov{M}_{I, d} \backslash \Gr_{M, I, W}]$. We define $\mc{F}_{M, N, I, W}: = \epsilon_{M, N, I, d}^* \mc S_{M, I, W}^{d}$ and $\mc{F}_{M, N, I, W}^\Xi: = (\epsilon_{M, N, I, d}^{\Xi})^* \mc S_{M, I, W}^{\mr{ad}, d}$.

%
%

\sssec{}
Applying \cite{var} Proposition 2.16 to $M$, we deduce that $\Cht_{M, I, N, W}$ is a Deligne-Mumford stack locally of finite type and that for $\lambda \in  \wh{\Lambda}_{M}^{+, \Q}$, the Deligne-Mumford stack $\Cht_{M, I, N, W}^{(\lambda)}$ (defined in Definition \ref{def-Cht-G-P-M-mu}) is of finite type.

Let $\mu \in \wh{\Lambda}_{G^{\mr{ad}}}^{+, \Q}$. We define $\Cht_{M, N, I, W}^{\leq \mu} :=  \Cht_{M, N, I, W} \cap \Cht_{M, N, I}^{\leq \mu}$, where $\Cht_{M, N, I}^{\leq \mu}$ is defined in Definition \ref{def-Cht-G-P-M-mu}. We define the quotient $\Cht_{M, N, I, W} / \Xi$ and $\Cht_{M, N, I, W}^{\leq \mu} / \Xi$. As in \ref{subsection-Xi-not-lattice-in-Z-M-for-Bun-M}, $\Xi$ is a lattice in $Z_G(F) \backslash Z_G(\mb A)$ but only a discrete subgroup in $Z_M(F) \backslash Z_M(\mb A)$. 
The decomposition (\ref{equation-Cht-M-leq-mu-decomp-en-nu}) induces a decomposition
\begin{equation}
\Cht_{M, N, I, W}^{\leq \mu} / \Xi = \underset{ \nu \in \wh{\Lambda}_{Z_M / Z_G}^{\mu}  } \bigsqcup \Cht_{M, N, I, W}^{\leq \mu,  \; \nu}  / \Xi .
\end{equation}
where each $\Cht_{M, I, N, W}^{\leq \mu, \; \nu} / \Xi$ is of finite type (just as in Lemma \ref{lem-Bun-M-leq-mu-nu-Xi-finite-type}).




\quad

Recall that we have the morphism of paws $\mf{p}_M: \Cht_{M, I, N} / \Xi \rightarrow (X \sm N)^I.$
\begin{defi}   \label{def-H-M-mu-nu}
For any $\mu \in \wh{\Lambda}_{G^{\mr{ad}}}^{+, \Q}$ and $\nu \in \wh{\Lambda}_{Z_M / Z_G}^{\Q}$, we define 
$$\mc H _{M, N, I, W}^{\leq  \mu, \; \nu}  := R(\mf{p}_M)_! ( \restr{ \mc{F}_{M, I, N, W}^\Xi    } {\Cht_{M, N, I, W}^{\leq  \mu, \; \nu} / \Xi} )  \in  D_c^b( (X \sm N)^I , \Ql ) ;$$
For any $j \in \Z$, we define degree $j$ cohomology sheaf
$$\mc H _{M, N, I, W}^{j, \; \leq \mu, \; \nu}  := R^j(\mf{p}_M)_! ( \restr{ \mc{F}_{M, I, N, W}^\Xi    } {\Cht_{M, N, I, W}^{\leq  \mu, \; \nu} / \Xi} ) .$$ 
\end{defi}

\sssec{}    \label{subsection-nu-not-in-cone-H-M-mu-nu-nulle}
If $\nu \notin \wh{\Lambda}_{Z_M / Z_G}^{\mu} $, by Lemma \ref{lem-Bun-M-leq-ad-mu-nu-non-vide}, $\Cht_{M, I, N, W}^{\leq \mu, \; \nu} / \Xi = \emptyset$. In this case $\mc H _{M, N, I, W}^{\leq  \mu, \; \nu}=0$. 


\begin{defi}    \label{def-H-M-j-leq-mu}
Let $\ov{\eta^I}$ be the geometric generic point of $X^I$ fixed in \ref{subsection-eta-I-bar}.
We define 
\begin{equation}    \label{equation-H-M-est-mc-H-M-restriction}
H _{M, N, I, W}^{j, \, \leq  \mu, \; \nu} := \restr{  \mc H _{M, N, I, W}^{j, \, \leq  \mu, \; \nu}  }{  \ov{\eta^I}  }.
\end{equation}
This is a finite dimensional $\Ql$-vector space.
We define 
\begin{equation}
H _{M, N, I, W}^{j, \, \leq  \mu} := \prod_{\nu \in \wh{\Lambda}_{Z_M / Z_G}^{\mu}} H _{M, N, I, W}^{j, \, \leq  \mu, \; \nu} .
\end{equation}
\end{defi}

\sssec{}
Let $\mu_1, \mu_2 \in \wh{\Lambda}_{G^{\mr{ad}}}^{+, \Q}$ and $\mu_1 \leq \mu_2$. We have an open immersion:
\begin{equation}    \label{equation-Cht-M-leq-mu-1-leq-mu-2}
\Cht_{M, N, I, W}^{\leq \mu_1} / \Xi \hookrightarrow \Cht_{M, N, I, W}^{\leq \mu_2} / \Xi.
\end{equation}
For any $j$, morphism (\ref{equation-Cht-M-leq-mu-1-leq-mu-2}) induces a morphism of vector spaces:
$$H _{M, N, I, W}^{j, \, \leq\mu_1} \rightarrow H _{M, N, I, W}^{j, \, \leq \mu_2}.$$

\begin{defi}   \label{def-H-M-j}
We define
$$   
H _{M, N, I, W}^j: = \varinjlim _{\mu}  H _{M, N, I, W}^{j, \; \leq\mu}
$$
as an inductive limit in the category of $\Ql$-vector spaces.
\end{defi}

\begin{defi}    \label{def-H-M-j-nu}
For any $\nu \in \wh{\Lambda}_{Z_M / Z_G}^{\Q}$, we define $H _{M, N, I, W}^{j, \, \nu}: = \varinjlim _{\mu}  H _{M, N, I, W}^{j, \, \leq\mu, \, \nu}$ as an inductive limit in the category of $\Ql$-vector spaces.
\end{defi}

\quad

\section{Constant term morphisms and cuspidal cohomology}

Let $P$ be a parabolic subgroup of $G$ and $M$ its Levi quotient. 
Let $W \in \on{Rep}_{\Ql}(\wh G^I)$. The goal of this section is to construct a constant term morphism from $H_{G, N, I, W}^j $ to $H_{M, N, I, W}^j$ (in fact, to a variant $H_{M, N, I, W}^{' \, j}$ of $H_{M, N, I, W}^j$ defined in Section \ref{subsection-more-on-H-M-prime} below). There are two steps.

Firstly, we will construct a commutative diagram 
\begin{equation}  \label{diagram-TC-Cht-Xi}
\xymatrixrowsep{1pc}
\xymatrixcolsep{1pc}
\xymatrix{
& \Cht_{P, N, I, W} / \Xi \ar[ld]_{i}  \ar[rd]^{\pi}  \ar[dd]_{\mf p_P} & \\
\Cht_{G, N, I, W} / \Xi   \ar[rd]_{\mf p_G}
& & \Cht_{M, N, I, W} / \Xi  \ar[ld]^{\mf p_M} \\
& (X \sm N)^I 
}
\end{equation}
where the morphism $\pi$ is of finite type. 
Therefore the complex $\pi_! i^* \mc{F}_{G, N, I, W}^{\Xi}$ on $\Cht_{M, N, I, W} / \Xi $ is well defined in $D_c^b(\Cht_{M, N, I, W} / \Xi, \Ql)$ (in the context of \ref{def-D-c-b-pour-chat-alg}). 
We will construct a canonical morphism of complexes on $\Cht_{M, N, I, W} / \Xi $:
\begin{equation}   \label{equation-F-G-vers-F-M}
\pi_! i^* \mc{F}_{G, N, I, W}^{\Xi} \rightarrow \mc{F}_{M, N, I, W}^{\Xi}.
\end{equation} 

\quad

Secondly, the cohomological correspondence given by (\ref{diagram-TC-Cht-Xi}) and (\ref{equation-F-G-vers-F-M}) will give a morphism from $H_{G, N, I, W}^j $ to $H_{M, N, I, W}^{j}$.


\subsection{Some geometry of the parabolic induction diagram}

Recall that we have morphisms over $X^I$ in (\ref{diagram-Gr-G-P-M}): $\Gr_{G, I} \xleftarrow{i^0} \Gr_{P, I} \xrightarrow{\pi^0} \Gr_{M, I}$.
\begin{prop}  \label{prop-Gr-G-W-dans-Gr-M-W}
We have $(i^0)^{-1} (\Gr_{G, I, W}) \subset (\pi^0)^{-1} (\Gr_{M, I, W}) ,$
where the inverse images are in the sense of reduced subschemes in $\Gr_{P, I}$.
\end{prop}
\dem
It is enough to prove the inclusion for each fiber over $X^I$. By \ref{subsection-usual-affine-grassmannian}, we reduce the case of the Beilinson-Drinfeld affine grassmannian with paws indexed by $I$ to the case of the usual affine grassmannian $\Gr_G = G_{\mc K} / G_{\mc O}$.

When $P=B$, the statement follows from Theorem 3.2 of \cite{mv}. More concretely, for $\omega$ a dominant coweight of $G$, we denote by $\Gr_{G, \omega}$ the Zariski closure of the Schubert cell defined by $\omega$ in $\Gr_G$. For $\nu$ a coweight of $T$, we denote by $\Gr_{T, \nu}$ the component of $\Gr_T$ (which is discrete) associated to $\nu$. We denote by $C_{\omega}$ the set of coweights of $G$ which are $W$-conjugated to a dominant coweight $\leq \omega$ (where the order is taken in the coweight lattice of $G$).
By Theorem 3.2 of {\it loc.cit.} the subscheme $(i^0)^{-1} (\Gr_{G, \omega}) \cap (\pi^0)^{-1} (\Gr_{T, \nu})$ in $\Gr_B$ is non-empty if and only if $\nu \in C_{\omega}$. Hence 
\begin{equation}   \label{equation-pi-i-Gr-G-W-to-Gr-T-W}
\pi^0 ((i^0)^{-1} (\Gr_{G, \omega})) =\underset{ \nu \in C_{\omega} } \bigsqcup \Gr_{T, \nu}.
\end{equation}

For general $P$ with Levi quotient $M$, we denote by $B'$ the Borel subgroup of $M$. We use the following diagram, where the square is Cartesian:
$$
\xymatrixrowsep{1pc}
\xymatrixcolsep{1pc}
\xymatrix{
 & & \Gr_{B}  \ar[rd]^{\pi_{B'}^B}  \ar[ld]_{i_P^B}  \\
 & \Gr_{P}   \ar[rd]^{\pi_M^P}  \ar[ld]_{i_G^P}    & & \Gr_{B'}   \ar[rd]^{\pi^{B'}_T}  \ar[ld]_{i^{B'}_M}    \\
\Gr_{G} & & \Gr_{M} & & \Gr_T
}$$
Since the square is Cartesian, we have 
\begin{equation}    \label{equation-image-of-Gr-G-omega-is-union-of-Gr-T-nu}
(\pi^{B'}_T) (i^{B'}_M)^{-1} (\pi^P_M) (i^P_G)^{-1} \Gr_{G, \omega} = (\pi^{B'}_T \circ \pi^B_{B'}) (i^P_G \circ i^B_P)^{-1} \Gr_{G, \omega}  \overset{  (  \ref{equation-pi-i-Gr-G-W-to-Gr-T-W} )  }{=} \underset{ \nu \in C_{\omega} } \bigsqcup \Gr_{T, \nu}.
\end{equation}

For any dominant coweight $\lambda$ of $M$, we denote by $\Gr_{M, \lambda}$ the Zariski closure of the Schubert cell defined by $\lambda$ in $\Gr_{M}$. Applying Theorem 3.2 of $loc.cit.$ to $\Gr_M \leftarrow \Gr_{B'} \rightarrow \Gr_T$, we have 
\begin{equation}    \label{equation-pi-i-Gr-M-lambda-to-Gr-T-nu}
(\pi^{B'}_T) (i^{B'}_M)^{-1} \Gr_{M, \lambda} = \underset{ \nu \in C_{\lambda} } \bigsqcup \Gr_{T, \nu}.
\end{equation}

The subscheme $(i^P_G)^{-1} \Gr_{G, \omega}$ in $\Gr_P$ is stable under the action of $P_{\mc O}$. The subscheme $(\pi^P_M) (i^P_G)^{-1} \Gr_{G, \omega}$ in $\Gr_M$ is stable under the action of $M_{\mc O}$, so is a union of strata in $\Gr_M$. 
We deduce from (\ref{equation-image-of-Gr-G-omega-is-union-of-Gr-T-nu}) and (\ref{equation-pi-i-Gr-M-lambda-to-Gr-T-nu}) that $\Gr_{M, \lambda}$ can be in $(\pi^P_M) (i^P_G)^{-1} \Gr_{G, \omega}$ only if $\lambda \in C_\omega$. Thus
$(\pi^P_M) (i^P_G)^{-1} \Gr_{G, \omega}  \subset  \underset{ \lambda \in C_\omega \cap \wh \Lambda_M^+}\bigcup \Gr_{M, \lambda}.$
\cqfd

\sssec{}
We define $\Gr_{P, I, W}:= (i^0)^{-1} (\Gr_{G, I, W})$. As a consequence of Proposition \ref{prop-Gr-G-W-dans-Gr-M-W}, morphisms (\ref{diagram-Gr-G-P-M}) induce morphisms over $X^I$:
\begin{equation}   \label{equation-Gr-G-P-M-I-W}
\Gr_{G, I, W} \xleftarrow{i^0} \Gr_{P, I, W} \xrightarrow{\pi^0} \Gr_{M, I, W}.
\end{equation}


\sssec{}
We deduce from the commutative diagram (\ref{diagram-TC-Cht-Gr}) that
$$i^{-1}(\Cht_{G, N, I, W}) \subset \pi^{-1}(\Cht_{M, N, I, W}),$$ where the inverse images are in the sense of reduced substacks in $\Cht_{P, N, I}$. We define $\Cht_{P, N, I, W}:= i^{-1}(\Cht_{G, N, I, W})$. 
Morphisms in (\ref{diagram-Cht-G-P-M-general}) induce morphisms over $(X \sm N)^I$:
\begin{equation}    \label{equation-Cht-G-P-M-I-W}
\Cht_{G, N, I, W}  \xleftarrow{i} \Cht_{P, N, I, W}  \xrightarrow{\pi} \Cht_{P, N, I, W} . 
\end{equation}

\quad

\sssec{}     \label{subsection-diagram-TC-Cht-Gr-W}
Let $d \in \Z_{\geq 0}$ large enough depending on $W$ as in Proposition \ref{prop-d-assez-grand} applied to $\Gr_{G, I, W}$ and to $\Gr_{M, I, W}$.
To simplify the notations, we write $\epsilon_{G, d}$ for $\epsilon_{G, N, I, d}$ defined in \ref{subsection-epsilon-G-N-I-d} and $\epsilon_{M, d}$ for $ \epsilon_{M, N, I, d}$ defined in \ref{subsection-epsilon-M-N-I-d}. Similarly we define $\epsilon_{P, d}$ to be the composition
$\Cht_{P, N, I, W} \rightarrow  [ P_{I, \infty} \backslash  \Gr_{P, I, W}] \rightarrow  [ P_{I, d} \backslash  \Gr_{P, I, W}].$
 
We deduce from the commutative diagram (\ref{diagram-TC-Cht-Gr}), morphisms (\ref{equation-Gr-G-P-M-I-W}) and (\ref{equation-Cht-G-P-M-I-W}) a commutative diagram of algebraic stacks: 
\begin{equation}   \label{diagram-TC-Cht-Gr-W}
\xymatrix{
\Cht_{G, N, I, W}  \ar[d]^{\epsilon_{G, d}} 
&\Cht_{P, N, I, W}    \ar[l]_i   \ar[d]^{\epsilon_{P, d}}  \ar[r]^{\pi} 
&  \Cht_{M, N, I, W}    \ar[d]^{\epsilon_{M, d}} \\
[G_{I, d} \backslash  \Gr_{G, I, W}]       
& [ P_{I, d} \backslash  \Gr_{P, I, W}]   \ar[l]_{ \ov{i^0_d} }  \ar[r]^{ \ov{\pi^0_d} } 
&  [  M_{I, d} \backslash  \Gr_{M, I, W}]
}
\end{equation} 

\sssec{}
The right square in (\ref{diagram-TC-Cht-Gr-W}) is not Cartesian. We have a commutative diagram, where the square is Cartesian: 
\begin{equation}  \label{diagram-cht-P-cht-M-pi-d}
\xymatrixrowsep{2pc}
\xymatrixcolsep{2pc}
\xymatrix{
\Cht_{P, N, I, W}  \ar[rdd]_{\epsilon_{P, d}}    \ar[rd]|-{\pi_{d}}   \ar[rrd]^{\pi}
& 
&  \\
& \widetilde{\Cht}_{M, N, I, W}     \ar[d]^{ \widetilde{ \epsilon_{M, d } } } \ar[r]_{ \widetilde{  \pi^0_{d}  } }
& \Cht_{M, N, I, W}  \ar[d]^{\epsilon_{M, d}}  \\
&  [ P_{I, d} \backslash  \Gr_{P, I, W}]  \ar[r]_{ \ov{\pi^0_d}  }
&  [  M_{I, d} \backslash  \Gr_{M, I, W}]
}
\end{equation}

\begin{rem}
Note that $\widetilde{\Cht}_{M, N, I, W}$ depends on the choice of $d$. We do not write $d$ in index to shorten the notation.
\end{rem}

\begin{defi}
Let $U$ be the unipotent radical of $P$. We have $P / U = M$. Applying Definition \ref{defi-G-I-infry} to $U$, we define the group scheme $U_{I, d}$ over $X^I$.
\end{defi}

\begin{lem}  \label{lem-pi-d-est-lisse}
The morphism $\pi_d$ is smooth of relative dimension $\dim_{X^I} U_{I, d}$.
\end{lem}

The following proof is suggested to the author by a referee.

\dem
Proposition \ref{prop-epsilon-d-smooth} works also for $P$ and $M$. So the morphism $\epsilon_{P, d}$ is smooth of relative dimension $\dim_{X^I} P_{I, d}$ and the morphism $\epsilon_{M, d}$ (hence $\widetilde{ \epsilon_{M, d} }$) is smooth of relative dimension $\dim_{X^I} M_{I, d}$. Thus to prove that $\pi_d$ is smooth, it is enough to show that it induces a surjective map between relative tangent spaces. 

For any closed point $x_P = ((x_i), \mc P \rightarrow \ta \mc P)$ of $\Cht_{P, N, I, W}$, let $x_M:=\pi_d(x_P)$. We have the canonical morphism 
\begin{equation}   \label{equation-T-P-to-T-M}
T_{\epsilon_{P, d}}(x_P) \rightarrow T_{\widetilde{ \epsilon_{M, d} }}(x_M)
\end{equation}
where $T_{\epsilon_{P, d}}(x_P)$ (resp. $T_{\widetilde{ \epsilon_{M, d} }}(x_M)$) is the tangent space of $\Cht_{P, N, I, W}$ (resp. $\widetilde{\Cht}_{M, N, I, W}$) at $x_P$ (resp. $x_M$) relative to $[ P_{I, d} \backslash  \Gr_{P, I, W}]$. 

Let $y = \epsilon_{P, d}(x_P)$.
By the proof of Proposition 2.8 in \cite{vincent}, we have a Cartesian square
\begin{equation}    \label{equation-Cht-P-cartesian-square}
\xymatrixrowsep{2pc}
\xymatrixcolsep{2pc}
\xymatrix{
\epsilon_{P, d}^{-1}(y)  \ar[r]  \ar[d]
& \Bun_{P, N + d \sum x_i}  \ar[d]^{(b_1^P, \; b_2^P)}  \\
\Bun_{P, N}  \ar[r]^{(\on{Id}, \on{Id}) \quad \quad \quad} 
& \Bun_{P, N} \times  \Bun_{P, N} 
}
\end{equation} 
where $b_1^P$ is a smooth morphism (which is the forgetful morphism of the level structure on $d \sum x_i$) and $b_2^P$ has zero differential (because it is the composition of the Frobenius morphism with some other morphism). We have $T_{\epsilon_{P, d}}(x_P) = T_{b_1^p}(x_P)$ (see for example \cite{laurent} I. 2. Proposition 1).
It is well-known that $\Bun_{P, N + d \sum x_i} \xrightarrow{b_1^P} \Bun_{P, N}$ is a $P_{d \sum x_i}$-torsor, where $P_{d \sum x_i}$ 
is defined in \ref{subsection-O-N}.
We deduce that $T_{b_1^p}(x_P) = \on{Lie}(P_{d \sum x_i})$.

Similarly, we have a Cartesian square (taking into account that $\wt{\epsilon_{M, d} }^{-1}(y) = \epsilon_{M, d}^{-1}( \ov{\pi^0_d}(y) )$)
\begin{equation}   \label{equation-Cht-M-cartesian-square}
\xymatrixrowsep{2pc}
\xymatrixcolsep{2pc}
\xymatrix{
\wt{\epsilon_{M, d} }^{-1}(y)  \ar[r]  \ar[d]
& \Bun_{M, N + d \sum x_i}    \ar[d]^{(b_1^M, \;  b_2^M)}  \\
\Bun_{M, N}  \ar[r]^{(\on{Id}, \on{Id}) \quad \quad \quad} 
& \Bun_{M, N} \times  \Bun_{M, N} 
}
\end{equation} 
where $b_1^M$ is a smooth morphism (which is the forgetful morphism of the level structure on $d \sum x_i$) and $b_2^M$ has zero differential.
We deduce that $T_{\widetilde{ \epsilon_{M, d} }}(x_M) = T_{b_1^M}(x_M) = \on{Lie}(M_{d \sum x_i}),$
where $M_{d \sum x_i}$ is defined in \ref{subsection-O-N}.

Morphism (\ref{equation-T-P-to-T-M}) is the canonical morphism $\on{Lie}(P_{d \sum x_i}) \rightarrow \on{Lie}(M_{d \sum x_i})$ induced by $P \twoheadrightarrow M$. So it is surjective. We deduce also that the relative tangent space of $\pi_d$ is $\on{Lie}(U_{d \sum x_i})$.
\cqfd

\subsection{Compatibility of the geometric Satake equivalence and parabolic induction}

The goal of this section is to recall (\ref{equation-TC-satake-geo}) and deduce (\ref{equation-inverse-image-S-M-W-equal-S-M-W-theta-0}), which is the key ingredient for the next section.

\sssec{}
We apply Definition \ref{def-Gr-G-I} to $\mb G_m$ and denote by $\Gr_{\mb G_m, I}$ the associated reduced ind-scheme.
We denote by $\rho_G$ (resp. $\rho_M$) the half sum of positive roots of $G$ (resp. $M$). 
Since $2(\rho_G - \rho_M)$ is a character of $M$, 
the morphism $2(\rho_G - \rho_M): M \rightarrow \mb G_m$ induces a morphism $\Gr_{M, I} \rightarrow \Gr_{\mb G_m, I}$ by sending a $M$-bundle $\mc M$ to the $\mb G_m$-bundle $\mc M \overset{M}{\times} \mb G_m$.
We have a morphism $\on{deg}: \Gr_{\mb G_m, I}  \rightarrow \Z$ by taking the degree of a $\mb G_m$-bundle. 
We have the composition of morphisms
\begin{equation}    \label{equation-Gr-M-I-to-Z}
\Gr_{M, I} \rightarrow \Gr_{\mb G_m, I}  \xrightarrow{\on{deg}}  \Z.
\end{equation}

We define $\Gr_{M, I}^n$ as the inverse image of $n \in \Z$. It is open and closed in $\Gr_{M, I}$.
We define $\Gr_{P, I}^n := (\pi^0)^{-1} \Gr_{M, I}^n$. Morphism (\ref{diagram-Gr-G-P-M}) induces a morphism
\begin{equation}
\Gr_{G, I} \xleftarrow{i^0_n} \Gr_{P, I}^n \xrightarrow{\pi^0_n} \Gr_{M, I}^n.
\end{equation}

\sssec{}
Recall that we have defined $\wh \Lambda_{G, P}$ in \ref{subsection-def-Lambda-G-P}. As in \cite{schieder} 2.1.2., we define $\Lambda_{G, P}:=\{ \lambda \in \Lambda_G | \langle  \check{\alpha} , \lambda \rangle =0 \text{ for all } \check{\alpha} \in \wh \Gamma_M\}$. The pairing $\langle \ , \ \rangle$ in \ref{subsection-def-Lambda-G} induces a pairing $\langle \ , \ \rangle: \wh \Lambda_{G, P} \times \Lambda_{G, P} \rightarrow \Z$.

\sssec{}    \label{subsection-W-theta-W-n}
We denote by $\on{Rep}_{\Ql}(\wh M^I)$ the category of finite dimensional $\Ql$-linear representations of $\wh M^I$. 
Let $W \in \on{Rep}_{\Ql}(\wh M^I)$. Then $Z_{\wh M}$ acts on $W$ via $Z_{\wh M} \hookrightarrow \wh M^I$ diagonally. We have the decomposition as $Z_{\wh M}$ representation: $W = \bigoplus_{ \theta \in \Lambda_{Z_{\wh M}} } W^{\theta} .$

Since $\theta \in \Lambda_{Z_{\wh M}} = \wh \Lambda_{G, P}$ and $2(\rho_G - \rho_M) \in \Lambda_{G, P}$, we can consider $\langle \theta , 2(\rho_G - \rho_M) \rangle$. Let $W_n = \bigoplus_{ \langle \theta , 2(\rho_G - \rho_M) \rangle =n } W^{\theta}$. We have $W = \bigoplus_{ n \in \Z } W_n$.

Let $\on{Rep}_{\Ql}(\wh M^I)^{\theta}$ be the category of finite dimensional $\Ql$-linear representations of $\wh M^I$ such that $Z_{\wh M}$ acts by $\theta$. We have $\on{Rep}_{\Ql}(\wh M^I) = \bigoplus_{ \theta \in \Lambda_{Z_{\wh M}} }   \on{Rep}_{\Ql}(\wh M^I)^{\theta}$.
Let $$\on{Rep}_{\Ql}(\wh M^I)_n = \bigoplus_{ \theta \in \Lambda_{Z_{\wh M}}, \, \langle \theta , 2(\rho_G - \rho_M) \rangle =n }  \on{Rep}_{\Ql}(\wh M^I)^{\theta} .$$ We have 
$$ \on{Rep}_{\Ql}(\wh M^I) = \bigoplus_{ n \in \Z }  \on{Rep}_{\Ql}(\wh M^I)_n ,$$

We define $(\on{Res}^{\wh G^I}_{\wh M^I})_n$ to be the composition of morphisms $\on{Rep}_{\Ql}(\wh G^I) \xrightarrow{ \on{Res}^{\wh G^I}_{\wh M^I}  }  \on{Rep}_{\Ql}(\wh M^I) \twoheadrightarrow \on{Rep}_{\Ql}(\wh M^I)_n$. 

\sssec{}
In morphism (\ref{equation-Gr-M-I-to-Z}), $\Gr_{M, I, W^{\theta}}$ is sent to $\langle \theta , 2(\rho_G - \rho_M) \rangle$.
We deduce that $\Gr_{M, I}^n \cap \Gr_{M, I, W} = \Gr_{M, I, W_n}$.

%
%
%
%

\sssec{}
In Corollary \ref{cor-Satake-functor-I}, we have defined a fully faithful functor $Sat_{G, I}: \on{Rep}_{\Ql}(\wh G^I)  \rightarrow \Perv_{G_{I, \infty}}(\Gr_{G, I}, \Ql)$ which sends $W$ to $\mc S_{G, I, W}$. We denote by $\Perv_{G_{I, \infty}}(\Gr_{G, I}, \Ql)^{\on{MV}}$ the subcategory of essential image of this functor. Similarly, we define $Sat_{M, I}: \on{Rep}_{\Ql}(\wh M^I)  \rightarrow \Perv_{M_{I, \infty}}(\Gr_{M, I}, \Ql)^{\on{MV}}$. Let $Sat_{M, I, n}$ be the restriction of $Sat_{M, I}$ to $\on{Rep}_{\Ql}(\wh M^I)_n$.

\begin{thm}  \label{thm-geo-satake-CT-I-paws}  (\cite{bd} 5.3.29, \cite{bg} Theorem 4.3.4, \cite{mv} Theorem 3.6 for $M=T$, \cite{BR} Proposition 15.2)

(a) For any $n \in \Z$, the complex 
$$(\pi_{n}^0)_! (i_{n}^0)^* \mc S_{G, I, W} \otimes \big( \Ql [1](\frac{1}{2}) \big)^{\otimes n }$$ is in $\Perv_{M_{I, \infty}}(\Gr_{M, I}^n, \Ql)^{\on{MV}}  .$

(b) We denote by $((\pi_{n}^0)_! (i_{n}^0)^*)^\sim$ the shifted functor $(\pi_{n}^0)_! (i_{n}^0)^* \otimes \big( \Ql [1](\frac{1}{2}) \big)^{\otimes n }$. Then there is a canonical isomorphism of fiber functors 
\begin{equation}   \label{equation-Sat-M-Res-equal-CT-Sat-G}
Sat_{M, I, n} \circ (\on{Res}^{\wh G^I}_{\wh M^I})_n = ((\pi_{n}^0)_! (i_{n}^0)^*)^\sim \circ Sat_{G, I}.
\end{equation} 
In other words, the following diagram of categories canonically commutes:
\begin{equation}
\xymatrixrowsep{2pc}
\xymatrixcolsep{6pc}
\xymatrix{
\Perv_{G_{I, \infty}}(\Gr_{G, I}, \Ql)^{\on{MV}}   \ar[r]^{  ( (\pi_{n}^0)_! (i_{n}^0)^* )^\sim  } 
& \Perv_{M_{I, \infty}}(\Gr_{M, I}^n, \Ql)^{\on{MV}}  \\
\on{Rep}_{\Ql}(\wh G^I)   \ar[u]^{Sat_{G, I}}   \ar[r]^{ (\on{Res}^{\wh G^I}_{\wh M^I})_n }
& \on{Rep}_{\Ql}(\wh M^I)_n   \ar[u]^{ Sat_{M, I, n} }
}
\end{equation}
\cqfd
\end{thm}


\begin{rem}
The references cited above in Theorem \ref{thm-geo-satake-CT-I-paws} are for the case of affine grassmannians (i.e. $I$ is a singleton). 
The general case (i.e. $I$ is arbitrary) can be deduced from the case of affine grassmannians using the fact that the constant term functor commutes with fusion (i.e. convolution). The proof for $I = \{1, 2\}$ is already included in the proof of Proposition 15.2 in \cite{BR}. For general $I$ the proof is similar. 
\end{rem}


\begin{cor}  \label{cor-satake-geo-TC}  There is a canonical isomorphism
\begin{equation}   \label{equation-TC-satake-geo}
\mc{S}_{M, I, W_{n }} \simeq (\pi_{n}^0)_! (i_{n}^0)^* \mc S_{G, I, W} [ n] (  n/2 ).
\end{equation}
\end{cor}
\dem
Applying (\ref{equation-Sat-M-Res-equal-CT-Sat-G}) to $W$ and taking into account that $\mc S_{M, I, W_n} = Sat_{M, I, n}(W_n)$ and $\mc S_{G, I, W} = Sat_{G, I}(W)$, we deduce (\ref{equation-TC-satake-geo}).
\cqfd

\quad

\sssec{}    \label{subsection-TC-satake-geo-en-quotient}
For any $n$, denote by $\Gr_{P, I, W}^n=\Gr_{P, I}^n \cap \Gr_{P, I, W}$. We have a commutative diagram, where the first line is induced by (\ref{equation-Gr-G-P-M-I-W}).
$$
\xymatrix{
\Gr_{G, I, W}  \ar[d]^{\xi_{G, d}} 
& \Gr_{P, I, W}^n  \ar[l]_{i^0_{n}}   \ar[d]^{\xi_{P, d}}  \ar[r]^{\pi^0_{n}} 
&  \Gr_{M, I, W_{n}}    \ar[d]^{\xi_{M, d}} \\
[G_{I, d} \backslash  \Gr_{G, I, W}]       
& [ P_{I, d} \backslash  \Gr_{P, I, W}^n]   \ar[l]_{ \ov{i^0_{d, n}} }  \ar[r]^{ \ov{\pi^0_{d, n}} } 
&  [  M_{I, d} \backslash  \Gr_{M, I, W_{n}}]
}$$
The right square is not Cartesian.
The morphism $$\Gr_{P, I, W}^n \rightarrow [ P_{I, d} \backslash  \Gr_{P, I, W}^n] \underset{ [  M_{I, d} \backslash  \Gr_{M, I, W_{n}}] }{\times} \Gr_{M, I, W_{n}} =  [ U_{I, d} \backslash  \Gr_{P, I, W}^n]$$ is a $U_{I, d}$-torsor. 
Since the group scheme $U_{I, d}$ is unipotent over $X^I$,
we deduce that
\begin{equation}
(\pi_{n}^0)_! (\xi_{P, d})^* \simeq (\xi_{M, d})^* (\ov{\pi_{d, n}^0})_! [-2m](-m),
\end{equation}
where $m= \dim \xi_{P, d} - \dim \xi_{M, d} = \dim_{X^I} U_{I, d}$. 

Corollary \ref{cor-satake-geo-TC} implies
\begin{equation}     \label{equation-TC-satake-geo-en-quotient}
\mc S_{M, I, W_{n }}^d \isom (\ov{\pi_{d, n}^0})_! (\ov{i_{d, n}^0})^* \mc S_{G, I, W}^d  [n-2m](n/2-m) .
\end{equation}

\quad

\sssec{}
Let $(\omega_i)_{i \in I} \in (\wh \Lambda_M^+)^I$. Let $V^{\omega_i}$ be the irreducible representation of $\wh M$ of highest weight $\omega_i$.
Note that $\wh \Lambda_{G, P} = \wh \Lambda_{M, M}$ (defined in \ref{subsection-def-Lambda-G-P}). By definition, it coincides with $\pi_1(M)$ defined in \cite{var} Lemma 2.2. We denote by $[\sum_{i \in I} \omega_i]$ the image of $\sum_{i \in I} \omega_i$ by the projection $\wh \Lambda_M \twoheadrightarrow \wh \Lambda_{M, M}$.

\begin{lem}  \label{lem-Cht-V-omega-i-non-empty} (\cite{var} Proposition 2.16 d))
$\Cht_{M, N, I, \boxtimes_{i \in I} V^{\omega_i}}$ is non-empty if and only if $[\sum_{i \in I} \omega_i]$ is zero.
\cqfd
\end{lem}

\sssec{}   \label{subsection-image-of-epsilon-M-d-only-in-theta-zero}
Let $W$ and $W^{\theta}$ as in \ref{subsection-W-theta-W-n}.
$W$ has a unique decomposition of the form $$W = \bigoplus_{ (\omega_i)_{i \in I} \in (\wh \Lambda_M^+)^I } \big( \boxtimes_{i \in I} V^{\omega_i} \big) \otimes_{\Ql} \mf M_{(\omega_i)_{i \in I}} ,$$
where $\mf M_{(\omega_i)_{i \in I}} $ are finite dimensional $\Ql$-vector spaces, all but a finite number of them are zero.
We have $$W^{\theta} = \bigoplus_{ (\omega_i)_{i \in I} \in (\wh \Lambda_M^+)^I, \,  [\sum_{i \in I} \omega_i] = \theta} \big( \boxtimes_{i \in I} V^{\omega_i} \big) \otimes_{\Ql} \mf M_{(\omega_i)_{i \in I}} .$$
Lemma \ref{lem-Cht-V-omega-i-non-empty} implies that $\Cht_{M, N, I, W^{\theta}}$ is non-empty if and only if $\theta$ is zero. 
For such $\theta$, we have $\langle \theta, 2(\rho_G - \rho_M)  \rangle =0$. We deduce that $\Cht_{M, N, I, W} = \cup_{n \in \Z} \Cht_{M, N, I, W_n} = \Cht_{M, N, I, W_0}$.
So the image of $$\epsilon_{M, d}: \Cht_{M, N, I, W} \rightarrow [M_{I, d} \backslash \Gr_{M, I, W} ]$$ is in $[M_{I, d} \backslash \Gr_{M, I, W_0} ] $. 

\sssec{}
With the notations of diagram (\ref{diagram-TC-Cht-Gr-W}), we have 
\begin{equation}    \label{equation-inverse-image-S-M-W-equal-S-M-W-theta-0}
\begin{aligned}
(\epsilon_{M, d})^* \mc S_{M, I, W}^d & = (\epsilon_{M, d})^* \mc S_{M, I, W_{0}}^d  \\
&  \isom  (\epsilon_{M, d} )^*   ( \ov{\pi^0_{d, 0}} )_!  ( \ov{i^0_{d, 0}} )^* \mc S_{G, I, W}^d[-2m](-m)  \\
& = (\epsilon_{M, d} )^*  ( \ov{\pi^0_{d}} )_!  ( \ov{i^0_{d}} )^* \mc S_{G, I, W}^d[-2m](-m)  .
\end{aligned}
\end{equation}
The first and third equality follows from \ref{subsection-image-of-epsilon-M-d-only-in-theta-zero}. The second isomorphism follows from (\ref{equation-TC-satake-geo-en-quotient}) applied to $n=0$.

\subsection{Construction of the morphism (\ref{equation-F-G-vers-F-M})}   \label{subsection-CT-complexe}

\sssec{}
Consider diagrams (\ref{diagram-TC-Cht-Gr-W}) and (\ref{diagram-cht-P-cht-M-pi-d}). 
Let 
$m=\dim_{X^I} U_{I, d} $ as in \ref{subsection-TC-satake-geo-en-quotient}. By Lemma \ref{lem-pi-d-est-lisse}, $m= \dim \pi_d $.
We construct a canonical map of functors from $D_c^b([ P_{I, d} \backslash  \Gr_{P, I, W}] , \Ql)$ to $D_c^b(\Cht_{M, N, I, W}, \Ql)$:
\begin{equation}  \label{fleche-de-foncteur-cle}
\pi_! (\epsilon_{P, d})^* \rightarrow  (\epsilon_{M, d} )^* ( \ov{\pi^0_d} )_![-2m](-m)
\end{equation}
as the composition:
\begin{equation}     \label{equation-base-change-with-trace} 
\begin{aligned}
\pi_! (\epsilon_{P, d})^*   & \simeq (\widetilde{  \pi^0_{d}  } )_! (\pi_d)_! (\pi_d)^* (\widetilde{ \epsilon_{M, d } })^*\\
& \rightarrow (\widetilde{  \pi^0_{d}  } )_! (\widetilde{ \epsilon_{M, d } })^* [-2m](-m) \overset{\thicksim} {\leftarrow}  (\epsilon_{M, d} )^* ( \ov{\pi^0_d} )_![-2m](-m)
\end{aligned}
\end{equation}

The second morphism in (\ref{equation-base-change-with-trace}) is induced by the isomorphism $(\pi_d)^* [2m](m) \simeq (\pi_d)^!$ (because $\pi_d$ is smooth) and the counit map $\on{Co}: (\pi_d)_! (\pi_d)^! \rightarrow \Id$. (The composition $(\pi_d)_! (\pi_d)^*  [2m](m) \isom (\pi_d)_! (\pi_d)^! \rightarrow \Id$ is the trace map in \cite{sga4}  XVIII 2.)

The third morphism is the proper base change (\cite{sga4} XVII 5, \cite{LO08} Theorem 12.1).

\sssec{}   \label{subsection-pi-i-F-G-to-F-M}
Now we construct a morphism of complexes in $D_c^b(\Cht_{M, N, I, W}, \Ql)$:
\begin{equation}  \label{equation-F-G-vers-F-M-long}
\begin{aligned}
\pi_! i^* \mc{F}_{G, N, I, W} & = \pi_! i^* (\epsilon_{G, d})^* \mc S_{G, I, W}^d \\ 
& ^{(a)}\!{\simeq}  \pi_! (\epsilon_{P, d})^* ( \ov{i^0_d} )^* \mc S_{G, I, W}^d \\
& ^{(b)}\!{\rightarrow}   (\epsilon_{M, d} )^* ( \ov{\pi^0_d} )_!  ( \ov{i^0_d} )^* \mc S_{G, I, W}^d[-2m](-m) \\
& ^{(c)}\!{ \overset{\thicksim} {\leftarrow} }  (\epsilon_{M, d})^* \mc S_{M, I, W}^d = \mc{F}_{M, N, I, W }
\end{aligned}
\end{equation}

(a) is induced by the commutativity of diagram (\ref{diagram-TC-Cht-Gr-W}).
(b) is induced by morphism (\ref{fleche-de-foncteur-cle}).
(c) is (\ref{equation-inverse-image-S-M-W-equal-S-M-W-theta-0}).

\quad

\sssec{}
All the constructions in 3.1-3.3 are compatible with the quotient by $\Xi$. In particular, just as in \ref{subsection-diagram-TC-Cht-Gr-W}, diagram (\ref{diagram-TC-Cht-Gr-Xi}) induces a commutative diagram:
\begin{equation}   \label{diagram-TC-Cht-Gr-W-Xi}
\xymatrix{
\Cht_{G, N, I, W} / \Xi \ar[d]^{\epsilon_{G, N, I, d}^{\Xi} } 
&\Cht_{P, N, I, W} / \Xi   \ar[l]_{i}   \ar[d]^{\epsilon_{P, N, I, d}^{\Xi}}  \ar[r]^{\pi} 
&  \Cht_{M, N, I, W} / \Xi   \ar[d]^{\epsilon_{M, N, I, d}^{\Xi}} \\
[G_{I, d}^{ad} \backslash  \Gr_{G, I, W}]       
& [ \ov P_{I, d} \backslash  \Gr_{P, I, W}]   \ar[l]_{ \ov{i^0_d} }  \ar[r]^{ \ov{\pi^0_d} } 
&  [  \ov M_{I, d} \backslash  \Gr_{M, I, W}]
}
\end{equation}

\begin{construction}  \label{constr-funtor-TC-for-complex}
Just as in \ref{subsection-pi-i-F-G-to-F-M} (using (\ref{diagram-TC-Cht-Gr-W-Xi}) instead of (\ref{diagram-TC-Cht-Gr-W})), we construct a canonical morphism of complexes in $D_c^b(\Cht_{M, N, I, W} / \Xi, \Ql)$:
\begin{equation}  
\pi_! i^* \mc{F}_{G, N, I, W}^{\Xi} \rightarrow \mc{F}_{M, N, I, W}^{\Xi}.
\end{equation} 
\end{construction}

\subsection{More on cohomology groups}     \label{subsection-more-on-H-M-prime}

When the level structure $N$ is non-empty, to construct the constant term morphism of cohomology groups, we need a variant of $H_{M, N, I, W}^j$.

\sssec{}
Let $\mc O_N$ be the ring of functions on $N$ as in \ref{subsection-O-N}. 
The finite group $G(\mc O_N)$ (resp. $P(\mc O_N)$ and $M(\mc O_N)$) acts on $\Cht_{G, N, I, W}$ (resp. $\Cht_{P, N, I, W}$ and $\Cht_{M, N, I, W}$) by changing the level structure on $N$: $g \in G(\mc O_N)$ sends a level structure $\psi_G$ to $g^{-1} \circ \psi_G$. 

By \cite{var} Proposition 2.16 b), $\Cht_{G, N, I, W}$ (resp. $\Cht_{P, N, I, W}$ and $\Cht_{M, N, I, W}$) is a finite étale Galois cover of 
$\restr{\Cht_{G, I, W}}{(X \sm N)^I}$ (resp. $\restr{\Cht_{P, I, W}}{(X \sm N)^I}$ and $\restr{\Cht_{M, I, W}}{(X \sm N)^I}$)
with Galois group $G(\mc O_N)$ (resp. $P(\mc O_N)$ and $M(\mc O_N)$).

\begin{defi}    \label{def-Cht-P-M-'}
We define
$$\Cht_{P, N, I, W}' := \Cht_{P, N, I, W} \overset{P(\mc O_N)} \times G(\mc O_N) , \quad
\Cht_{M, N, I, W}' := \Cht_{M, N, I, W} \overset{P(\mc O_N)} \times G(\mc O_N) ,$$
where $P(\mc O_N)$ acts on $G(\mc O_N)$ by left action (by left multiplication) and $P(\mc O_N)$ acts on $\Cht_{M, N, I, W}$ via the quotient $P(\mc O_N) \twoheadrightarrow M(\mc O_N)$.
\end{defi}

\sssec{}   \label{subsection-Cht-G-P-M-prime-morphisms}
Morphisms (\ref{equation-Cht-G-P-M-I-W}) induce morphisms
\begin{equation}    \label{equation-Cht-G-P-M-prime}
\Cht_{G, N, I, W} \xleftarrow{i'} \Cht_{P, N, I, W}' \xrightarrow{\pi'} \Cht_{M, N, I, W}'.
\end{equation}
Indeed, the morphism $i'$ is giving by 
$$\big( (\mc P, \psi_P) \rightarrow (\ta \mc P, \ta \psi_P) , g \in G(\mc O_N)  \big) \mapsto \big( (\mc G, g^{-1} \circ \psi_G) \rightarrow  (\ta \mc G, g^{-1} \circ \ta \psi_G) \big) $$
where $\mc G = \mc P \overset{P}{\times} G$ and $\psi_G = \psi_P \overset{P}{\times} G$.
The morphism $\pi'$ is induced by $\pi$, which is $P(\mc O_N)$-equivariant (because $P(\mc O_N)$ acts on $\Cht_{P, N, I, W}$ and $\Cht_{M, N, I, W}$ by changing the level structure on $N$).

\begin{rem}  \label{rem-Cht-P-prime-is-fiber-product}
The morphism $\Cht_{P, N, I, W}' \rightarrow \Cht_{P, I, W} \underset{\Cht_{G, I, W}}\times \Cht_{G, N, I, W}$ is a $G(\mc O_N)$-equivariant morphism of $G(\mc O_N)$-torsors over $\Cht_{P, I, W}$, thus it is an isomorphism.
In \cite{var} 2.28, the stack $\Cht_{P, N, I, W}'$ is denoted by $FBun_{P, D, n, \bar{\omega}}$. The reason why we will need $\Cht_{P, N, I, W}'$ instead of $\Cht_{P, N, I, W}$ is justified in Example \ref{exemple-TC-cht-sans-patte} and Theorem \ref{thm-i-contract-pour-mu-grand}.
\end{rem}

\begin{defi}
We define $$\Cht_{P, N, I, W}^{' \, \leq \mu} := \Cht_{P, N, I, W}^{\leq \mu} \overset{P(\mc O_N)} \times G(\mc O_N) , $$
$$\Cht_{M, N, I, W}^{' \, \leq \mu} := \Cht_{M, N, I, W}^{\leq \mu} \overset{P(\mc O_N)} \times G(\mc O_N) , \quad \Cht_{M, N, I, W}^{' \, \leq \mu, \, \nu} := \Cht_{M, N, I, W}^{\leq \mu, \, \nu} \overset{P(\mc O_N)} \times G(\mc O_N)   .$$
\end{defi}

\sssec{}
We have a commutative diagram of algebraic stacks:
\begin{equation}      
\xymatrixrowsep{1pc}
\xymatrixcolsep{2pc}
\xymatrix{
&\Cht_{P, N, I, W}  \ar[ld]_{i}   \ar[r]^{\pi}   \ar[d]
&  \Cht_{M, N, I, W}  \ar[d]  \\
\Cht_{G, N, I, W} \ar[d]
&\Cht_{P, N, I, W}'  \ar[l]_{i'}   \ar[d] \ar[r]^{\pi'} 
&  \Cht_{M, N, I, W}'  \ar[d] \\
\Bun_{G} 
&\Bun_{P}  \ar[l]   \ar[r]
&  \Bun_{M} 
}
\end{equation}
We deduce that $\Cht_{P, N, I, W}^{' \, \leq \mu}$ is also the inverse image of $\Bun_P^{\leq \mu}$ by $\Cht_{P, N, I, W}' \rightarrow \Bun_P$ and $\Cht_{M, N, I, W}^{' \, \leq \mu}$ (resp. $\Cht_{M, N, I, W}^{' \, \leq \mu, \, \nu}$) is also the inverse image of $\Bun_M^{\leq \mu}$ (resp. $\Bun_M^{\leq \mu, \, \nu}$) by $\Cht_{M, N, I, W}' \rightarrow \Bun_M$.

\begin{defi}    \label{def-mc-H-M-prime}
Just as in Section \ref{subsection-coho-Cht-M}, 
we construct a morphism
$\epsilon_{M, d}^{\Xi \, '}: \Cht_{M, N, I, W}' / \Xi \rightarrow  [  \ov{M}_{I, d} \backslash  \Gr_{M, I, W}]$
and we define $\mc F_{M, N, I, W}^{ \; ' \; \Xi}$ to be the inverse image of $\mc S_{M, I, W}^{\on{ad}, d}$. 
We define $\mc H _{M, N, I, W}^{' \, \leq \mu, \; \nu}  := R(\mf{p}_M)_! ( \restr{ \mc{F}_{M, I, N, W}^{ \; ' \, \Xi}    } {\Cht_{M, N, I, W}^{' \, \leq  \mu, \; \nu} / \Xi} ) $,  
$\mc H _{M, N, I, W}^{' \, j, \; \leq \mu, \; \nu}  := H^j (\mc H _{M, N, I, W}^{' \, \leq \mu, \; \nu}  )$ and
$H_{M, N, I, W}^{' \, j, \, \leq \mu, \, \nu} := \restr{  \mc H _{M, N, I, W}^{' \; j, \; \leq \mu, \; \nu}    }{  \ov{\eta^I}  }. $
\end{defi}

\sssec{} \label{subsection-H-M-leq-mu-support-on-cone}
Just as in \ref{subsection-nu-not-in-cone-H-M-mu-nu-nulle}, if $\nu \notin \wh{\Lambda}_{Z_M / Z_G}^{\mu}$ (defined in \ref{subsection-def-translated-cone-ad}), then $\Cht_{M, N, I, W}^{' \, \leq \mu, \, \nu}$ is empty and $H_{M, N, I, W}^{' \, j, \, \leq \mu, \, \nu}=0$.

\begin{defi}   \label{def-H-M-prime-j}
Just as in Definition \ref{def-H-M-j-leq-mu} and Definition \ref{def-H-M-j}, we define
$$H_{M, N, I, W}^{' \, j, \, \leq \mu} := \prod_{\nu \in \wh{\Lambda}_{Z_M / Z_G}^{\mu}} H_{M, N, I, W}^{' \, j, \, \leq \mu, \, \nu} ; \quad  {H'}_{\! \! \! M, N, I, W}^{\; j} : = \varinjlim_{\mu} {H'}_{\! \! \! M, N, I, W}^{\; j, \; \leq\mu}   .$$
\end{defi}

\begin{defi}    \label{def-H-M-j-nu-prime}
For any $\nu \in \wh{\Lambda}_{Z_M / Z_G}^{\Q}$, we define $H _{M, N, I, W}^{' \, j, \, \nu}: = \varinjlim _{\mu}  H _{M, N, I, W}^{' \, j, \, \leq\mu, \, \nu}$. 
\end{defi}

\subsection{Constant term morphism for cohomology groups}

\sssec{}

Morphisms (\ref{diagram-TC-Cht-Gr-W-Xi}) induce morphisms over $(X \sm N)^I$:
\begin{equation}     \label{diagram-TC-Cht-Gr-W-Xi-prime} 
\xymatrix{
\Cht_{G, N, I, W} / \Xi \ar[d]^{\epsilon_{G, N, I, d}^{\Xi} } 
&\Cht_{P, N, I, W}' / \Xi   \ar[l]_{i'}   \ar[d]^{\epsilon_{P, N, I, d}^{' \, \Xi}}  \ar[r]^{\pi'} 
&  \Cht_{M, N, I, W}' / \Xi   \ar[d]^{\epsilon_{M, N, I, d}^{' \, \Xi}} \\
[G_{I, d}^{ad} \backslash  \Gr_{G, I, W}]       
& [ \ov P_{I, d} \backslash  \Gr_{P, I, W}]   \ar[l]_{ \ov{i^0_d} }  \ar[r]^{ \ov{\pi^0_d} } 
&  [  \ov M_{I, d} \backslash  \Gr_{M, I, W}]
}
\end{equation}


\sssec{}   \label{subsection-Cht-G-P-M-prime-leq-mu}

For any $\mu\in \wh{\Lambda}_{G^{\mr{ad}}}^{+, \Q}$ and any $\nu \in \wh \Lambda_{Z_M / Z_G}^{\Q}$, the first line of morphisms (\ref{diagram-TC-Cht-Gr-W-Xi-prime}) induces morphisms over $(X \sm N)^I$:
\begin{equation}     \label{equation-Cht-G-P-M-prime-leq-mu}
\Cht_{G, N, I, W}^{\leq \mu} / \Xi   \xleftarrow{i'}     \Cht_{P, N, I, W}^{' \, \leq \mu, \, \nu} / \Xi    \xrightarrow{\pi'}   \Cht_{M, N, I, W}^{' \, \leq \mu, \, \nu} / \Xi 
\end{equation}


\quad

The proof of \cite{var} Proposition 5.7 in fact proves the following:
\begin{prop}  \label{prop-Varshavsky-proper} (\cite{var} Proposition 5.7)
For any $\mu\in \wh{\Lambda}_{G^{\mr{ad}}}^{+, \Q}$ and any $\nu \in \wh \Lambda_{Z_M / Z_G}^{\mu}$ (defined in \ref{subsection-def-translated-cone-ad}), there exists an open dense subscheme $\Omega^{\leq \mu, \nu}$ of $(X \sm N)^I$ such that the restriction of the morphism $i'$ on $\restr{\Cht_{P, N, I, W}^{' \, \leq \mu, \, \nu} / \Xi}{ \Omega^{\leq \mu, \nu}  }$ is proper. In particular, the restriction of the morphism $i'$ on $\restr{\Cht_{P, N, I, W}^{' \, \leq \mu, \, \nu} / \Xi}{ \eta^I  }$ is proper.
\cqfd
\end{prop}

\begin{rem}   \label{rem-i-proper-schematic}
In $loc.cit.$, the level is denoted by $D$, the paws are indexed by $n$, 
the index $d$ is related to our $\nu$, the index $\ov{k}$ is related to our $W$, the index $[g]$ is in $G(\mc O_N) / P(\mc O_N)$. The open subscheme $\Omega^{\leq \mu, \nu} $ is of the form
$$\Omega(m) = \{  (x_i)_{i \in I} \in (X \sm N)^I, \; x_i \neq {}^{\tau^r} x_j \text{ for all } i, j \text{ and } r = 1, 2, \cdots, m \},$$
where ${}^{\tau^r} x$ is the image of $x$ by $\Frob^r: X \rightarrow X$ and $m$ is some positive integer. 

In the proof of $loc.cit.$, $\Bun_G^{\leq \mu}$ is denoted by $V$ and $\Omega(m)$ is denoted by $U$.
Varshavsky shows that for fixed $\mu$ and $\nu$, there exists a level $D$ large enough and an integer $m$ large enough (both depending on $\mu$ and $\nu$), such that over $\Bun_G^{\leq \mu} \times \Omega(m) \subset \Bun_G \times (X \sm N)^I$, the morphism
$$\restr{ \Cht_{P, D, I, W}^{' \, \leq \mu, \, \nu} }{ \Omega(m)  } \rightarrow \restr{ \Cht_{G, D, I, W}^{\leq \mu} }{ \Omega(m)  } $$
is a closed embedding. In particular, it is proper. Then we descend to level $N$.

Note that $i'$ is schematic (i.e. representable). This is implied by the well-known fact that $\Bun_P \rightarrow \Bun_G$ is schematic (a $P$-structure of a $G$-bundle $\mc G$ over $X \times S$ is a section of the fibration $\mc G / P \rightarrow X \times S$).
\end{rem}


\sssec{}
Now consider the following commutative diagram:
\begin{equation}   \label{diagram-Cht-G-P-M-leq-mu-on-open}
\xymatrixrowsep{1pc}
\xymatrixcolsep{1pc}
\xymatrix{
& \restr{\Cht_{P, N, I, W}^{' \, \leq \mu, \, \nu} / \Xi}{\Omega^{\leq \mu, \, \nu}} \ar[dl]_{i'}  \ar[dr]^{\pi'}  \ar[dd]^{\mf p_P}
&  \\
\restr{\Cht_{G, N, I, W}^{\leq \mu} / \Xi}{\Omega^{\leq \mu, \, \nu}} \ar[dr]_{\mf{p}_G}
& 
& \restr{\Cht_{M, N, I, W}^{' \, \leq \mu, \, \nu} / \Xi}{\Omega^{\leq \mu, \, \nu}} \ar[dl]^{\mf{p}_M} \\
& \Omega^{\leq \mu, \, \nu}
&
}
\end{equation}

To simplify the notations, we denote by $\mc F_{G, N, \Omega^{\leq \mu, \, \nu}, W}^{\; \Xi} $ the restriction of $\mc F_{G, I, N, W}^{\; \Xi} $ to $\restr{\Cht_{G, N, I, W}^{\leq \mu} / \Xi}{\Omega^{\leq \mu, \, \nu}} $ 
and by $\mc F_{M, N, \Omega^{\leq \mu, \, \nu}, W}^{\; ' \; \Xi, \, \nu} $ the restriction of $\mc F_{M, I, N, W}^{\; ' \; \Xi} $ to $\restr{\Cht_{M, N, I, W}^{ ' \; \leq \mu, \, \nu} / \Xi}{\Omega^{\leq \mu, \, \nu}} $.

The commutative diagram (\ref{diagram-TC-Cht-Gr-W-Xi-prime}) is compatible with the Harder-Narasimhan stratification. Just as in Construction \ref{constr-funtor-TC-for-complex}, we construct a canonical morphism of complexes in $D_c^b(\restr{\Cht_{M, N, I, W}^{' \, \leq \mu, \, \nu} / \Xi}{\Omega^{\leq \mu, \, \nu}}, \Ql)$:
\begin{equation}   \label{morphism-TC-F-G-F-M-leq-mu-on-U}
 (\pi')_! (i')^*  \mc F_{G, N, \Omega^{\leq \mu, \, \nu}, W}^{\; \Xi}  \rightarrow  \mc F_{M, N, \Omega^{\leq \mu, \, \nu}, W}^{ \; ' \; \Xi, \, \nu} . 
\end{equation}

\quad

\sssec{}
Thanks to Proposition \ref{prop-Varshavsky-proper}, we can apply \cite{sga5} III 3 to diagram (\ref{diagram-Cht-G-P-M-leq-mu-on-open}) and the cohomological correspondence (\ref{morphism-TC-F-G-F-M-leq-mu-on-U}). 

Concretely, firstly we have morphisms of functors from $D_c^b(\restr{\Cht_{G, N, I, W}^{\leq \mu} / \Xi}{\Omega^{\leq \mu, \, \nu}}, \Ql)$ to $D_c^b(\Omega^{\leq \mu, \, \nu}, \Ql)$ (all functors are considered as derived functors):
\begin{equation}    \label{morphism-TC-complex-adj}
(\mf{p}_G)_!  \overset{(a)} \rightarrow (\mf{p}_G)_! (i')_* (i')^* \overset{(b)} \simeq   (\mf{p}_G)_! (i')_! (i')^*  \overset{(c)} \simeq   (\mf{p}_M)_! (\pi')_! (i')^* 
\end{equation}
where (a) is the adjunction morphism, (b) is induced by $i'_! \isom i'_*$ which is because that $i'$ is schematic and proper (Proposition \ref{prop-Varshavsky-proper}), (c) is induced by the commutativity of diagram (\ref{diagram-Cht-G-P-M-leq-mu-on-open}).

Secondly we combine (\ref{morphism-TC-complex-adj}) with (\ref{morphism-TC-F-G-F-M-leq-mu-on-U}). We obtain a composition of morphisms of complexes in $D_c^b(\Omega^{\leq \mu, \, \nu}, \Ql)$:
\begin{equation}   \label{morphism-TC-complex-p-G-F-G-vers-p-M-F-M}
(\mf{p}_G)_!  \mc F_{G, N, \Omega^{\leq \mu, \, \nu}, W}^{\; \Xi}   \xrightarrow{(\ref{morphism-TC-complex-adj})}  (\mf{p}_M)_!  (\pi')_! (i')^*  \mc F_{G, N, \Omega^{\leq \mu, \, \nu}, W}^{\; \Xi}   \xrightarrow{(\ref{morphism-TC-F-G-F-M-leq-mu-on-U})}  (\mf{p}_M)_!  \mc F_{M, N, \Omega^{\leq \mu, \, \nu}, W}^{\; ' \; \Xi, \, \nu} .
\end{equation}
By Definition \ref{def-mc-H-G-leq-mu} and Definition \ref{def-mc-H-M-prime}, (\ref{morphism-TC-complex-p-G-F-G-vers-p-M-F-M}) is also written as:
\begin{equation}     \label{equation-mc-H-G-to-mc-H-M-prime-leq-mu-nu}
\mc C_G^{P, \, \leq \mu, \, \nu}: \restr{ \mc H_{G, N, I, W}^{\leq \mu} }{\Omega^{\leq \mu, \, \nu}} \rightarrow \restr{ {\mc H '}_{\! \! \! M, N, I, W}^{ \; \leq \mu, \, \nu}   }{\Omega^{\leq \mu, \, \nu}}.
\end{equation}

\quad

\sssec{}
From now on, we restrict everything to the geometric generic point $\ov{\eta^I}$ of $X^I$ fixed in \ref{subsection-eta-I-bar}.
Recall that we have defined $H_{G, N, I, W}^{j, \, \leq \mu} = \restr{  \mc H _{G, N, I, W}^{j, \; \leq \mu}    }{  \ov{\eta^I}  } $ in Definition \ref{def-H-G-j} and $H _{M, N, I, W}^{' \, j, \, \leq  \mu, \; \nu} = \restr{  \mc H _{M, N, I, W}^{' \, j, \, \leq  \mu, \; \nu}  }{  \ov{\eta^I}  }$ in Definition \ref{def-mc-H-M-prime}.

For any $j \in \Z$, morphism (\ref{equation-mc-H-G-to-mc-H-M-prime-leq-mu-nu}) induces a morphism of cohomology groups
\begin{equation}     \label{equation-def-CT-leq-mu-nu}
C_G^{P, \, j, \, \leq \mu, \, \nu}: H _{G, N, I, W}^{j, \; \leq \mu}  \rightarrow    {H'}_{\! \! \! M, N, I, W}^{\; j, \; \leq \mu, \; \nu} . 
\end{equation}
By \ref{subsection-H-M-leq-mu-support-on-cone}, for $\nu \notin \wh{\Lambda}_{Z_M / Z_G}^{\mu}$, the morphism $C_G^{P, \, j, \, \leq \mu, \, \nu}$ is the zero morphism.

\sssec{}
We define a morphism:
\begin{equation}    \label{equation-def-CT-leq-mu}
C_G^{P, \, j, \, \leq \mu} = \prod_{\nu \in \wh{\Lambda}_{Z_M / Z_G}^{\mu} } C_G^{P, \, j, \, \leq \mu, \, \nu}: H _{G, N, I, W}^{j, \; \leq \mu}  \rightarrow    {H'}_{\! \! \! M, N, I, W}^{\; j, \; \leq \mu}  
\end{equation}
where ${H'}_{\! \! \! M, N, I, W}^{\; j, \; \leq \mu}  $ is defined in Definition \ref{def-H-M-prime-j}.

\sssec{}
Let $\mu_1, \mu_2 \in \wh{\Lambda}_{G^{\mr{ad}}}^{+, \Q}$ with $\mu_1 \leq \mu_2$. By Lemma \ref{lem-TC-suite-exacte-longue}, the commutative diagram of stacks
\begin{equation}
\xymatrix{
\Cht_{G, N, I, W}^{\leq \mu_2}  / \Xi
&  \Cht_{P, N, I, W}^{' \, \leq \mu_2}  / \Xi  \ar[l]_{i'}  \ar[r]^{\pi'}  
& \Cht_{M, N, I, W}^{' \, \leq \mu_2}  / \Xi \\
\Cht_{G, N, I, W}^{\leq \mu_1} / \Xi \ar@{^{(}->}[u]
& \Cht_{P, N, I, W}^{' \, \leq \mu_1}    / \Xi  \ar[l]_{i'}  \ar[r]^{\pi'}  \ar@{^{(}->}[u]
& \Cht_{M, N, I, W}^{' \, \leq \mu_1}  / \Xi  \ar@{^{(}->}[u]
}
\end{equation}
induces a commutative diagram of cohomology groups:
\begin{equation}  \label{equation-TC-mu-1-mu-2}
\xymatrix{
H _{G, N, I, W}^{j, \; \leq \mu_1} \ar[r]  \ar[d]^{  C_G^{P, \, j, \, \leq \mu_1} } 
&  H _{G, N, I, W}^{j, \; \leq \mu_2} \ar[d]^{  C_G^{P, \, j, \, \leq \mu_2} }   \\
  {H'}_{\! \! \! M, N, I, W}^{\; j, \; \leq \mu_1}    \ar[r]   
&    {H'}_{\! \! \! M, N, I, W}^{\; j, \; \leq \mu_2}  . 
}
\end{equation}

We have defined $H _{G, N, I, W}^{j}  = \varinjlim_{\mu}  H _{G, N, I, W}^{j, \; \leq\mu} $ in Definition \ref{def-H-G-j} and ${H'}_{\! \! \! M, N, I, W}^{\; j}  = \varinjlim_{\mu} {H'}_{\! \! \! M, N, I, W}^{\; j, \; \leq\mu}$ in Definition \ref{def-H-M-prime-j}.
The commutative diagram (\ref{equation-TC-mu-1-mu-2}) induces a morphism between inductive limits:

\begin{defi}   \label{def-CT-cohomology}
For all parabolic subgroups $P$, for all degrees $j \in \Z$, we define the constant term morphism of cohomology groups:
\begin{equation}   \label{equation-C-G-N-P-j}
 C_{G, N}^{P, \, j}: H _{G, N, I, W}^{j}  \rightarrow    {H'}_{\! \! \! M, N, I, W}^{\; j}  . 
\end{equation}
\end{defi}

\quad

\begin{rem}   \label{rem-lim-prod-vs-prod-lim}
The morphisms ${H'}_{\! \! \! M, N, I, W}^{\; j, \; \leq\mu, \; \nu}  \rightarrow \varinjlim_{\mu'}  {H'}_{\! \! \! M, N, I, W}^{\; j, \; \leq\mu', \; \nu} $ for each $\nu \in \wh{\Lambda}_{Z_M / Z_G}^{\Q}$ induce a morphism
\begin{equation}    \label{equation-H-M-lim-prod-to-prod-lim}
\varinjlim_{\mu} \prod_{\nu \in \wh{\Lambda}_{Z_M / Z_G}^{\Q} }  {H'}_{\! \! \! M, N, I, W}^{\; j, \; \leq\mu, \; \nu}  \rightarrow  \prod_{\nu \in \wh{\Lambda}_{Z_M / Z_G}^{\Q} }  \varinjlim_{\mu'}  {H'}_{\! \! \! M, N, I, W}^{\; j, \; \leq\mu', \; \nu}
\end{equation}
With the notations in Definition \ref{def-H-M-prime-j} and Definition \ref{def-H-M-j-nu-prime}, morphism (\ref{equation-H-M-lim-prod-to-prod-lim}) is the natural map
\begin{equation}    \label{equation-H-M-to-prod-H-M-nu}
H_{M, N, I, W}^{' \, j} \rightarrow \prod_{\nu \in \wh{\Lambda}_{Z_M / Z_G}^{\Q} }  {H'}_{\! \! \! M, N, I, W}^{\; j, \; \nu}
\end{equation}
For each $\nu \in \wh{\Lambda}_{Z_M / Z_G}^{\Q}$, taking inductive limit over $\mu$ of (\ref{equation-def-CT-leq-mu-nu}), we define $C_{G, N}^{P, \, j, \, \nu}: H _{G, N, I, W}^{j}  \rightarrow  {H'}_{\! \! \! M, N, I, W}^{\; j, \; \nu}$. We form a morphism
\begin{equation}   \label{equation-prod-nu-C-G-P-nu}
\prod_{\nu \in \wh{\Lambda}_{Z_M / Z_G}^{\Q} } C_{G, N}^{P, \, j, \, \nu} :  H _{G, N, I, W}^{j}    \rightarrow  \prod_{\nu \in \wh{\Lambda}_{Z_M / Z_G}^{\Q} }  {H'}_{\! \! \! M, N, I, W}^{\; j, \; \nu}     .
\end{equation}
It is equal to the composition of (\ref{equation-C-G-N-P-j}) and (\ref{equation-H-M-to-prod-H-M-nu}).

In Lemma \ref{lem-H-M-leq-mu-nu-to-H-M-nu-Xi-G-inj} below, we will prove that for $\mu$ large enough, $H_{M, N, I, W}^{' \, j, \, \leq \mu} \rightarrow \prod_{\nu \in \wh{\Lambda}_{Z_M / Z_G}^{\Q} }  {H'}_{\! \! \! M, N, I, W}^{\; j, \; \nu}$ is injective. This implies that (\ref{equation-H-M-to-prod-H-M-nu}) is injective. Thus the kernel of (\ref{equation-prod-nu-C-G-P-nu}) is the same as the kernel of (\ref{equation-C-G-N-P-j}).
\end{rem}


\begin{rem}
Now consider all parabolic subgroups (not only the standard ones). If $P_1$ and $P_2$ are conjugated, then the conjugation induces an isomorphism $M_1 \simeq M_2$. This induces for any $j$ an isomorphism ${H'}_{\! \! \! M_1, N, I, W}^{\; j}  \simeq  {H'}_{\! \! \! M_2, N, I, W}^{\; j} $. The following diagram commutes
$$\xymatrixrowsep{1pc}
\xymatrixcolsep{3pc}
\xymatrix{
H _{G, N, I, W}^{j}  \ar[r]^{C_{G, N}^{P_1, \, j}}  \ar[rd]_{C_{G, N}^{P_2, \, j}}
&  {H'}_{\! \! \! M_1, N, I, W}^{\; j}  \ar[d]_{\simeq}  \\
&  {H'}_{\! \! \! M_2, N, I, W}^{\; j} 
}
$$
thus we have $\Ker C^{P_1, \,  j }_{G, N} = \Ker C^{P_2, \,  j }_{G, N} $ in $H _{G, N, I, W}^{j}$.

However, we don't know how to compare the constant term morphism along different parabolic subgroups which have a common Levi subgroup. It is perhaps possible to do that, but quite difficult because it would be a generalization of the functional equation for Eisenstein series.
\end{rem}

\begin{defi}  \label{def-cusp-coho}
For any degree $j \in \Z$, we define the cuspidal cohomology group:
\begin{equation}
 H_{G, N, I, W}^{j, \; \on{cusp} }:= \underset{P \varsubsetneq G} \bigcap   \Ker C^{P, \,  j }_{G, N} .
\end{equation}
This is a $\Ql$-vector subspace of $H _{G, N, I, W}^{j} $.
\end{defi}

\begin{rem} 
For
$$  
\xymatrixrowsep{0.5pc}
\xymatrixcolsep{0.5pc}
\xymatrix{
& & P_1 \ar[rd]  \ar[ld] & & \\
& P_2  \ar[rd]  \ar[ld]  & & P_{1, 2}  \ar[rd]  \ar[ld] & \\
G & & M_2 & & M_1
}
$$
we have $C_{M_2, N}^{P_{1, 2}, \, j} \circ C_{G, N}^{P_2, \, j} = C_{G, N}^{P_1, \, j} .$
Thus we have an equivalent definition:
$$H_{G, N, I, W}^{j, \; \on{cusp} } = \underset{P \text{ maximal parabolic} } \bigcap   \Ker C^{P, \,  j }_{G, N} .$$
\end{rem}

\begin{example}   \label{exemple-TC-cht-sans-patte}
(Shtukas without paws)

When $I = \emptyset$ and $W=\bf 1$, we have $\Cht_{G, N, \emptyset, \bf 1} = G(F) \backslash G(\mb A) / K_{G, N}$. (Note that $G$ is split. See \cite{vincent} (0.5) and Remarque 8.21 for more details.) Moreover, let $K_{P, N} := K_{G, N} \cap P(\mb O)$, $K_{U, N} := K_{G, N} \cap U(\mb O)$ and $K_{M, N} := K_{P, N} / K_{U, N}$. We write $\overset{\on{set}}{=}$ for equalities of sets which are not equalities of groupoids.
\begin{align*}
\Cht_{P, N, \emptyset, \bf 1}' & = ( P(F) \backslash P(\mb A) / K_{P, N} ) \overset{P(\mc O_N)} \times G(\mc O_N) 
 = P(F) \backslash \big(   P(\mb A) \overset{P(\mb O)} \times G(\mb O)   \big) / K_{G, N} \\
& = P(F) \backslash G(\mb A) / K_{G, N}   \\
\Cht_{M, N, \emptyset, \bf 1}' & = ( M(F) \backslash M(\mb A) / K_{M, N} ) \overset{P(\mc O_N)} \times G(\mc O_N)  \\
& = M(F)  \backslash \big(   M(\mb A) \overset{P(\mb O) / K_{U, N}} \times G(\mb O) / K_{G, N}  \big)   \overset{\on{set}}{=} M(F)  \backslash \big(   M(\mb A) \overset{P(\mb O)} \times G(\mb O)   \big) / K_{G, N} \\
& =  M(F)  U(\mb A)  \backslash \big(   P(\mb A) \overset{P(\mb O)} \times G(\mb O)   \big) / K_{G, N} = M(F) U(\mb A) \backslash G(\mb A) / K_{G, N} .
\end{align*}


In this case, $\Gr_{P, \emptyset, \bf 1} = \Gr_{M, \emptyset, \bf 1} = \on{Spec} \Fqbar$. We can choose $d=0$ in (\ref{diagram-cht-P-cht-M-pi-d}). Thus $\wt{\Cht}_{M, N, \emptyset, 1}' = \Cht_{M, N, \emptyset, 1}' $. 
The constant term morphism $C_{G, N}^{P, \, j}$ in Definition \ref{def-CT-cohomology} coincides (up to constants depending on 
$\nu \in \wh{\Lambda}_{Z_M / Z_G}^{\Q}$ component by component) with the classicial constant term morphism:
\begin{equation}   \label{equation-classical-CT}
\begin{aligned}
C_c(G(F) \backslash G(\mb A) / K_{G, N} \Xi, \Ql ) & \rightarrow C( U(\mb A) M(F) \backslash G(\mb A) / K_{G, N}  \Xi, \Ql ) \\
f \quad \quad & \mapsto \quad  f^P: g \mapsto \int_{ {U(F) \backslash  U( \mb A)}  } f(ug)du .
\end{aligned}
\end{equation}
Therefore $H_{G, N, \emptyset, \bf 1}^{0, \, \on{cusp}} = C_c^{\on{cusp}}(G(F) \backslash G(\mb A) / K_{G, N} \Xi, \Ql )$.
\end{example}

\begin{rem}
When $I = \emptyset$, $W=\bf 1$ and $N=\emptyset$ (without level), for any $\mu \in \wh \Lambda_{G^{\mr{ad}}}^{+, \Q}$, $H_{M, N, I, W}^{' \, 0, \, \leq\mu}$ is included in the subspace of $C( U(\mb A) M(F) \backslash G(\mb A) / G(\mb O)  \Xi, \Ql )$ of functions supported on the components of $U(\mb A) M(F) \backslash G(\mb A) / G(\mb O)  \Xi$ indexed by a translated cone $\wh \Lambda_{Z_M / Z_G}^{\mu}$ in $\wh \Lambda_{Z_M / Z_G}^{\Q}$. 
The image of the constant term morphism is included in $H_{M, N, I, W}^{' \, 0} = \varinjlim_{\mu} H_{M, N, I, W}^{' \, 0, \, \leq\mu}$.
This space is already defined independently by Jonathan Wang in \cite{wang} Section 5.1 and is denoted by $\mc C_{P, \, -}$ in $loc.cit$.
\end{rem}

\quad
  
\section{Contractibility of deep enough horospheres}     \label{section-contractibility}

In this section, let $P$ be a parabolic subgroup of $G$ and $M$ its Levi quotient. 
The goal is to prove Proposition \ref{prop-TC-isom-mu-grand}, which will be a consequence of Theorem \ref{thm-i-contract-pour-mu-grand} and Theorem \ref{thm-pi-contract-pour-mu-grand}.

\subsection{More on Harder-Narasimhan stratification}    \label{subsection-more-on-HN}

To state Theorem \ref{thm-i-contract-pour-mu-grand} and Theorem \ref{thm-pi-contract-pour-mu-grand}, we need to introduce some locally closed substacks of $\Cht_{G, N, I, W}$.

\begin{defi}   \label{def-S-M-mu}
Let $\mu \in \wh{\Lambda}_{G^{\mr{ad}}}^{+, \Q}$. We define a set
$$
\begin{aligned}
S_M(\mu) : & =\{ \lambda \in \wh{\Lambda}_{G^{\mr{ad}}}^{+, \Q} \; | \; \lambda \leq^{\ov M} \mu \} \\
& = \{ \lambda \in \wh{\Lambda}_{G^{\mr{ad}}}^{+, \Q} \; | \; \lambda \leq \mu \} \bigcap \{  \lambda \in \wh{\Lambda}_{G^{\mr{ad}}}^{+, \Q} \; | \;  \on{pr}_P^{ad}(\lambda) = \on{pr}_P^{ad}(\mu) \} ,
\end{aligned}
$$ 
where the second equality follows from \ref{subsection-pr-p-ad-gamma-geq-0} (taking into account Notation \ref{notation-leq-means-leq-G-ad}).
The set $S_M(\mu)$ is bounded.
\end{defi}

\begin{rem}
The set $S_M(\mu)$ is the same as the one (modulo $\wh{\Lambda}_{Z_G}^{\Q}$) used in \cite{DG15} Sections 8 and 9 .
\end{rem}

\begin{defi}
$$\Bun_G^{=\mu}:= \underset{ \lambda \in \wh{\Lambda}_{G}^{+, \Q}, \;  \Upsilon_G(\lambda) =\mu} \bigcup \Bun_G^{(\lambda)}, \quad \Bun_M^{=\mu}:= \underset{ \lambda \in \wh{\Lambda}_{M}^{+, \Q}, \;    \Upsilon_G(\lambda) =\mu} \bigcup \Bun_M^{(\lambda)},$$
where $\Bun_G^{(\lambda)}$ (resp. $\Bun_M^{(\lambda)}$) is defined in Definition \ref{def-Bun-G-equal-mu} (resp. Definition \ref{def-Bun-M-equal-lambda}).
\end{defi}

\begin{defi}
$$\Bun_G^{S_M(\mu)}:= \underset{\lambda \in S_M(\mu)} \bigcup \Bun_G^{=\lambda}; \quad \Bun_M^{ S_M(\mu)}:= \underset{\lambda \in S_M(\mu)} \bigcup \Bun_M^{=\lambda}.$$
\end{defi}

\sssec{}   \label{subsection-S-M-mu-close-in-leq-mu}
If $\lambda \in S_M(\mu)$, $\lambda' \in \wh{\Lambda}_{G^{\mr{ad}}}^{+, \Q} $ and $\lambda \leq \lambda' \leq \mu$, then $\on{pr}_P^{ad}(\lambda) = \on{pr}_P^{ad}(\lambda') = \on{pr}_P^{ad}(\mu)$. This implies that $\lambda' \in S_M(\mu)$. Using \cite{DG15} Corollary 7.4.11, we deduce that: 
\begin{lem}   \label{lem-Bun-G-S-M-mu-locally-closed}
The substack $\Bun_G^{S_M(\mu)}$ is closed in $\Bun_G^{\leq \mu}$.
\cqfd
\end{lem}

\sssec{}   
We deduce from the definition of $S_M(\mu)$ and \ref{subsection-Bun-M-leq-ad-mu-nu-union-of-Bun-M-lambda} that
\begin{equation}   \label{equation-Cht-M-S-M-mu-Cht-M-leq-mu-nu}
 \Bun_M^{S_M(\mu)} = \Bun_M^{\leq\mu, \; \on{pr}_P^{ad}(\mu) } .
\end{equation}

Recall that $\Bun_{M}^{\leq \mu}$ is open in $\Bun_M$ (see Lemma \ref{lem-Bun-M-leq-ad-mu-open}) and $\Bun_M^{\on{pr}_P^{ad}(\mu)}$ is open and closed in $\Bun_M$ (see \ref{subsection-def-Bun-M-nu}). 
\begin{lem} \label{lem-Cht-M-S-M-mu-open-closed}  (\cite{DG15} Corollary 7.4.11, Lemma 8.2.6)
The substack $\Bun_M^{S_M(\mu)}$ is open and closed in $\Bun_M^{\leq \mu}$, and is open in $\Bun_M^{\on{pr}_P^{ad}(\mu)}$ and in $\Bun_M$.
\cqfd
\end{lem}


We define $\Bun_P^{S_M(\mu)}:=\Bun_P^{\leq \mu} \cap \; \pi^{-1} (\Bun_M^{S_M(\mu)})$. By Lemma \ref{lem-Cht-M-S-M-mu-open-closed}, it is open and closed in $\Bun_P^{\leq \mu}$, and is open in $\Bun_P$. So it is reduced. 
\begin{lem}    \label{lem-Bun-G-P-M-S-M-mu}  
Morphisms (\ref{diagram-Bun-G-Bun-P-Bun-M-leq-ad-mu}) induce morphisms
\begin{equation}   \label{equation-Bun-G-P-M-S-M-mu}
\Bun_G^{S_M(\mu)} \leftarrow \Bun_P^{S_M(\mu)} \rightarrow \Bun_M^{S_M(\mu)}.
\end{equation}
\end{lem}
\dem
We need to verify that the image of $\Bun_P^{S_M(\mu)} \rightarrow \Bun_G^{\leq \mu}$ is in the closed substack $\Bun_G^{S_M(\mu)}$.
Since $\Bun_P^{S_M(\mu)} $ is reduced, it is enough to consider geometric points. Let $\mc P \in \Bun_P^{S_M(\mu)}$ be a geometric point. Let $\mc M$ be its image in $\Bun_M^{S_M(\mu)}$. By definition of $\Bun_M^{S_M(\mu)}$, there exists $\lambda \in S_M(\mu)$ such that $\mc M \in \Bun_M^{=\lambda}$. 

Let $\mc G$ be the image of $\mc P$ in $\Bun_G^{\leq \mu}$. By  \ref{subsection-Bun-G-leq-ad-mu-union-strate}, there exists $\lambda' \leq \mu$ such that $\mc G \in \Bun_G^{=\lambda'}$. 
Taking into account that $\Bun_G^{=\lambda'} \subset \Bun_G^{\leq \lambda'}$,
by Lemma \ref{lem-Bun-P-leq-mu-dans-Bun-M-leq-mu}, we deduce that $\mc M \in \Bun_M^{\leq \lambda'}$. Hence $\lambda \leq \lambda'$.
By \ref{subsection-S-M-mu-close-in-leq-mu}, this implies that $\lambda' \in S_M(\mu)$. Thus $\mc G \in \Bun_G^{S_M(\mu)}$.
\cqfd

\begin{defi}    \label{def-Cht-G=mu-Cht-M=mu}
We define $\Cht_{G, N, I, W}^{= \mu}$ (resp. $\Cht_{G, N, I, W}^{S_M(\mu)}$) as the inverse image of $\Bun_{G}^{= \, \mu}$ (resp. $\Bun_G^{S_M(\mu)}$) by the morphism $$\Cht_{G, N, I, W} \rightarrow \Bun_{G}, \quad \left( (x_i)_{i \in I}, (\mc G, \psi) \xrightarrow{\phi}  (\ta \mc G, \ta \psi) \right) \mapsto \mc G.$$ 
Similarly, we define $\Cht_{P, N, I, W}^{ S_M(\mu)}$, $\Cht_{M, N, I, W}^{= \mu}$ and $\Cht_{M, N, I, W}^{ S_M(\mu)}$.
\end{defi}

\sssec{}   \label{subsection-Cht-M-S-M-mu-open-in-Cht-M}
We deduce from Lemma \ref{lem-Bun-G-S-M-mu-locally-closed} that $\Cht_{G, N, I, W}^{ S_M(\mu)}$ is closed in $\Cht_{G, N, I, W}^{\leq \mu}$. We deduce from Lemma \ref{lem-Cht-M-S-M-mu-open-closed} that $\Cht_{M, N, I, W}^{S_M(\mu)}$ is open and closed in $\Cht_{M, N, I, W}^{\leq \mu}$, and is open in $\Cht_{M, N, I, W}^{\on{pr}_P^{ad}(\mu)}$ and in $\Cht_{M, N, I, W}$.

\sssec{}
The commutativity of diagram (\ref{equation-Cht-G-P-M-Bun-G-P-M}) and Lemma \ref{lem-Bun-G-P-M-S-M-mu} imply that $\Cht_{P, N, I, W}^{S_M(\mu)} = \Cht_{P, N, I, W}^{\leq \mu}  \cap \;  \pi^{-1} ( \Cht_{M, N, I, W}^{S_M(\mu)} ). $
Morphisms (\ref{equation-Cht-G-P-M-leq-mu-Xi}) induce morphisms:
$$\Cht_{G, N, I, W}^{S_M(\mu)} \xleftarrow{i^{S_M(\mu)}} \Cht_{P, N, I, W}^{S_M(\mu)} \xrightarrow{\pi^{S_M(\mu)}} \Cht_{M, N, I, W}^{ S_M(\mu)} . $$


\sssec{}
As in Definition \ref{def-Cht-P-M-'}, we define 
$$\Cht_{P, N, I, W}^{' \, S_M(\mu)}  := \Cht_{P, N, I, W}^{S_M(\mu)} \overset{P(\mc O_N)} \times G(\mc O_N) , \quad
\Cht_{M, N, I, W}^{' \, S_M(\mu)}  := \Cht_{M, N, I, W}^{S_M(\mu)}  \overset{P(\mc O_N)} \times G(\mc O_N) .$$
Morphisms (\ref{equation-Cht-G-P-M-prime}) induce morphisms
\begin{equation}
\Cht_{G, N, I, W}^{S_M(\mu)} \xleftarrow{i^{' \, S_M(\mu)}} \Cht_{P, N, I, W}^{' \, S_M(\mu)} \xrightarrow{\pi^{' \, S_M(\mu)}} \Cht_{M, N, I, W}^{' \, S_M(\mu)}
\end{equation}

\subsection{Geometric statements}    \label{subsection-geo-statements}

Firstly consider the morphism $i^{' \, S_M(\mu)}$. 
\begin{thm}  (\cite{var} Theorem 2.25 and Proposition 5.7, \cite{DG15} Proposition 9.2.2)  \label{thm-i-contract-pour-mu-grand}
There exists a constant $C'(G, X, N, W)$, such that if $\mu \in \wh{\Lambda}_{G^{\mr{ad}}}^{+, \Q}$ and $\langle  \mu, \alpha \rangle >  C'(G, X, N, W)$ for all $\alpha \in \Gamma_G - \Gamma_M$, then
the morphism $i^{' \, S_M(\mu)}$
is a schematic finite universal homeomorphism.
\end{thm}

\dem
(1) Schematic and finite follows from Proposition 5.7 of \cite{var} (recalled in Proposition \ref{prop-Varshavsky-proper} and Remark \ref{rem-i-proper-schematic})

(2) Surjectivity is implied by Theorem 2.25 of \cite{var}. 

(3) Universally injectivity is implied by the fact that $\Bun_{P}^{S_M(\mu)} \rightarrow \Bun_{G}^{S_M(\mu)}$ is an isomorphism for $\mu$ satisfying the assumption of Theorem \ref{thm-i-contract-pour-mu-grand} (see \cite{DG15} Proposition 9.2.2) and the well-known fact that $\Gr_{P, I, W} \rightarrow \Gr_{G, I, W}$ is bijective. 
(More concretely, it is enough to prove that for any algebraically closed field $k$ containing $\Fqbar$, the map $\Cht_{P, I, W}^{S_M(\mu)}(k) \rightarrow \Cht_{G, I, W}^{S_M(\mu)}(k)$ is injective. Let $((x_i), \mc G \xrightarrow{\phi_G} \ta\mc G) \in \Cht_{G, I, W}^{S_M(\mu)}(k)$. By (3), there exists $((x_i), \mc P \xrightarrow{\phi_P} \ta\mc P) \in \Cht_{P, I, W}^{S_M(\mu)}(k)$ such that $\mc P \overset{P}{\times} G \simeq \mc G$ and $\phi_P \overset{P}{\times} G \simeq \phi_G$.
Since $\Bun_{P}^{S_M(\mu)}(k) \rightarrow \Bun_{G}^{S_M(\mu)}(k)$ is injective, $\mc P$ is unique. Choosing a trivialisation of $\mc P$ over $\Gamma_{\sum \infty x_i}$, we deduce from the injectivity of $\Gr_{P, I, W}(k) \rightarrow \Gr_{G, I, W}(k)$ that $\phi_P$ is unique.)
\cqfd


\sssec{}
Now we consider the morphism $\pi^{' \, S_M(\mu)}$. For all $d$ large enough, similar to diagram (\ref{diagram-cht-P-cht-M-pi-d}), we have a commutative diagram
\begin{equation}    \label{equation-Cht-P-to-Cht-M-tilte-prime}
\xymatrixrowsep{1.5pc}
\xymatrixcolsep{3pc}
\xymatrix{
\Cht_{P, N, I, W}^{' \, S_M(\mu)}  \ar[rdd]   \ar[rd]|-{\pi_{d}^{' \, S_M(\mu)}}   \ar[rrd]^{\pi^{' \, S_M(\mu)}}
& 
&  \\
& \wt{\Cht}_{M, N, I, W}^{' \, S_M(\mu)}     \ar[d]  \ar[r]
& \Cht_{M, N, I, W}^{' \, S_M(\mu)}  \ar[d] \\
&  [ P_{I, d} \backslash  \Gr_{P, I, W}]  \ar[r]
&  [  M_{I, d} \backslash  \Gr_{M, I, W}]
}
\end{equation}
where $\widetilde{\Cht}_{M, N, I, W}^{' \, S_M(\mu)} $ is the fiber product, which depends on $d$.
By \ref{subsection-Cht-M-S-M-mu-open-in-Cht-M}, $\Cht_{M, N, I, W}^{' \, S_M(\mu)} $ is open in $\Cht_{M, N, I, W}^{'}$ and $\Cht_{P, N, I, W}^{' \, S_M(\mu)}$ is open in $\Cht_{P, N, I, W}^{'}$. By Lemma \ref{lem-pi-d-est-lisse}, the morphism $\pi_{d}^{' \, S_M(\mu)}$ is smooth of relative dimension $\dim_{X^I} U_{I, d}$.

\quad

We now introduce a notation of unipotent group scheme (which should rather be called "elementary unipotent group scheme").

\begin{defi} \label{def-sch-en-gp-unipotent}

(a) Let $H$ be a group scheme of finite dimension over a scheme $S$. 
We say that $H$ is a unipotent group scheme if 
$H$ admits a filtration $H=H^{(0)} \supset H^{(1)} \supset \dots \supset H^{(m)} \supset H^{(m+1)}=0$ such that for every $j$, the quotient $H^{(j)} / H^{(j+1)}$ is an additive group scheme (i.e. isomorphic to $\mb G_{a, S}^n$ for some $n$ locally for the étale topology) over $S$.

(b) A morphism of algebraic stacks $f: \mc X \rightarrow \mc Y$ is called unipotent if for any scheme $S$ and any morphism $S \rightarrow \mc Y$, the fiber product $S \underset{\mc Y}{\times} \mc X$ is locally for the smooth topology on $S$ isomorphic to a quotient stack $[H_1 / H_2]$, where $H_1$ and $H_2$ are unipotent group schemes over $S$ and $H_2$ acts on $H_1$ as a group scheme over S acting on a scheme over S. 
\end{defi}

\begin{thm} \label{thm-pi-contract-pour-mu-grand}
There exists a constant $C(G, X, N, I, d)$, such that if $\mu \in \wh{\Lambda}_{G^{\mr{ad}}}^{+, \Q}$ and $\langle  \mu, \alpha \rangle > C(G, X, N, I, d)$ for all $\alpha \in \Gamma_G - \Gamma_M$, then the morphism $\pi_{d}^{' \, S_M(\mu)}$ is unipotent in the sense of Definition \ref{def-sch-en-gp-unipotent}.

\end{thm}

The proof will be given in Sections \ref{subsection-proof-pi-step-1}-\ref{subsection-proof-pi-step-3}.

\begin{rem}
Theorem \ref{thm-pi-contract-pour-mu-grand} will be used to prove Proposition \ref{prop-TC-isom-mu-grand}, where only the statement for the geometric fibers of $\pi_{d}^{' \, S_M(\mu)}$ is needed. Since the proof is the same for a geometric fiber or a fiber over a general base, we prove it over a general base.
\end{rem}


\subsection{Proof of Theorem \ref{thm-pi-contract-pour-mu-grand}: step 1}    \label{subsection-proof-pi-step-1}

\sssec{}
We have a similar diagram as (\ref{equation-Cht-P-to-Cht-M-tilte-prime}) without index $'$. 
The morphism $\pi_{d}^{S_M(\mu)}: \Cht_{P, N, I, W}^{S_M(\mu)}  \rightarrow \wt{\Cht}_{M, N, I, W}^{S_M(\mu)}$ is $P(\mc O_N)$-equivariant and the morphism $\pi_{d}^{' \, S_M(\mu)}$ is induced by $\pi_{d}^{S_M(\mu)} \overset{P(\mc O_N)}{\times} G(\mc O_N)$.
So to prove Theorem \ref{thm-pi-contract-pour-mu-grand}, it is enough to prove the statement for $\pi_{d}^{S_M(\mu)}$ instead of $\pi_{d}^{' \, S_M(\mu)}$.

The problem is local for the smooth topology. So it is enough to prove the statement for the base change by $\Gr_{P, I, W} \rightarrow [P_{I, d} \backslash \Gr_{P, I, W}]$:
$$\pi_d^{\times, \, S_M(\mu)}: \Cht_{P, N, I, W}^{S_M(\mu)} \underset{ [P_{I, d} \backslash \Gr_{P, I, W}] }{   \times  }\Gr_{P, I, W}   \rightarrow \widetilde{\Cht}_{M, N, I, W}^{S_M(\mu)} \underset{ [P_{I, d} \backslash \Gr_{P, I, W}] }{ \times }  \Gr_{P, I, W} .$$


\sssec{}   \label{subsection-cube-Cht}
Note that $ \widetilde{\Cht}_{M, N, I, W}^{S_M(\mu)} \underset{ [P_{I, d} \backslash \Gr_{P, I, W}] }{ \times }  \Gr_{P, I, W}  \simeq \Cht_{M, N, I, W}^{S_M(\mu)} \underset{ [M_{I, d} \backslash \Gr_{M, I, W}] }{   \times  }\Gr_{P, I, W} $.
We have the following commutative diagram, where the front and back Cartesian squares are defined in the proof of Proposition 2.8 in \cite{vincent} (replace $G$ by $P$ and $M$, respectively). We have already used these Cartesian squares in (\ref{equation-Cht-P-cartesian-square}) and (\ref{equation-Cht-M-cartesian-square}).

{  \resizebox{15cm}{!}{ $$ \!\!\!\!\!\!\!\!\!\!\!\!\!\!\!\!\!\!\!
\xymatrixrowsep{1pc}
\xymatrixcolsep{0.5pc}
\xymatrix{
\underset{ [P_{I, d} \backslash \Gr_{P, I, W}] } { \Cht_{P, N, I, W}^{S_M(\mu)}  \times \Gr_{P, I, W} } \ar[rr] \ar[dd]  \ar[rd]^{ \pi_d^{\times, \, S_M(\mu)}    }
& &  \underset{(X \sm N)^I} { \Bun_{P, N, I, d}^{S_M(\mu)}  \times \Gr_{P, I, W}  } \ar@{.>}[dd]^{(b_1^P, b_2^P)}   \ar[rd]^{ ( \pi^{Bun}, \; \Id ) }
& \\
&  \underset{ [P_{I, d} \backslash \Gr_{P, I, W}] }{ \widetilde{\Cht}_{M, N, I, W}^{S_M(\mu)}  \times \Gr_{P, I, W} }  \ar[rr]  \ar[dd]
&  & \underset{(X \sm N)^I} { \Bun_{M, N, I, d}^{S_M(\mu)}   \times \Gr_{P, I, W}  }   \ar[dd]^{(b_1^M, b_2^M)}  \\
\Bun_{P, N}^{S_M(\mu)}  \ar@{.>}[rr]^{(\on{Id}, \on{Id}) \quad \quad \quad \quad \quad}   \ar[rd]^{\pi^{Bun}}
&  & \Bun_{P, N}^{S_M(\mu)}  \times \Bun_{P, N}  \ar@{.>}[rd]^{(\pi^{Bun}, \pi^{Bun} )}
& \\
& \Bun_{M, N}^{S_M(\mu)}   \ar[rr]^{(\Id, \Id)} 
&  &  \Bun_{M, N}^{S_M(\mu)}  \times \Bun_{M, N}
}
$$
} }

\sssec{}
Now let $S$ be an affine scheme over $\Fqbar$ and let $$\big( (x_i), (\mc M, \psi) \xrightarrow{\phi} (\ta \mc M, \ta \psi), s \big): S \rightarrow \widetilde{\Cht}_{M, N, I, W}^{S_M(\mu)} \underset{ [P_{I, d} \backslash \Gr_{P, I, W}] }\times \Gr_{P, I, W}$$ be an $S$-point. Consider 
$$ S \rightarrow \widetilde{\Cht}_{M, N, I, W}^{S_M(\mu)} \underset{ [P_{I, d} \backslash \Gr_{P, I, W}] }\times \Gr_{P, I, W} \rightarrow \underset{(X \sm N)^I} { \Bun_{M, N, I, d}^{S_M(\mu)}   \times \Gr_{P, I, W}  }  $$
and
$$ S \rightarrow \widetilde{\Cht}_{M, N, I, W}^{S_M(\mu)} \underset{ [P_{I, d} \backslash \Gr_{P, I, W}] }\times \Gr_{P, I, W} \rightarrow \Bun_{M, N}^{S_M(\mu)}  $$

We define 
$\ms Z$, $\ms Y_{N, d}$ and $\ms Y_{N}$
to be the following fiber products:
{  \resizebox{15cm}{!}{ $$ \!\!\!\!\!\!\!\!\!
\xymatrixrowsep{2pc}
\xymatrixcolsep{2pc}
\xymatrix{
\ms Z \ar[r] \ar[d]
&\Cht_{P, N, I, W}^{S_M(\mu)} \underset{ [P_{I, d} \backslash \Gr_{P, I, W}] }\times \Gr_{P, I, W}  \ar[d]^{\pi_d^{\times, \, S_M(\mu)}} 
& \ms Y_{N, d}  \ar[r] \ar[d]
&\Bun_{P, N, I, d}^{S_M(\mu)} \underset{(X \sm N)^I}\times \Gr_{P, I, W} \ar[d]^{(\pi^{Bun}, \, \Id)}
&\ms Y_N  \ar[r] \ar[d]
&\Bun_{P, N}^{S_M(\mu)} \ar[d]^{\pi^{Bun}} \\
S \ar[r]
& \widetilde{\Cht}_{M, N, I, W}^{S_M(\mu)} \underset{ [P_{I, d} \backslash \Gr_{P, I, W}] }\times \Gr_{P, I, W} 
& S \ar[r]
&\Bun_{M, N, I, d}^{S_M(\mu)} \underset{(X \sm N)^I}\times \Gr_{P, I, W}  
&S \ar[r]
&\Bun_{M, N}^{S_M(\mu)}
}
$$
} }

Applying Lemma \ref{lem-cube} to the diagram in \ref{subsection-cube-Cht}, we deduce a Cartesian square 
\begin{equation}    \label{diagram-Cht-est-egalisateur-de-Bun-egal-mu}
\xymatrixrowsep{2pc}
\xymatrixcolsep{3pc}
\xymatrix{
\ms Z \ar[r] \ar[d]
&\ms Y_{N, d}  \ar[d]^{(b_1, b_2)} \\
\ms Y_N \ar[r]^{(Id, Id) \quad \quad}
& \ms Y_N \underset{S}\times \ms Y_N .
}
\end{equation}
where $b_1$ (resp. $b_2$) is induced by $b_1^P$ (resp. $b_2^P$).


\begin{rem}
By the proof of Proposition 2.8 in \cite{vincent}, $b_1^P$ (resp. $b_1^M$) is the forgetful morphism of the level structure on $I$ (thus smooth) and $b_2^P$ (resp. $b_2^M$) is the composition of the Frobenius morphism with some other morphism. We deduce that $b_1$ is smooth and $b_2$ has zero differential. Moreover, the morphism $\Bun_{P, N} \rightarrow \Bun_{M, N}$ is smooth, thus $\ms Y_N$ is smooth over $S$. Similarly $\ms Y_{N, d}$ is smooth over $S$. We deduce that $\ms Z$ is smooth over $S$. Note that the same argument without $S_M(\mu)$ would give another proof of Lemma \ref{lem-pi-d-est-lisse}.
\end{rem}

\subsection{Proof of Theorem \ref{thm-pi-contract-pour-mu-grand}: step 2}     \label{subsection-proof-pi-step-2}

The goal of this subsection is to describe $\ms Y_N$ and $\ms Y_{N, d}$.

(1) Firstly we describe the fiber of $\Bun_P^{S_M(\mu)} \rightarrow \Bun_M^{S_M(\mu)}$ in Proposition \ref{prop-fiber-Bun-P-Bun-M}.

\sssec{}   \label{subsection-fix-T-then-M-includ-in-P}
We fix a maximal torus $T \subset B$. This allows us to view the Levi quotient $M$ of a standard parabolic subgroup $P$ as a subgroup $M \subset P$ (the unique splitting that contains $T$).
Then $P = M \rtimes U$, where $M$ acts on $U$ via the embedding $M \subset P$ and $P$ acts on $U$ by the adjoint action.

\sssec{}
Let $S \rightarrow \Bun_{M}^{S_M(\mu)}$ be a morphism and $\mc M$ the corresponding $M$-bundle over $X \times S$.
We define the fiber space $U_{\mc M}: = (U \times \mc M) / M$. It is easy to check that it is a group scheme over $X \times S$ (see \cite{these} C.2 for more details).

\begin{defi} 
Let $S$ be an affine scheme over $\Fqbar$. Let $A$ be a sheaf of groups on $X \times S$. We denote by $\on{pr}_S: X \times S \rightarrow S$ the second projection. 

(a) We define $R^0(\on{pr}_S)_*A$ as the sheaf of groups on $S$:
$$(S' \rightarrow S)  \mapsto \Hom_{X \times S}(X \times S', A'),$$
where $A'$ is the inverse image of $A$ by $X \times S' \rightarrow X \times S$.

(b) (\cite{giraud} V.2.1) We define $R^1(\on{pr}_S)_*A$ as the sheaf of sets on $S$ associated to the presheaf:
$$(S' \rightarrow S)  \mapsto H^1(X \times S', A').$$
Indeed $R^1(\on{pr}_S)_*A$ is a sheaf of pointed sets with a canonical section which corresponds to the trivial $A$-torsor.
\end{defi}

\begin{prop}   \label{prop-fiber-Bun-P-Bun-M}
There exists a constant $C(G, X) \in \Q_{\geq 0}$, such that if $\langle  \mu, \alpha \rangle > C(G, X) \text{ for all }  \alpha \in \Gamma_G - \Gamma_M$, then $R^0(\on{pr}_S)_* U_{\mc M}$ is a unipotent group scheme over $S$ and the fiber of $\Bun_{P}^{S_M(\mu)} \rightarrow \Bun_{M}^{S_M(\mu)}$ over $S$ is the classifying stack $[S / R^0(\on{pr}_S)_* U_{\mc M} ]$. 
\end{prop}
\dem
We denote by $\ms Y$ the fiber of $\Bun_{P}^{S_M(\mu)} \rightarrow \Bun_{M}^{S_M(\mu)}$ over $S$. For any scheme $S' \rightarrow S$, 
the groupoid $\ms Y(S')$ classifies the 
$\restr{U_{\mc M}}{X \times S'}$-bundle over $X \times S'$ (see \cite{these} Lemme C.3.2 for more details).

By Lemma \ref{lem-R^0-U-unipotent-R^1-U-trivial} (b) below, all $U_{\mc M}$-bundles are trivial. Taking into account that $R^0(pr_S)_*  U_{\mc M}(S')$ is the group of automorphisms of the trivial $\restr{U_{\mc M}}{X \times S'}$-bundle on $X \times S'$ and Lemma \ref{lem-R^0-U-unipotent-R^1-U-trivial} (a), 
we deduce the proposition.
\cqfd

%

\begin{lem}   \label{lem-R^0-U-unipotent-R^1-U-trivial}
There exists a constant $C(G, X) \in \Q_{\geq 0}$, such that if $\langle  \mu, \alpha \rangle > C(G, X) \text{ for all }  \alpha \in \Gamma_G - \Gamma_M$, then

(a) the sheaf of groups $ R^0(pr_S)_* U_{\mc M} $ is a unipotent group scheme; 

(b) the sheaf of pointed sets $ R^1(pr_S)_* U_{\mc M} $ is trivial.
\end{lem}

\begin{rem}    \label{rem-U-comm-vs-non-comm}
If $U$ is commutative, then $U_{\mc M}$ is an additive group scheme over $X \times S$ (in the sense of Definition \ref{def-sch-en-gp-unipotent}).
Lemma \ref{lem-R^0-U-unipotent-R^1-U-trivial} (a) is automatic and (b) follows directly from \cite{DG15} Proposition 10.4.5. 

The difficulty is that in general, $U$ is not commutative. 
To prove Lemma \ref{lem-R^0-U-unipotent-R^1-U-trivial}, we will need to use a filtration of $U$ where the graded are commutative groups.
\end{rem}

\sssec{}
We have a canonical filtration of $U$ (see the proof of Proposition 11.1.4 (c) in \cite{DG15} for more details): 
\begin{equation}  \label{equation-filtration-U}
U=U^{(0)} \supset U^{(1)} \supset \dots \supset U^{(m)} \supset U^{(m+1)}=0 ,
\end{equation}
where $U^{(j)}$ is the subgroup generated by the root subgroups corresponding to the positive roots $\alpha$ of $G$, such that $$\sum_{\beta \in \Gamma_G - \Gamma_M} \on{coeff}_{\beta} (\alpha) \geq j+1.$$ (Here $\on{coeff}_{\beta} (\alpha)$ denotes the coefficient of $\alpha$ in simple root $\beta$.)
For each $j$, the subgroup $U^{(j+1)}$ of $U^{(j)}$ is normal and the quotient is equipped with an isomorphism $\vartheta^{(j)} : \mb G_a^{n_j} \isom U^{(j)} / U^{(j+1)} $ for some $n_j \in \N$.

\sssec{}
The filtration (\ref{equation-filtration-U}) induces for every $j \in \{1, \cdots, m+1\}$ an exact sequence of groups:
\begin{equation}  \label{equation-suite-exacte-U^j}
0 \rightarrow U^{(j-1)} / U^{(j)} \rightarrow U/ U^{(j)} \rightarrow U/ U^{(j-1)} \rightarrow 0 .
\end{equation}

For every $j$, the subgroup $U^{(j)}$ of $P$ is normal. Then $P$ acts on $U^{(j)}$ by the adjoint action and $M$ acts on $U^{(j)}$ via $M \hookrightarrow P$. We deduce that $M$ acts on $U^{(j)} / U^{(j+1)}$ and $U / U^{(j)}$.

We define the fiber spaces $(U^{(j)} / U^{(j+1)})_{\mc M} := \big( \mc M  \times  U^{(j)} / U^{(j+1)} \big) / M$, it is an additive group scheme over $X \times S$. We define the fiber space $ (U / U^{(j)} )_{\mc M} := \big( \mc M  \times U/ U^{(j)}  \big) / M $, it is a group scheme over $X \times S$.  
(see \cite{these} C.2 for more details). 

\begin{prop}  \label{prop-annulation-H-1-sans-niveau}
There exists a constant $C(G, X)$ such that for $\mu \in \wh{\Lambda}_{G^{\mr{ad}}}^{+, \Q}$, if $\langle  \mu, \alpha \rangle > C(G, X) \text{ for all }  \alpha \in \Gamma_G - \Gamma_M$,
then for any $\mc M \in \Bun_{M}^{S_M(\mu)}(S)$ and any $j$, the sheaf $R^1(pr_S)_* \big( (U^{(j)} / U^{(j+1)})_{\mc M}  \big) $ is trivial.
\cqfd
\end{prop}
\dem 
This is \cite{DG15} Proposition 10.4.5 (a). We take $C(G, X):=\on{max}_{i} \{c_i'\}$, where $c_i'$ are the constants in $loc.cit.$.
\cqfd

\begin{lem}  \label{lem-H^1-A-B-C}
Let $0 \rightarrow A \rightarrow B \rightarrow C \rightarrow 0 $ be an exact sequence of sheaves of groups on $X \times S$.

(a) If the sheaf of pointed sets $R^1(pr_S)_* A$ is trivial, then we have an exact sequence of sheaves of groups:
$$0 \rightarrow R^0(pr_S)_*A \rightarrow R^0(pr_S)_* B \rightarrow R^0(pr_S)_*C \rightarrow 0. $$

(b) If moreover the sheaf of pointed sets $R^1(pr_S)_* C$ is also trivial, then the sheaf of pointed sets $R^1(pr_S)_* B$ is trivial.
\end{lem}
\dem
By \cite{giraud} V Proposition 2.3, the exact sequence $0 \rightarrow A \rightarrow B \rightarrow C \rightarrow 0 $ induces an exact sequence of sheaves of pointed sets on $S$:
$$
\begin{aligned}
0 \rightarrow  R^0(pr_S)_* A \rightarrow R^0(pr_S)_* B \rightarrow R^0(pr_S)_* C \rightarrow  \\
R^1(pr_S)_*A \rightarrow R^1(pr_S)_* B \rightarrow R^1(pr_S)_*C .
\end{aligned}
$$
We deduce the lemma.
\cqfd

\noindent {\bf Proof of Lemma \ref{lem-R^0-U-unipotent-R^1-U-trivial}.} 
For each $j$, the exact sequence (\ref{equation-suite-exacte-U^j}) induces an exact sequence of group schemes over $X \times S$:
\begin{equation}  \label{equation-suite-exacte-U^j-F-M}
0 \rightarrow (U^{(j-1)} / U^{(j)})_{\mc M} \rightarrow (U/ U^{(j)})_{\mc M} \rightarrow (U/ U^{(j-1)})_{\mc M} \rightarrow 0 .
\end{equation}
We apply Lemma \ref{lem-H^1-A-B-C} to (\ref{equation-suite-exacte-U^j-F-M}) successively for $j=1, j=2, \cdots ,$ until $j=m+1$. Taking into account that $R^1(pr_S)_* \big( (U^{(j)} / U^{(j+1)})_{\mc M}  \big) $ is trivial (by Proposition \ref{prop-annulation-H-1-sans-niveau}) and that $R^0(pr_S)_*(U^{(j)}/ U^{(j+1)})_{\mc M}$ is additive in the sense of Definition \ref{def-sch-en-gp-unipotent} (because $U^{(j)} / U^{(j+1)} \isom \mb G_a^{n_j} $), we deduce Lemma \ref{lem-R^0-U-unipotent-R^1-U-trivial}.
\cqfd

\quad

(2) Now we add level structure on $N \times S + \Gamma_{\sum d x_i}$ to the argument in (1), i.e. we describe the fiber of $\Bun_{P, N, I, d}^{S_M(\mu)} \rightarrow \Bun_{M, N, I, d}^{S_M(\mu)}$ in Proposition \ref{prop-Bun-P-Bun-M-avec-N+d}.

\sssec{}  \label{subsection-Ker-faiseau-en-groupes}
%
%
Let $V$ be a group scheme on $X \times S$. For any divisor $i_D: D \hookrightarrow X \times S$, we denote by $\restr{V}{D}$ the fiber product $D \underset{X \times S} \times V$. We denote by $\mc V$ and $\restr{\mc V}{D}$ the associated sheaves of groups.
We define the sheaf of groups $\mc Ker_{V, D}$ on $X \times S$ as the kernel of the morphism $\mc V \rightarrow (i_D)_* ( \restr{\mc V}{D} ). $
If $V$ is smooth, the morphism $\mc V \rightarrow (i_D)_* ( \restr{\mc V}{D} )$ is surjective. 

\sssec{}
Let $S$ be an affine scheme over $\Fqbar$. Let $\big( (x_i)_{i \in I}, \mc M, \psi_M \big)$ be an $S$-point of $\Bun_{M, N, I, d}^{S_M(\mu)}$. Let $D := N \times S + \Gamma_{\sum d x_i} .$ 
Applying \ref{subsection-Ker-faiseau-en-groupes} to the group scheme $U_{\mc M}$ on $X \times S$, we obtain an exact sequence of sheaves of groups:
\begin{equation}   \label{equation-Ker-U-M-D-exact-seq}
0 \rightarrow \mc Ker_{U_{\mc M}, D} \rightarrow \mc U_{\mc M} \rightarrow (i_D)_* \restr{ \mc U_{\mc M} }{D} \rightarrow 0.
\end{equation}

\begin{prop}   \label{prop-Bun-P-Bun-M-avec-N+d}
There exists a constant $C(G, X, N, I, d) \in \Q_{\geq 0}$, such that if $\langle  \mu, \alpha \rangle > C(G, X, N, I, d) \text{ for all }  \alpha \in \Gamma_G - \Gamma_M$, then $R^0(pr_S)_* \mc Ker_{U_{\mc M}, D}$ is a unipotent group scheme over $S$ and the fiber of $\Bun_{P, N, I, d}^{S_M(\mu)} \rightarrow \Bun_{M, N, I, d}^{S_M(\mu)}$ over $S$ is the classifying stack $[S / R^0(pr_S)_* \mc Ker_{U_{\mc M}, D} ]$.
\end{prop}

\dem
We recall that $\ms Y_{N, d}$ denotes the fiber of $\Bun_{P, N, I, d}^{S_M(\mu)} \rightarrow \Bun_{M, N, I, d}^{S_M(\mu)}$ over $S$. For any scheme $S' \rightarrow S$, the groupoid $\ms Y_{N, d}(S')$ classifies the data of $(\mc F, \beta)$, where $\mc F$ is a $\mc U_{\mc M}$-bundle on $X \times S'$ and $\beta$ is an isomorphism of $\mc U_{\mc M}$-bundles $\restr{  \mc F  }{ D }  \isom  \restr{ \mc U_{\mc M} } {D} $. 
By (\ref{equation-Ker-U-M-D-exact-seq}), this groupoid is equivalent to the groupoid of $ \mc Ker_{U_{\mc M}, D}$-bundles on $X \times S'$.

Similarly to the case without level, Proposition \ref{prop-Bun-P-Bun-M-avec-N+d} follows from Lemma \ref{lem-R-0-Ker-R-1-Ker} below.
\cqfd

\begin{lem}   \label{lem-R-0-Ker-R-1-Ker}
There exists a constant $C(G, X, N, I, d) \in \Q_{\geq 0}$, such that if $\langle  \mu, \alpha \rangle > C(G, X, N, I, d) \text{ for all }  \alpha \in \Gamma_G - \Gamma_M$, then 

(a) the sheaf of groups $ R^0(pr_S)_* \mc Ker_{U_{\mc M}, D} $ is a unipotent group scheme;

(b) the sheaf of pointed sets $ R^1(pr_S)_* \mc Ker_{U_{\mc M}, D} $ is trivial.
\end{lem}
\dem
The proof is the same as Lemma \ref{lem-R^0-U-unipotent-R^1-U-trivial}, except that we replace $(U^{(j-1)} / U^{(j)})_{\mc M}$ by $\mc Ker_{ (U^{(j-1)} / U^{(j)})_{\mc M} , D }$, and that we use Lemma \ref{lem-annulation-H-1-niveau-N+d} below instead of Proposition \ref{prop-annulation-H-1-sans-niveau}.
\cqfd

\begin{lem}     \label{lem-annulation-H-1-niveau-N+d}
There exists a constant $C(G, X, N, I, d) \in \Q_{\geq 0}$ such that for $\mu \in \wh{\Lambda}_{G^{\mr{ad}}}^{+, \Q}$,
if $\langle  \mu, \alpha \rangle > C(G, X, N, I, d) , \text{ for all }  \alpha \in \Gamma_G - \Gamma_M$,
then for any $((x_i), \mc M, \psi) \in \Bun_{M, N, I, d}^{S_M(\mu)}(S)$ and any $j$, the sheaf $R^1(pr_S)_* \big( (U^{(j)} / U^{(j+1)})_{\mc M} (- N \times S - \Gamma_{\sum d x_i}   )  \big) $ is trivial.
\end{lem}
\dem
Let $C(G, X, N, I, d):=C(G, X)+\deg N + |I| \cdot d$, where $C(G, X)$ is the constant in Proposition \ref{prop-annulation-H-1-sans-niveau}. 
We repeat the argument in \cite{DG15} Proposition 10.4.5, except that in $loc.cit.$ Remark 10.3.5 we replace the reductive group $\wt G$ by $\wt G \times \mb G_m$ and the $\wt G$-bundle $\mc F_{\wt G}$ by the $\wt G \times \mb G_m$-bundle $\mc F_{\wt G} \times \mc O( -N \times S - \Gamma_{\sum d x_i} )$.

\cqfd

\subsection{Proof of Theorem \ref{thm-pi-contract-pour-mu-grand}: step 3}      \label{subsection-proof-pi-step-3}

\sssec{}
Let $S$ be a scheme over $\Fq$. Let $H_S$ and $H_S'$ be two group schemes over $S$. Let $f: H_S' \rightarrow H_S$ be a morphism of group schemes over $S$. We denote by $[S/H_S']$ the classifying stack of $H_S'$ on $S$. Similarly for $[S/H_S]$.
Then $f$ induces a morphism of stacks: $\ov f: [S/ H_S'] \rightarrow [S/H_S].$

\begin{lem} \label{lem-egalisateur-est-unipotent}
Let $f, g: H_S' \rightarrow H_S$ be two morphisms of connected group schemes. Let $[H_S / H_S'] $ be the quotient stack where $H_S'$ acts on $H_S$ by $h' \cdot h = f(h') h g(h')^{-1}$. Then the following diagram is Cartesian:
\begin{equation}
\xymatrix{
[H_S / H_S'] \ar[r]  \ar[d]
& [S / H_S']    \ar[d]^{(\ov f, \ov g)} \\
[S / H_S]  \ar[r]^{(id, id) \quad \quad}
& [S / H_S] \underset{S}\times [S / H_S]  
}
\end{equation}
where the morphism $[H_S / H_S'] \rightarrow [S / H_S] $ is induced by $H_S \rightarrow S$ and $H_S' \xrightarrow{g} H_S$.
\end{lem}


\dem
The fiber product is $[H_S \times_S H_S / H_S' \times_S H_S]$, where $H_S'$ acts on $H_S \times_S H_S$ by $(f, g)$ (from the left) and $H_S$ acts on $H_S \times_S H_S$ by diagonal action (from the right). The morphism $ \alpha:  [H_S \times_S H_S / H_S' \times_S H_S] \rightarrow [S / H_S]$ (resp. $\beta: [H_S \times_S H_S / H_S' \times_S H_S] \rightarrow  [S / H_S']$) is given by $H_S \times_S H_S \rightarrow S$ and the second projection $H_S' \times_S H_S \rightarrow H_S$ (resp. the first projection $H_S' \times_S H_S \rightarrow H_S'$).

The morphism of group schemes over $S$: 
\begin{equation}    \label{equation-H-times-H-isom-H-times-H}
H_S \times_S H_S \rightarrow H_S \times_S H_S, \quad (x, y) \mapsto (xy^{-1}, y)
\end{equation}
is an isomorphism. Moreover, it is $H_S' \times_S H_S$-equivariant for the action of $H_S' \times_S H_S$ on the LHS as above and the action of $H_S' \times_S H_S$ on the RHS given by $(h', h)(z, t) = ( f(h') z g(h')^{-1} , g(h') t h^{-1}  )$. 
The isomorphism (\ref{equation-H-times-H-isom-H-times-H}) induces an isomorphism of quotient stacks 
\begin{equation}    \label{equation-H-times-H-quotient-isom-H-quotient}
[H_S \times_S H_S / H_S' \times_S H_S] \underset{(\ref{equation-H-times-H-isom-H-times-H})}{\isom} [H_S \times_S H_S / H_S' \times_S H_S] \simeq [H_S / H_S'],
\end{equation}
where $H_S'$ acts on $H_S$ by $h' \cdot x = f(h') x g(h')^{-1}$. The morphism $[H_S / H_S'] \rightarrow [S / H_S]$ is the composition of the inverse of (\ref{equation-H-times-H-quotient-isom-H-quotient}) and $\alpha$.
\cqfd

\begin{lem}    \label{lem-BH1-BH2-induce-H1-H2}
Let $S$ be an affine scheme. Let $H_1$ and $H_2$ be two unipotent group schemes over $S$. Let $\varphi:  [S / H_1] \rightarrow [S/ H_2]$ be a morphism of stacks. Then there exists $f: H_1 \rightarrow H_2$ a morphism of group schemes over $S$ such that $\varphi = \ov{f}$.
\end{lem}
\dem
Since $S$ is affine and $H_2$ is unipotent, all $H_2$-torsors on $S$ are trivial. The morphism $\varphi$ is given by a $H_2$-torsor $\ms H$ on $S$ which is $H_1$-equivariant. We trivialize $\ms H$ as a $H_2$-torsor. Then the action of $H_1$ on $\ms H$ gives the morphism $f$.

\cqfd

\noindent {\bf End of the proof of Theorem \ref{thm-pi-contract-pour-mu-grand}.} 
Let $\mu$ satisfying the hypothesis in Proposition \ref{prop-Bun-P-Bun-M-avec-N+d}, then 
$\ms Y_{N, d} = [S / H_{N, d}]$ (resp. $\ms Y_N = [S / H_N]$), where $H_{N, d} := R^0(pr_S)_* \mc Ker_{U_{\mc M}, N \times S + \Gamma_{\sum d x_i}}$ (resp. $H_N:= R^0(pr_S)_* \mc Ker_{U_{\mc M}, N \times S}$) is a unipotent group scheme over $S$. 

By Lemma \ref{lem-BH1-BH2-induce-H1-H2}, the two morphisms $b_1$ and $b_2$ in diagram (\ref{diagram-Cht-est-egalisateur-de-Bun-egal-mu}) are induced by two morphisms of group schemes $f_1, f_2: H_{N, d} \rightarrow H_{N}$. By Lemma \ref{lem-egalisateur-est-unipotent}, $\ms Z$ is isomorphic to $[H_N / H_{N, d}],$ where $H_{N, d}$ acts on $H_N$ by $h' \cdot h = f_1(h') h f_2(h')^{-1}$. 
\cqfd

\subsection{Cohomological statements}     \label{subsection-coho-statements}

\begin{defi}  \label{def-constant-C-tilde}
Let $d_W$ be the smallest integer in Proposition \ref{prop-d-assez-grand} such that the action of $G_{I, \infty}$ on $\Gr_{G, I, W}$ factors through $G_{I, d_W}$. We have defined the constants $C'(G, X, N, W)$ and $C(G, X, N, I, d_W)$ in Theorem \ref{thm-i-contract-pour-mu-grand} and Theorem \ref{thm-pi-contract-pour-mu-grand} respectively. We take:
$$\widetilde{C} (G, X, N, W):= Max\{  C'(G, X, N, W) , \; C(G, X, N, I, d_W)  \}$$
\end{defi}

\begin{defi} Let $\mu \in \wh{\Lambda}_{G^{\mr{ad}}}^{+, \Q}$. 
For any $j \in \Z$, we define degree $j$ cohomology sheaves
$$\mc H_{G, N, I, W}^{j, \, S_M(\mu)} = R^j(\mf{p}_G)_! ( \restr{\mc{F}_{G, N, I, W}^\Xi} {\Cht_{G, N, I, W}^{S_M(\mu)} / \Xi } )$$
$$\mc H_{M, N, I, W}^{' \, j, \, S_M(\mu)} =  R^j(\mf{p}_M ')_! ( \restr{\mc{F}_{M, N, I, W}^{ ' \,\Xi} } {\Cht_{M, N, I, W}^{' S_M(\mu)} / \Xi } ).$$
\end{defi}


\sssec{}
If $\langle  \mu, \alpha \rangle > \wt C(G, X, N, W) \text{ for all }  \alpha \in \Gamma_G - \Gamma_M$, then by Theorem \ref{thm-i-contract-pour-mu-grand}, the morphism $i^{' \, S_M(\mu)}: \Cht_{P, N, I, W}^{' \, S_M(\mu)} \rightarrow \Cht_{G, N, I, W}^{S_M(\mu)}$ is proper and schematic. Applying the construction in Section 3 to the truncation $S_M(\mu)$, we obtain a constant term morphism (in $D_c^b((X \sm N)^I, \Ql)$):
\begin{equation}    \label{equation-mc-C-G-P-S-M-mu}
\mc C_{G, N}^{P, \,  j, \, S_M(\mu)}: \mc H_{G, N, I, W}^{j, \, S_M(\mu)} \rightarrow \mc H_{M, N, I, W}^{' \, j, \, S_M(\mu)} .
\end{equation}

Here is the main result of Section 4:
\begin{prop}  \label{prop-TC-isom-mu-grand}
Let $P$ be a parabolic subgroup of $G$ and $M$ its Levi quotient. For $\mu \in \wh{\Lambda}_{G^{\mr{ad}}}^{+, \Q}$, if $\langle  \mu, \alpha \rangle > \wt C(G, X, N, W) \text{ for all }  \alpha \in \Gamma_G - \Gamma_M$, then for any $j$, morphism (\ref{equation-mc-C-G-P-S-M-mu}) is an isomorphism.
\end{prop}

\dem
By (\ref{morphism-TC-complex-p-G-F-G-vers-p-M-F-M}), $ \mc C_{G, N}^{P, \,  j, \, S_M(\mu)}$ is the composition of two morphisms:
$$ 
\begin{aligned}
R^j(\mf{p}_G)_!  (\restr{ \mc{F}_{G, N, I, W}^{\Xi} }  {\Cht_{G, N, I, W}^{S_M(\mu)}  } )  &  \xrightarrow{(1)}   
 R^j(\mf{p}_M')_!  (\pi^{' \, S_M(\mu)})_! (i^{' \, S_M(\mu)})^* ( \restr{ \mc{F}_{G, N, I, W}^{\Xi} } {\Cht_{G, N, I, W}^{S_M(\mu)}  } )  \\
& \xrightarrow{(2)} R^j(\mf{p}_M')_!  (\restr{ \mc{F}_{M, N, I, W}^{' \, \Xi} }  {\Cht_{M, N, I, W}^{' \, S_M(\mu)} } ) .
\end{aligned}
$$

The morphism (1) is induced by the composition of functors
$$ \begin{aligned}
R^j(\mf{p}_G)_!  \rightarrow  R^j(\mf{p}_G)_! (i^{' \, S_M(\mu)})_* (i^{' \, S_M(\mu)})^* & \simeq  R^j(\mf{p}_G)_! (i^{' \, S_M(\mu)})_! (i^{' \, S_M(\mu)})^*  \\ 
& \isom  R^j(\mf{p}_M')_!  (\pi^{' \, S_M(\mu)})_! (i^{' \, S_M(\mu)})^*
\end{aligned}
$$
defined in (\ref{morphism-TC-complex-adj}). By Theorem \ref{thm-i-contract-pour-mu-grand} and Lemma \ref{lem-coho-horosphere-i} below applied to $i^{' \, S_M(\mu)}$, the morphism (1) is an isomorphism.

The morphism (2) is induced by the morphism $$(\pi^{' \, S_M(\mu)})_! (i^{' \, S_M(\mu)})^* \mc{F}_{G, N, I, W}^{\Xi} \rightarrow \mc{F}_{M, N, I, W}^{' \, \Xi}$$ defined in (\ref{equation-F-G-vers-F-M-long}), which is a composition of the counit map $$\text{Co}: (\pi_d^{' \, S_M(\mu)})_! (\pi_d^{' \, S_M(\mu)})^! \rightarrow \Id$$ and some isomorphisms. By Theorem \ref{thm-pi-contract-pour-mu-grand} and
Lemma \ref{lem-Tr-est-iso-mu-assez-grand} below applied to $\pi_d^{' \, S_M(\mu)}$, the morphism (2) is an isomorphism.
\cqfd

\quad

\begin{lem}   \label{lem-coho-horosphere-i}
Let $f: \ms X \rightarrow \ms Y$ be a schematic finite universal homeomorphism of algebraic stacks, then the unit map $\Id \rightarrow f_*f^*$ is an isomorphism.
\cqfd
\end{lem}


\begin{lem}  \label{lem-Tr-est-iso-mu-assez-grand}
Let $f: \ms X \rightarrow \ms Y$ be an unipotent morphism of algebraic stacks (see Definition \ref{def-sch-en-gp-unipotent}), then the counit map $ f_!f^! \rightarrow \Id$ is an isomorphism.
\end{lem}
\dem
The proof consists of 4 steps.

(i) Using proper base change and the fact that $f$ is smooth, we reduce to the case when $\ms Y = \on{Spec} k$ is a point, thus $\ms  X = U_1 / U_2$ is a quotient of unipotent group schemes $U_1$ and $U_2$ over $k$.

Indeed, to prove the lemma, it is enough to prove that for any geometric point $i_y: y \rightarrow\ms Y$, the morphism $(i_y)^* f_! f^! \rightarrow (i_y)^*$ is an isomorphism. Form the following Cartesian square
\begin{equation}
\xymatrix{
f^{-1}(y)  \ar[d]^{ \wt{f} }  \ar[r]^{\wt{i_y} }
& \ms X \ar[d]^{f} \\
y  \ar[r]^{ i_y}
& \ms Y 
}
\end{equation}
Since $f$ is smooth, we have $f^! \simeq  f^*[2n](n)$ and $(\wt f)^! \simeq  (\wt f)^*[2n](n)$, where $n$ is the dimension of $f$.
We deduce that
\begin{equation}
(i_y)^* f_! f^! \simeq (\wt{f})_! (\wt i_y)^* f^! \simeq (\wt{f})_! (\wt i_y)^* (f)^* [2n](n) \simeq (\wt{f})_! (\wt{f})^* (i_y)^*[2n](n) \simeq  (\wt{f})_! (\wt{f})^! (i_y)^*
\end{equation}
where the first isomorphism is the proper base change (\cite{LO09} Theorem 12.1).
Thus it is enough to prove that $(\wt{f})_! (\wt{f})^! (i_y)^* \rightarrow (i_y)^*$ is an isomorphism.

(ii) We denote by $BU_2$ the classifying stack of $U_2$ over $k$.
Let $f_1: U_1 / U_2 \rightarrow BU_2$ and $f_2: BU_2 \rightarrow \on{Spec} k$ be the canonical morphisms. Then $f =f_2 \circ f_1$. We have a commutative diagram of functors:
\begin{equation}
\xymatrix{
f_!f^! = (f_2)_! (f_1)_!(f_1)^!(f_2)^! \ar[r]   \ar[d]
& \Id \\
(f_2)_! (f_2)^! \ar[ru]
}
\end{equation}
Thus it is enough to prove that the counit maps $(f_1)_! (f_1)^! \rightarrow \Id$ and $(f_2)_! (f_2)^! \rightarrow \Id$ are isomorphisms. 

(iii) Note that $f_1$ is a $U_1$-torsor over $BU_2$. By Definition \ref{def-sch-en-gp-unipotent}, we reduce to the case of $\mb A^1$-torsor. Using (i) again, we reduce to the case when $f_1$ is the map $\mb A^1 \rightarrow \on{Spec} k$, where it is clear that $(f_1)_! (f_1)^! \rightarrow \Id$ is an isomorphism.

(iv) Let $g_2: \on{Spec} k \rightarrow BU_2$ be the canonical morphism. Then $f_2 \circ g_2 \simeq \Id$. We have a commutative diagram of functors:
\begin{equation}
\xymatrix{
(f_2)_! (g_2)_!(g_2)^!(f_2)^! \ar[r]^{\quad \quad \quad \simeq}   \ar[d]
& \Id \\
(f_2)_! (f_2)^! \ar[ru]
}
\end{equation}
We deduce that to prove that $(f_2)_! (f_2)^! \rightarrow \Id$ is an isomorphism, it is enough to prove that $(g_2)_! (g_2)^! \rightarrow \Id$ is an isomorphism. Note that $g_2$ is a $U_2$-torsor over $BU_2$. Just like in (iii), we prove that $(g_2)_! (g_2)^! \rightarrow \Id$ is an isomorphism.
\cqfd

\quad

\begin{rem} 
In fact, to prove that the morphism (2) in Proposition \ref{prop-TC-isom-mu-grand} is an isomorphism, 
it is enough to write $\pi_d^{' \, S_M(\mu)}$ as the tower
$$\Cht_P^{S_M(\mu)} \xrightarrow{\pi_{d, m}} \cdots \rightarrow \wt{\Cht}_{P / U^{(j+1)}} \xrightarrow{\pi_{d, j}}  \wt{\Cht}_{P / U^{(j)}} \rightarrow \cdots \xrightarrow{\pi_{d, 0}} \wt{\Cht}_{M} $$
and prove that for each $j$, the morphism $\on{Co}: (\pi_{d, j})_! (\pi_{d, j})^! \rightarrow \Id$ is an isomorphism. For this, we only need the statement of Theorem \ref{thm-pi-contract-pour-mu-grand} for each $\pi_{d, j}$ (and replace unipotent group scheme by additive group scheme). The proof of such a statement still uses the three steps, but in step 2 Remark \ref{rem-U-comm-vs-non-comm} we only need to consider the case of commutative groups.
\end{rem}

\quad

\section{Finiteness of the cuspidal cohomology}

The goal of this section is to prove:
\begin{thm}  \label{thm-cusp-dim-fini-second}
The $\Ql$-vector space $H_{G, N, I, W}^{j, \; \on{cusp}}$ (defined in Definition \ref{def-cusp-coho}) has finite dimension.
\end{thm}

Theorem \ref{thm-cusp-dim-fini-second} will be a direct consequence of the following proposition.
\begin{prop}  \label{prop-H-cusp-inclus-dans-H-leq-mu}
Let $G, X, N, I, W$ as before.
There exists $\mu_0 \in \wh{\Lambda}_{G^{\mr{ad}}}^{+, \Q}$ (depending on $G, X, N, W, j$) such that
$$H_{G, N, I, W}^{j, \; \on{cusp} } \subset \mathrm{Im}(  H_{G, N, I, W}^{j, \; \leq \mu_0 }  \rightarrow H_{G, N, I, W}^j ).$$
\end{prop}
The proof of this proposition is essentially based on Proposition \ref{prop-TC-isom-mu-grand} and an induction argument on the semisimple rank of the group $G$.
We will present our strategy in Section 5.1 and give the proof in Sections \ref{subsection-begin-of-induction}-\ref{subsection-injectivity}.

\begin{nota}
In the remaining part of this section, to simplify the notations, we will omit the indices $N, I, W$.
\end{nota}

\subsection{Strategy of the proof}

\sssec{}   \label{subsection-choose-r}
We denote by $\wh R_{G^{\mr{ad}}}$ the coroot lattice of $G^{\mr{ad}}$. We have $\wh R_{G^{\mr{ad}}} \subset \wh\Lambda_{G^{\mr{ad}}}$. Let $\wh R_{G^{\mr{ad}}}^{+}:=\wh \Lambda_{G^{\mr{ad}}}^{+} \cap \wh R_{G^{\mr{ad}}}$. For any $r \in \N$, we have $\frac{1}{r} \wh R_{G^{\mr{ad}}}^+ \subset \wh \Lambda_{G^{\mr{ad}}}^{+, \Q}$ and 
$$\varinjlim _{\mu \in \wh \Lambda_{G^{\mr{ad}}}^{+, \Q}} H_{G}^{j, \, \leq \mu} = \varinjlim _{\mu \in \frac{1}{r} \wh R_{G^{\mr{ad}}}^+ } H_{G}^{j, \, \leq \mu} .$$

Let $\iota: \wh \Lambda_{Z_M / Z_G}^{\Q} \subset \wh\Lambda_{G^{\mr{ad}}}^{\Q}$ be the inclusion.
We fix $r$ such that $\bigcup_{P \subsetneq G} \iota \circ \on{pr}_P^{\mr{ad}}( \wh \Lambda_{G^{\mr{ad}}}^+ ) \subset \frac{1}{r} \wh R_{G^{\mr{ad}}}^+$, where $\on{pr}_P^{\mr{ad}}: \wh \Lambda_{G^{\mr{ad}}}^{\Q} \rightarrow\wh \Lambda_{Z_M / Z_G}^{\Q}$ is defined in (\ref{equation-pr-p-ad}).

\sssec{}   \label{subsection-def-P-alpha-M-alpha}
For any $\alpha \in \Gamma_G$, we denote by $\check{\alpha} \in \wh \Gamma_G$ the corresponding coroot, and vice versa.
Let $P_{\alpha}$ be the maximal parabolic subgroup with Levi quotient $M_{\alpha}$ such that $\Gamma_G - \Gamma_{M_{\alpha}} = \{ \alpha \}$.

In this section, for $\mu \in \wh \Lambda_{G^{\mr{ad}}}^{+, \Q}$, we will write $\mu - \frac{1}{r}\check{\alpha}$ instead of $\mu - \frac{1}{r}\Upsilon_G(\check{\alpha})$, where $\Upsilon_G: \wh{\Lambda}_G^{\Q} \rightarrow \wh{\Lambda}_{G^{\mr{ad}}}^{\Q}$ is defined in \ref{equation-Upsilon-G}.

\sssec{}     \label{subsection-H-leq-mu-to-H-j}
We have defined the inductive limits $H_G^j$ in Definition \ref{def-H-G-j} and $H_{M_\alpha}^{' \, j}$ in Definition \ref{def-H-M-prime-j}.
For any $\lambda \in \frac{1}{r} \wh R_{G^{\mr{ad}}}^+$, let $\mc I_{\lambda}: H_G^{j, \, \leq \lambda } \rightarrow H_G^j$ be the morphism to the inductive limit. Let $H_G^{j, \, \leq \lambda } \rightarrow H_{M_\alpha}^{' \, j}$ be the composition of morphisms $H_G^{j, \, \leq \lambda } \xrightarrow{\mc I_{\lambda} } H_G^j \xrightarrow{ C_G^{P_\alpha, \, j}} H_{M_\alpha}^{' \, j}$, where the second morphism is defined in Definition \ref{def-CT-cohomology}.

\sssec{}
Since for every $c \in H_{G}^{j } $, there exists $\lambda \in \wh R_{G^{\mr{ad}}}^+$ large enough such that $c \in \mathrm{Im}(  H_{G}^{j, \; \leq \lambda }  \rightarrow H_{G}^j )$, Proposition \ref{prop-H-cusp-inclus-dans-H-leq-mu} will be a direct consequence of (b) in the following proposition:

\begin{prop}   \label{prop-finiteness-a-b-c}
Let $G$ be a connected split reductive group. There exists a constant $C_G^0 \in \Q^{\geq 0}$ (depending on $G, X, N, W, j$), such that the following properties hold:

(a) Let $\mu \in  \frac{1}{r} \wh R_{G^{\mr{ad}}}^+ $ such that $\langle \mu , \gamma \rangle \geq C_G^0$ for all $\gamma \in \Gamma_G$.
Then for any $\alpha \in \Gamma_G$ such that $\mu - \frac{1}{r} \check{\alpha} \in \frac{1}{r} \wh R_{G^{\mr{ad}}}^+$ (which is automatic if $C_G^0 > \frac{2}{r}$), the morphism
$$\Ker( H_G^{j, \, \leq \mu - \frac{1}{r} \check{\alpha} } \rightarrow H_{M_\alpha}^{' \, j}) \rightarrow \Ker( H_G^{j, \, \leq \mu } \rightarrow H_{M_\alpha}^{' \, j})$$ is surjective.

(b) There exists $\mu_0 \in \frac{1}{r} \wh R_{G^{\mr{ad}}}^+ $ (depending on $C_G^0$), such that for any $\lambda \in \frac{1}{r} \wh R_{G^{\mr{ad}}}^+$ satisfying $\lambda \geq \mu_0$ and $\langle \lambda , \gamma \rangle \geq C_G^0$ for all $\gamma \in \Gamma_G$, the morphism
$$\Ker( H_G^{j, \, \leq \mu_0  } \rightarrow \prod_{P \subsetneq G} H_{M}^{' \, j}) \rightarrow \Ker( H_G^{j, \, \leq \lambda } \rightarrow \prod_{P \subsetneq G} H_{M}^{' \, j})$$ is surjective.

(c) There exists a constant $C_G \geq C_G^0$, such that for any $\lambda \in \frac{1}{r} \wh R_{G^{\mr{ad}}}^+ $ satisfying $\langle \lambda , \gamma \rangle \geq C_G$ for all $\gamma \in \Gamma_G$, the morphism $\mc I_{\lambda}: H_G^{j, \, \leq \lambda} \rightarrow H_G^{j}$ is injective.
\end{prop}

\sssec{}
The proof of Proposition \ref{prop-finiteness-a-b-c} uses an induction argument on the semisimple rank of the group $G$: firstly we prove the statements (a), (b) and (c) for every Levi subgroup of $G$ of rank $0$. Secondly we prove the key step: for $n \geq 1$, if (c) is true for all Levi subgroups of rank $n-1$, then (a) is true for all Levi subgroups of rank $n$. Then we deduce easily (a) $\Rightarrow$ (b) and (b) $\Rightarrow$ (c) for all Levi subgroups of rank $n$. 


\sssec{}   \label{subsection-compatibility-Xi-M-Xi-G}
As in \ref{subsection-fix-T-then-M-includ-in-P},
we fix a maximal torus $T \subset B$ and view the Levi quotient $M$ of a (standard) parabolic $P$ as a subgroup $M \subset P$. 

Recall that we have fixed $\Xi = \Xi_G \subset Z_G(\mb A)$ in \ref{subsection-def-Xi}.
Applying \ref{subsection-def-Xi} to each Levi subgroup $M$ of $G$, we fix $\Xi_M \subset Z_M(\mb A)$. Moreover, we choose $\Xi_M$ for different Levi subgroups in a compatible way: if $M_2$ is a Levi subgroup of $M_1$, then we have $\Xi_G \subset \Xi_{M_1} \subset \Xi_{M_2} \subset T(\mb A)$.

\subsection{Beginning of the induction: semisimple rank $0$}    \label{subsection-begin-of-induction}

\sssec{}
The only Levi subgroup of semisimple rank $0$ is the maximal torus $T$. Then $T^{\mr{ad}}$ is trivial and $\wh{\Lambda}_{T^{\mr{ad}}}^{+} = \wh{\Lambda}_{T^{\mr{ad}}}$ has only one element: $0$. 

The algebraic stack $\Cht_T / \Xi_T$ is of finite type. There is only one term in the inductive limit $H_T^j$, which is of finite dimension.

There is no constant term morphism for $T$. 
So we have $H_T^{j, \, \on{cusp} } = H_T^{j}$

\begin{lem}  \label{lem-recurrence-debut}
Take $C_T^0 =C_T =0$ and $\mu_0=0$.
Proposition \ref{prop-finiteness-a-b-c} is true for $T$.

\cqfd
\end{lem}

\subsection{From semisimple rank $n-1$ to $n$}

\begin{lem}   \label{lem-recurrence-a}
Let $G$ be a connected split reductive group of semisimple rank $n$. Suppose that Proposition \ref{prop-finiteness-a-b-c} (c) is true for every Levi quotient $M$ of $G$ of semisimple rank $n-1$, with a constant $C_M$. We take $C_G^0: =$
$$ Max \left\lbrace  \{  C_M \; | \; M \text{ Levi quotient of semisimple rank } n-1 \text{ of }  G  \}, \; \widetilde{C}(G, X, N, W) \right\rbrace $$
where $\wt C(G, X, N, W)$ is the constant defined in Definition \ref{def-constant-C-tilde}. Then for this constant $C_G^0$ Proposition \ref{prop-finiteness-a-b-c} (a) is true for $G$.
\end{lem}

We need some preparations before the proof of Lemma \ref{lem-recurrence-a}.

\sssec{}   \label{subsection-setting-P-alpha-M-alpha}
Let $\mu \in  \frac{1}{r} \wh R_{G^{\mr{ad}}}^+ $ such that $\langle \mu , \gamma \rangle \geq C_G^0$ for all $\gamma \in \Gamma_G$.
Let $\alpha \in \Gamma_G$ such that $\mu - \frac{1}{r} \check{\alpha} \in \frac{1}{r} \wh R_{G^{\mr{ad}}}^+$.
Let $P:=P_\alpha$ and $M:=M_\alpha$ as in \ref{subsection-def-P-alpha-M-alpha}. Note that $\Gamma_G - \Gamma_M = \{ \alpha \}$. 

\begin{lem}    \label{lem-S-2-minus-S-1-is-S-M-mu}
Let $S_1 = \{ \lambda \in \frac{1}{r} \wh R_{G^{\mr{ad}}}^+ | \lambda \leq \mu -  \frac{1}{r} \check{\alpha}  \}$ and $S_2 = \{ \lambda \in \frac{1}{r} \wh R_{G^{\mr{ad}}}^+  | \lambda \leq \mu \}$. Then
\begin{equation}    \label{equation-S-2-S-1-is-S-M-mu}
S_2 - S_1 = S_M(\mu) \cap ( \frac{1}{r} \wh R_{G^{\mr{ad}}}^+) ,
\end{equation}
where $S_M(\mu)$ is defined in Definition \ref{def-S-M-mu}.
\end{lem}
\dem
For any $\lambda \in S_2$, we have $\mu - \lambda = \sum_{\check{\gamma} \in \wh \Gamma_G} \frac{c_\gamma}{r} \check{\gamma}$ for some $c_\gamma \in \Z_{\geq 0}.$ Thus
\begin{equation}   \label{equation-mu-minus-alpha-lamba}
(\mu - \frac{1}{r} \check{\alpha} ) - \lambda = (\frac{c_\alpha}{r}  - \frac{1}{r}) \check{\alpha} + \sum_{\check{\gamma} \in \wh \Gamma_G, \, \check{\gamma} \neq \check{\alpha}} \frac{c_\gamma}{r} \check{\gamma} , \quad c_\gamma \in \Z_{\geq 0}.
\end{equation} 
If moreover $\lambda \notin S_1$, then in (\ref{equation-mu-minus-alpha-lamba}), there should be at least one coefficient strictly negative. So we must have $c_\alpha - 1 <0$. Since $c_\alpha \in \Z_{\geq 0}$, we must have $c_\alpha = 0$. We deduce that
$$\mu - \lambda = \sum_{\check{\gamma} \in \wh \Gamma_G, \, \check{\gamma} \neq \check{\alpha}} \frac{c_\gamma}{r} \check{\gamma} = \sum_{\check{\gamma} \in \wh \Gamma_M} \frac{c_\gamma}{r} \check{\gamma}, \quad c_\gamma \in \Z_{\geq 0}.$$
By Definition \ref{def-S-M-mu}, we have $\lambda \in S_M(\mu)$.
\cqfd

\begin{lem}    \label{lem-H-M-leq-mu-nu-to-H-M-nu-Xi-G-inj}
Let $\mu$ and $M$ as in \ref{subsection-setting-P-alpha-M-alpha}. Suppose that Proposition \ref{prop-finiteness-a-b-c} (c) is true for $M$. 
Then for any $j \in \Z$, the morphism $H _{M}^{' \, j, \, \leq \mu} \rightarrow H_M^{' \, j}$ is injective.
\end{lem}

The point of the proof of this lemma is to replace the quotient by $\Xi_M$ in (\ref{equation-H-M-leq-lamba-to-H-M-injective}) by the quotient by $\Xi_G$ in (\ref{equation-H-M-leq-lamba-to-H-M-injective-Xi-G}).

\dem
By Proposition \ref{prop-finiteness-a-b-c} (c) for $M$, for any $\lambda \in \frac{1}{r} \wh R_{M^{\mr{ad}}}^+ $ satisfying $\langle \lambda , \gamma \rangle > C_M$ for all $\gamma \in \Gamma_M$, the morphism 
\begin{equation}   \label{equation-H-M-leq-lamba-to-H-M-injective}
H_c^j ( \Cht_{M, \ov{\eta^I}}^{\leq^{M^{\mr{ad}}} \lambda} / \Xi_M, \mc F_M ) \rightarrow H_c^j ( \Cht_{M, \ov{\eta^I}} / \Xi_M, \mc F_M ) 
\end{equation}
is injective, where everything is defined as in Section \ref{subsection-coho-Cht-G} by replacing $G$ by $M$.

We can assume that $\Xi_M$ in \ref{subsection-compatibility-Xi-M-Xi-G} is small enough (containing $\Xi_G$). Then for any $\nu \in A_M$ (defined in \ref{subsection-def-Bun-M-nu}), the composition of morphisms
\begin{equation}    \label{equation-Cht-M-Xi-G-to-Cht-M-Xi-M}
\Cht_M^{\nu} / \Xi_G \rightarrow \Cht_M / \Xi_G \rightarrow \Cht_M/ \Xi_M
\end{equation}
is an open and closed immersion. 

(For the following discussion, see \cite{these} Illustration 7.4.4 for an example for $G = GL_3$.)
Let $\nu \leq \on{pr}_P^{\mr{ad}}(\mu)$. We use a special case of \ref{subsection-leq-mu-fix-nu-equal-leq-M-mu-nu}: 
By \ref{subsection-leq-in-Lambda-Z-M-Z-G}, we have $\on{pr}_P^{\mr{ad}} \circ \Upsilon_G ( \check{\alpha} ) > 0.$
Let $c_\alpha \in \Q_{\geq 0}$ 
such that $\on{pr}_P^{\mr{ad}} (\mu ) -  c_\alpha \on{pr}_P^{\mr{ad}} \circ \Upsilon_G ( \check{\alpha} ) = \nu.$
Let $\mu_{\nu} := \mu - c_\alpha \check{\alpha} $. 
For any $\lambda \in \frac{1}{r} \wh R_{G^{\mr{ad}}}$, the condition $\lambda \leq \mu$ and $\on{pr}_P^{ad}(\lambda) = \nu $ is equivalent to $\lambda \leq^{\ov{M}} \mu_{\nu}$.
We deduce that $\Cht_M^{\leq \mu, \, \nu} = \Cht_M^{\leq^{\ov M}  \mu_{\nu} }$.

Let $\Psi : \ov{M} \twoheadrightarrow M^{\mr{ad}}$. If $\mu_1 \leq \mu_2$, then $\mu_{1, \nu} \leq^{\ov M} \mu_{2, \nu}$ and $\Psi (\mu_{1, \nu} ) \leq^{M^{\mr{ad}} }  \Psi (\mu_{2, \nu} )$. 
For all $\gamma \in \Gamma_M$, since $\langle \check{\alpha} , \gamma \rangle \leq 0$, we have $\langle \mu_{\nu} , \gamma \rangle \geq \langle \mu , \gamma \rangle$. By hypothesis $\langle \mu , \gamma \rangle \geq C_G^0 \geq C_M$, so $\langle \mu_{\nu} , \gamma \rangle \geq C_M$. 
Then the injectivity of (\ref{equation-H-M-leq-lamba-to-H-M-injective}) with $\lambda = \Psi(\mu_{\nu})$ implies that the morphism 
\begin{equation}    \label{equation-H-M-leq-lamba-to-H-M-injective-Xi-G}
H_c^j ( \Cht_{M, \ov{\eta^I}}^{\leq \mu, \, \nu} / \Xi_G, \mc F_M ) \rightarrow H_c^j ( \Cht_{M, \ov{\eta^I}}^{\nu} / \Xi_G, \mc F_M ) 
\end{equation}
is injective. Note that we have defined $H_M^{j, \, \leq \mu, \, \nu} = H_c^j ( \Cht_{M, \ov{\eta^I}}^{\leq \mu, \, \nu} / \Xi_G, \mc F_M )$ in Definition \ref{def-H-M-j-leq-mu} and $H_M^{j, \, \nu} = H_c^j ( \Cht_{M, \ov{\eta^I}}^{\nu} / \Xi_G, \mc F_M )$ in Definition \ref{def-H-M-j-nu}.

Moreover, since $\Cht_M' = \Cht_M \overset{P(\mc O_N)} \times G(\mc O_N)$ is a disjoint union of copies of $\Cht_M$, we deduce that the morphism $H _{M}^{' \, j, \, \leq \mu, \, \nu } \rightarrow H_M^{' \, j, \, \nu}$ is also injective, where $H _{M}^{' \, j, \, \leq \mu, \, \nu }$ is defined in Definition \ref{def-mc-H-M-prime} and $H_M^{' \, j, \, \nu}$ is defined in Definition \ref{def-H-M-j-nu-prime}.

Note that by Lemma \ref{lem-Bun-M-leq-ad-mu-nu-non-vide}, for $\nu \notin A_M$ or $\nu \nleq \on{pr}_P^{\mr{ad}}(\mu)$, the cohomology group ${H}_{M}^{' \; j,  \; \leq \mu, \, \nu} =0$.
By Remark \ref{rem-lim-prod-vs-prod-lim}, we have a commutative diagram, 
$$
\xymatrix{
\varinjlim_{\mu} {H'}_{M}^{\; j, \; \leq\mu}  \ar[r]^{f}
&  \prod_{\nu \in \wh{\Lambda}_{Z_M / Z_G}^{\Q}}   {H'}_{M}^{\; j,  \; \nu} \\
{H}_{M}^{' \; j, \; \leq\mu} = \prod_{\nu \in \wh{\Lambda}_{Z_M / Z_G}^{\mu}}   {H'}_{M}^{\; j,  \; \leq \mu, \, \nu}   \ar[u]^{g}  \ar[ur]^{h}
}
$$
where $f$ is (\ref{equation-H-M-to-prod-H-M-nu}) and $h$ is induced component by component by $H _{M}^{' \, j, \, \leq \mu, \, \nu } \rightarrow H_M^{' \, j, \, \nu}$. By the above discussion, $h$ is injective. We deduce that the morphism $g$ is injective. 
\cqfd

\noindent {\bf Proof of Lemma \ref{lem-recurrence-a}.} 
The proof consists of 4 steps.

(1)
Let $S_1$ and $S_2$ as in Lemma \ref{lem-S-2-minus-S-1-is-S-M-mu}.
We define $\Cht_G^{S_2}$ and $\Cht_G^{S_1}$ as in \ref{subsection-appendix-def-Cht-G-S-Cht-M-S} (taking into account \ref{subsection-appendix-replace-Lambda-by-R}). We deduce from Lemma \ref{lem-S-2-minus-S-1-is-S-M-mu} that $\Cht_G^{S_2} - \Cht_G^{S_1}  = \Cht_G^{ S_M(\mu)}$ and $\Cht_M^{' \, S_2} - \Cht_M^{' \, S_1} = \Cht_M^{' \, S_M(\mu) }$.

We deduce from \ref{subsection-Bun-G-leq-mu-union-de-strate} that
$$\Bun_G^{=\lambda} \neq \emptyset   \Rightarrow   \Upsilon_G(\lambda) \in \bigcup_{P \subsetneq G} \iota \circ \on{pr}_P^{\mr{ad}}( \wh \Lambda_{G^{\mr{ad}}}^+ ) \subset \frac{1}{r} \wh R_{G^{\mr{ad}}}^+,$$
where the last inclusion follows from the choice of $r$ in \ref{subsection-choose-r}. We deduce that $\Cht_G^{=\lambda} = \emptyset$ if $\lambda \notin \frac{1}{r} \wh R_{G^{\mr{ad}}}^+$. Thus $\Cht_G^{S_2} = \Cht_G^{\leq \mu}$, $\Cht_G^{S_1} = \Cht_G^{ \leq \mu -  \frac{1}{r} \check{\alpha}} $, $\Cht_M^{' \, S_2} = \Cht_M^{' \, \leq \mu}$ and $\Cht_M^{' \, S_1} = \Cht_M^{' \,  \leq \mu -  \frac{1}{r} \check{\alpha}} $.

Applying Lemma \ref{lem-TC-suite-exacte-longue} to $S_1$ and $S_2$, we obtain a commutative diagram of cohomology groups, where the upper and lower lines are part of the long exact sequences in (\ref{diagram-TC-suite-exacte-longue-general}):
\begin{equation} \label{diagram-TC-suite-exacte-longue-S-M-mu}
\xymatrix{
H _{G}^{j, \, \leq \mu -  \frac{1}{r} \check{\alpha} }  \ar[r]   \ar[d]^{C_G^{P, \, j, \, \leq \mu -  \frac{1}{r} \check{\alpha}} }
& H _{G}^{j, \, \leq \mu}   \ar[r]   \ar[d]^{C_G^{P, \, j, \, \leq \mu} }
& H _{G}^{j, \,  S_M(\mu) }  \ar[d]^{ C_G^{P, \, j, \, S_M(\mu)} }  \\
H _{M}^{' \, j, \, \leq \mu -  \frac{1}{r} \check{\alpha}  }  \ar[r]   
& H _{M}^{' \, j, \, \leq \mu }   \ar[r]   
&  H _{M}^{' \, j, \, S_M(\mu) }  
}
\end{equation}
Note that if $\Cht_G^{S_M(\mu)}=\emptyset$, then the proof is finished.

\quad

(2) By the hypothesis of Lemma \ref{lem-recurrence-a}, $\langle \mu, \alpha \rangle \geq C_G^0 \geq \wt C(G, X, N, W)$. By Proposition \ref{prop-TC-isom-mu-grand}, for any $j$, the morphism $C_G^{P, \, j, \, S_M(\mu)}: H_G^{j, \, S_M(\mu)} \rightarrow H_M^{' \, j, \, S_M(\mu)}$ is an isomorphism. 

\quad

(3) We deduce from (\ref{equation-TC-mu-1-mu-2}) a commutative diagram:
\begin{equation}   \label{diagram-TC-Cht-mu-vers-Cht-limit-nu} 
\xymatrix{
 H _{G}^{j, \, \leq \mu}   \ar[r]^{\mc I_G}   \ar[d]_{   C_G^{P, \, j, \, \leq\mu}  }
&  H _{G}^{j }  \ar[d]^{  C_G^{P, \, j}  }  \\
 H _{M}^{' \, j,  \, \leq \mu}   \ar[r]^{\mc I_M}   
&  H _{M}^{' \, j} 
}
\end{equation}
By Lemma \ref{lem-H-M-leq-mu-nu-to-H-M-nu-Xi-G-inj}, the morphism $\mc I_M$ in (\ref{diagram-TC-Cht-mu-vers-Cht-limit-nu}) is injective.


\quad

(4) Let $a \in \Ker( H_G^{j, \leq \mu } \rightarrow H_{M}^{' \, j})$. 
%
By the commutativity of (\ref{diagram-TC-Cht-mu-vers-Cht-limit-nu}), $\mc I_M \circ C_G^{P, \, j,  \, \leq \mu } (a)  = C_G^{P, \, j} \circ \mc I_G (a)=0$. By step (3), $\mc I_M$ is injective. So $C_G^{P, \, j,  \, \leq \mu } (a)=0$. 

By the commutativity of (\ref{diagram-TC-suite-exacte-longue-S-M-mu}) and the isomorphism in step (2), we deduce that the image of $a$ in $H _{G}^{j, \, S_M(\mu) }$ is zero. So there exists $a' \in H_G^{j, \leq \mu -  \frac{1}{r} \check{\alpha}}$ whose image in $H_G^{j, \leq \mu}$ is $a$. 
\cqfd

\begin{rem} In fact, we have
$$
H_M^{' \, j, \, \leq \mu}  = \big( \prod_{\nu < \on{pr}_P^{\mr{ad}}(\mu) } H _{M}^{' \, j, \, \leq \mu, \, \nu } \big)  \oplus  H _{M}^{' \, j, \, \leq \mu, \, \on{pr}_P^{\mr{ad}}(\mu) }   = H _{M}^{' \, j, \, \leq \mu -  \frac{1}{r} \check{\alpha}  }  \oplus   H _{M}^{' \, j, \, S_M(\mu) }  
$$
Thus the bottom line of (\ref{diagram-TC-suite-exacte-longue-S-M-mu}) was canonically split.
\end{rem}

%

\quad

\begin{lem}  \label{lem-recurrence-b}
If the property (a) of Proposition \ref{prop-finiteness-a-b-c} is true for $G$, then the property (b) of Proposition \ref{prop-finiteness-a-b-c} is true for $G$.
\end{lem}
\dem
Let $\nabla(C_G^0)$ be the set of $\mu \in  \frac{1}{r} \wh R_{G^{\mr{ad}}}^+ $ such that $\langle \mu , \gamma \rangle > C_G^0$ for all $\gamma \in \Gamma_G$.  
Let $\Omega(C_G^0)$ be the set of $\mu \in \nabla(C_G^0)$ such that $\mu - \frac{1}{r} \check{\alpha} \notin \nabla(C_G^0)$ for all $\check{\alpha} \in \wh \Gamma_G$.
The set $\Omega(C_G^0)$ is bounded, thus is finite.
Let $\mu_0 \in  \frac{1}{r} \wh R_{G^{\mr{ad}}}^+ $ such that $\mu_0 > \mu$ for all $\mu \in \Omega(C_G^0)$. 

For any $\lambda \in \nabla(C_G^0)$, there exists a (zigzag) chain $\lambda = \lambda^{(0)} > \lambda^{(1)} > \cdots > \lambda^{(m-1)} > \lambda^{(m)} $ in $\frac{1}{r} \wh R_{G^{\mr{ad}}}^+$ for some $m \in \Z_{\geq 0}$ such that 

(i) for any $j$, we have $\lambda^{(j)} \in \nabla(C_G^0)$,

(ii) for any $j$, we have $\lambda^{(j)} - \lambda^{(j+1)} = \frac{1}{r} \check{\alpha}$ for some simple coroot $\check{\alpha} \in \wh \Gamma_G$.

(iii) $\lambda^{(m)} \in \Omega(C_G^0)$.

(Indeed, $\lambda^{(0)}$ satisfies (i). Suppose that we have already constructed a chain until $\lambda^{(j)}$ which satisfies (i) and (ii). If $\lambda^{(j)}$ satisfies (iii), we are done. If not, then there exists some $\check{\alpha} \in \wh \Gamma_G$ such that $\lambda^{(j)} - \frac{1}{r} \check{\alpha} \in \nabla(C_G^0)$. We define $\lambda^{(j+1)}:=\lambda^{(j)} - \frac{1}{r} \check{\alpha} $ and continue the process.)


Applying successively the property (a) of Proposition \ref{prop-finiteness-a-b-c} to $\lambda^{(0)} $, $\lambda^{(1)} $, $\cdots$,  until $\lambda^{(m)} $, we deduce that the morphism
$$\Ker( H_G^{j, \, \leq \lambda^{(m)}   } \rightarrow \prod_{P \subsetneq G} H_{M}^{' \, j}) \rightarrow \Ker( H_G^{j, \, \leq \lambda } \rightarrow \prod_{P \subsetneq G} H_{M}^{' \, j})$$ is surjective. Assume in addition that $\lambda \geq \mu_0$, then 
the morphism $H_G^{j, \, \leq \lambda^{(m)} } \rightarrow H_G^{j, \, \leq \lambda}$ factors through $H_G^{j, \, \leq \mu_0} $. We deduce the lemma.
\cqfd

\subsection{Injectivity}    \label{subsection-injectivity}


\begin{lem}  \label{lem-recurrence-c}
If the property (b) of Proposition \ref{prop-finiteness-a-b-c} is true for $G$, then the property (c) of Proposition \ref{prop-finiteness-a-b-c} is true for $G$.
\end{lem}

We need some preparations before the proof of Lemma \ref{lem-recurrence-c}.


\sssec{}
For $\mu \in \frac{1}{r} \wh R_{G^{\mr{ad}}}^+ $, let $\mc I_{\mu}: H_G^{j, \, \leq \mu } \rightarrow H_G^j$ be the morphism to the inductive limit as in \ref{subsection-H-leq-mu-to-H-j}.
For $\lambda \in \frac{1}{r} \wh R_{G^{\mr{ad}}}^+ $ such that $\lambda \geq \mu$, we denote by $\mc I_{\mu}^{\lambda}: H_G^{j, \leq \mu} \rightarrow H_G^{j, \leq \lambda}$ the morphism defined in \ref{subsection-H-G-leq-mu-1-to-leq-mu-2}.
We have $\Ker(\mc I_{\mu}^{\lambda}) \subset \Ker(\mc I_{\mu}) \subset H_G^{j, \leq \mu}.$

For $\lambda_2 \geq \lambda_1 \geq \mu$, we have
$\Ker(\mc I_{\mu}^{\lambda_1}) \subset \Ker(\mc I_{\mu}^{\lambda_2}).$

\begin{lem} \label{lem-limite-stationne-Ker}
Let $\mu \in \frac{1}{r} \wh R_{G^{\mr{ad}}}^+ $. There exists $\mu^{\sharp} \in \wh R_{G^{\mr{ad}}}^+$ such that $\mu^{\sharp} \geq \mu$ and $\Ker(\mc I_{\mu}^{\mu^{\sharp}}) = \Ker(\mc I_{\mu})$.
\end{lem}
\dem
We have the filtered system $\{  \Ker(\mc I_{\mu}^{\lambda}) \; | \; \lambda \in \frac{1}{r} \wh R_{G^{\mr{ad}}}^+ , \lambda \geq \mu   \}$ in $\Ker(\mc I_{\mu})$ and $\Ker(\mc I_{\mu}) =\varinjlim _{\lambda} \Ker(\mc I_{\mu}^{\lambda})$. Since $\Ker(\mc I_{\mu})$ is of finite dimension, the result is clear.

\cqfd



\begin{construction} \label{construction-C-G}
Let $\mu_0$ be the one in the property (b) of Proposition \ref{prop-finiteness-a-b-c}.
Choose $\mu_0^{\sharp} \in \frac{1}{r} \wh R_{G^{\mr{ad}}}^+ $ which satisfies Lemma \ref{lem-limite-stationne-Ker} for $\mu_0$. 
Let $C_G = \on{max} \{ C_G^0, \on{max}_{\gamma \in \Gamma_G} \{  \langle \mu_0^{\sharp}  , \gamma  \rangle \}  \}$. 
\end{construction}


\noindent {\bf Proof of Lemma \ref{lem-recurrence-c}. } 
Let $\lambda \in \frac{1}{r} \wh R_{G^{\mr{ad}}}^+ $ such that $\langle \lambda , \gamma \rangle \geq C_G$ for all $\gamma \in \Gamma_G$. 
By Construction \ref{construction-C-G}, $ \langle \lambda - \mu_0^{\sharp} , \gamma \rangle = \langle \lambda , \gamma \rangle - \langle \mu_0^{\sharp} , \gamma \rangle \geq C_G - \langle \mu_0^{\sharp} , \gamma \rangle \geq 0$ for all $\gamma \in \Gamma_G$. Thus $\mu_0^{\sharp} \leq \lambda$.
Consider the morphisms:
$$H_G^{j, \, \leq \mu_0} \rightarrow  H_G^{j, \, \leq \mu_0^{\sharp}} \rightarrow  H_G^{j, \, \leq \lambda} \rightarrow H_G^j.$$
We have $\Ker(\mc I_{\mu_0}^{\mu_0^{\sharp}}) \subset \Ker(\mc I_{\mu_0}^{\lambda}) \subset \Ker(\mc I_{\mu_0})$. By Lemma \ref{lem-limite-stationne-Ker}, $\Ker(\mc I_{\mu_0}^{\mu_0^{\sharp}}) = \Ker(\mc I_{\mu_0})$, hence $\Ker(\mc I_{\mu_0}^{\lambda}) = \Ker(\mc I_{\mu_0})$. 

For any element $b \in \Ker (H_G^{j, \leq \lambda} \rightarrow H_G^j )$, we have $b \in \Ker( H_G^{j, \leq \lambda } \rightarrow \prod H_{M}^{' \, j})$. By the property (b) of Proposition \ref{prop-finiteness-a-b-c}, $b$ is the image of an element $b_0 \in \Ker( H_G^{j, \leq \mu_0 } \rightarrow \prod H_{M}^{' \, j})$. We have $b_0 \in \Ker(\mc I_{\mu_0}) = \Ker(\mc I_{\mu_0}^{\lambda})$, so its image $b$ in $H_G^{j, \leq \lambda}$ is zero. This implies that the morphism $H_G^{j, \leq \lambda} \rightarrow H_G^j$ is injective.
\cqfd

\quad

\section{Rational Hecke-finite cohomology}

In this section, we will define a subspace $H_{G, N, I, W}^{j, \; \on{Hf-rat}}$ of $H_{G, N, I, W}^{j}$ and prove:
\begin{prop}   \label{prop-cusp-egal-Hfrat-section-6}
The two $\Ql$-vector subspaces $H_{G, N, I, W}^{j, \; \on{cusp}}$ and $H_{G, N, I, W}^{j, \; \on{Hf-rat}}$ of $H_{G, N, I, W}^{j}$ are equal.
\end{prop} 

In Section 6.1 we give some preparations. In Section 6.2 we show that the constant term morphisms commute with the action of the Hecke algebra. Using this, in Section 6.3 we prove Proposition \ref{prop-cusp-egal-Hfrat-section-6}. 

\quad

In Section 6, all the stacks are restricted to $\ov{\eta^I}$.

\subsection{Compatibility of constant term morphisms and level change}

%
%

\sssec{}   \label{subsection-Cht-G-K-I-W}
Let $K$ be a compact open subgroup of $G(\mb O)$. 
Let $N$ be a level such that $K_N \subset K$. We define $$\Cht_{G, K, I, W} := \Cht_{G, N, I, W} / \big(   K / K_N  \big) .$$ It is independent of the choice of $N$. 

Let $d \in \N$ large enough as in Proposition \ref{prop-d-assez-grand}, we define $\mc F_{G, K, I, W}$ to be the inverse image of $\mc S_{G, I, W}^d$ by $\epsilon_{K, d}: \Cht_{G, K, I, W} \rightarrow [ G_{I, d} \backslash \Gr_{G, I, W}] $. Just as in Remark \ref{rem-F-G-I-W-independent-of-d}, $\mc F_{G, K, I, W}$ is independent of $d$. Similarly we define $\mc F_{G, K, I, W}^{\Xi}$ over $\Cht_{G, K, I, W} / \Xi$. We define $H_{G, K, I, W}^j:=\varinjlim _{\mu} H_c^j(\Cht_{G, K, I, W}^{\leq \mu} / \Xi, \mc F_{G, K, I, W}^{\Xi})$.

\sssec{}    \label{subsection-adj-pr-G-K-prime-K}
Let $K' \subset K$ be two compact open subgroups of $G(\mb O)$. The inclusion $K' / K_N \hookrightarrow K / K_N$ induces a morphism $\on{pr}^G_{K', K}: \Cht_{G, K', I, W}  \rightarrow \Cht_{G, K, I, W} $. Note that all the stacks are restricted to $\ov{\eta^I}$. Morphism $\on{pr}^G_{K', K}$ is finite étale of degree the cardinality of $K/K'$. The following diagram is commutative:
$$
\xymatrixrowsep{1pc}
\xymatrixcolsep{1pc}
\xymatrix{
\Cht_{G, K', I, W}   \ar[rr]^{\on{pr}^G_{K', K}}   \ar[rd]_{\epsilon_{K', d}}
& & \Cht_{G, K, I, W}  \ar[ld]^{\epsilon_{K, d}} \\
&  [ G_{I, d} \backslash \Gr_{G, I, W}]  
}
$$

Note that $(\on{pr}^G_{K', K})_* = (\on{pr}^G_{K', K})_!$ and
\begin{equation*}  \label{equation-inverse-image-of-F-n-is-F-n'}
(\on{pr}^G_{K', K})^* \mc F_{G, K, I, W} = (\on{pr}^G_{K', K})^* (\epsilon_{K, d})^* \mc S_{G, I, W}^d = (\epsilon_{K', d})^* \mc S_{G, I, W}^d =\mc F_{G, K', I, W} .
\end{equation*}
The adjunction morphism $\on{adj}(\on{pr}^G_{K', K}): \Id \rightarrow (\on{pr}^G_{K', K})_* (\on{pr}^G_{K', K})^*$ induces an (injective) morphism of cohomology groups, which we still denote by 
\begin{equation*}
\on{adj}(\on{pr}^G_{K', K}): H _{G, K, I, W}^j \rightarrow H _{G, K', I, W}^j.
\end{equation*}
Note that $(\on{pr}^G_{K', K})^!  = (\on{pr}^G_{K', K})^* $. The counit morphism (in this case equal to the trace map) $\on{Co}(\on{pr}^G_{K', K}): (\on{pr}^G_{K', K})_! (\on{pr}^G_{K', K})^! \rightarrow \Id$ induces a (surjective) morphism of cohomology groups, which we still denote by $$\on{Co}(\on{pr}^G_{K', K}): H _{G, K', I, W}^j \rightarrow H _{G, K, I, W}^j .$$

\sssec{}   \label{subsection-G-A-acts-on-Cht-G-infty}
Let $v$ be a place in $X$. 
Let $N = N^v + nv$. Taking projective limit over $n$, we define $$\varprojlim_{n} \Cht_{G, N^v + nv , I, W} $$
Let $g \in G(F_v)$. The right action of $g$ (by left multiplication by $g^{-1}$) induces an isomorphism 
\begin{equation*}   
\varprojlim_{n} \Cht_{G, N^v + nv , I, W}   \isom \varprojlim_{n} \Cht_{G, N^v + nv , I, W}  \quad (\mc G \rightarrow \ta \mc G, \psi^v, \psi_v) \mapsto (\mc G' \rightarrow \ta \mc G', \psi^v, \psi'_v )
\end{equation*}
where $\psi^v$ (resp. $\psi_v$) is the level structure outside $v$ (resp. on $v$). The $G$-bundle $\mc G'$ is defined by gluing $\restr{G}{\Gamma_{\infty v} }$ and $\restr{\mc G}{ X -v }$ by $\restr{G}{\Gamma_{\infty v} -v} \xrightarrow{g} \restr{G}{\Gamma_{\infty v} -v} \underset{\sim}{\xleftarrow{\psi_v}} \restr{\mc G}{ \Gamma_{\infty v} -v }$. We have $\psi_v' = g^{-1} \circ \psi_v$.

Let $$\Cht_{G, \infty, I, W}:= \varprojlim_{N} \Cht_{G, N, I, W} . $$
Similarly, $\Cht_{G, \infty, I, W}$ is equipped with an action of $G(\mb A)$.

\sssec{}   \label{subsection-Cht-P-K-'}
Let $P$ be a parabolic subgroup of $G$ and $M$ its Levi quotient. 
We define $$\Cht_{P, \infty, I, W}:= \varprojlim_{N} \Cht_{P, N, I, W} . $$
Just as in \ref{subsection-G-A-acts-on-Cht-G-infty}, $\Cht_{P, \infty, I, W}$ is equipped with an action of $P(\mb A)$. 
For any compact open subgroup $K \subset G(\mb O)$, we define 
\begin{equation}   \label{equation-Cht-P-K-I-W-'}
\Cht_{P, K, I, W}' := \Cht_{P, \infty, I,  W} \overset{P(\mb O) } \times G(\mb O) / K .
\end{equation}

We have a morphism
\begin{equation}
\Cht_{P, \infty, I,  W} \overset{P(\mb O) } \times G(\mb O)  \rightarrow   \Cht_{G, \infty, I, W} 
\end{equation}
by sending $\big( (\mc P, \psi_P) \rightarrow (\ta \mc P, \ta \psi_P), g \in G(\mb O)  \big)$ to $\big( (\mc G, g^{-1} \circ \psi_G) \rightarrow (\ta \mc G, g^{-1} \circ \ta \psi_G)  \big)$, where $\mc G = \mc P \overset{P}{\times} G$ and $\psi_G = \psi_P \overset{P}{\times} G$. It induces a morphism
\begin{equation}   \label{equation-Cht-P-infty-'-to-fiber-prod}
\Cht_{P, \infty, I,  W} \overset{P(\mb O) } \times G(\mb O)  \rightarrow  \Cht_{P, I, W} \underset{\Cht_{G, I, W}} \times \Cht_{G, \infty, I, W} .
\end{equation}
This is a $G(\mb O)$-equivariant morphism of $G(\mb O)$-torsors over $\Cht_{P, I, W}$, where $G(\mb O)$ acts on the LHS of (\ref{equation-Cht-P-infty-'-to-fiber-prod}) by right action (right multiplication) on $G(\mb O)$ and acts on the RHS of (\ref{equation-Cht-P-infty-'-to-fiber-prod}) by the right action on $\Cht_{G, \infty, I, W}$ defined in \ref{subsection-G-A-acts-on-Cht-G-infty}. Thus (\ref{equation-Cht-P-infty-'-to-fiber-prod}) is an isomorphism. We have 
$$ \Cht_{P, \infty, I,  W} \overset{P(\mb O) } \times G(\mb O) / K  \isom  \Cht_{P, I, W} \underset{\Cht_{G, I, W}} \times \Cht_{G, \infty, I, W} /K ,$$
i.e. 
\begin{equation}   \label{equation-Cht-P-K-'-is-fiber-prod}
\Cht_{P, K, I, W}' = \Cht_{P, I, W} \underset{\Cht_{G, I, W}} \times \Cht_{G, K, I, W}.
\end{equation}

When $K=K_N$ for some level $N$, we have $\Cht_{P, N, I, W} = \Cht_{P, \infty, I,  W} / K_{P, N} $, where $K_{P, N}:=K_N \cap P(\mb O)$. We deduce that $\Cht_{P, K_N, I, W}'$ defined in (\ref{equation-Cht-P-K-I-W-'}) coincides with $\Cht_{P, N, I, W}'$ defined in Definition \ref{def-Cht-P-M-'}.

\sssec{}   \label{subsection-Cht-M-N-'-rewrite}
We define $$\Cht_{M, \infty, I, W}:= \varprojlim_{N} \Cht_{M, N, I, W} . $$ Just as in \ref{subsection-G-A-acts-on-Cht-G-infty}, $\Cht_{M, \infty, I, W}$ is equipped with an action of $M(\mb A)$. 
Recall that for any level $N$, in Definition \ref{def-Cht-P-M-'}, we defined $\Cht_{M, N, I, W}' = \Cht_{M, N, I, W} \overset{P(\mc O_N)  } \times G(\mc O_N) $. 
Let $K_{U, N}:=K_N \cap U(\mb O)$ and $K_{M, N}:=K_{P, N} / K_{U, N}$. Taking into account that $\Cht_{M, N, I, W} = \Cht_{M, \infty, I, W} / K_{M, N}$, we deduce
\begin{equation}   \label{equation-Cht-M-infty-P-O-quotient-K-U-N}
\Cht_{M, N, I, W}' = \Cht_{M, \infty, I, W} \overset{P(\mb O) / K_{U, N} } \times G(\mb O) / K_N .
\end{equation}
When we consider the action of the Hecke algebras in \ref{subsection-Hecke-G-acts-on-H-M-'} in the next section, we will need some functoriality on $K_N$. For this reason, we rewrite (\ref{equation-Cht-M-infty-P-O-quotient-K-U-N}) in the following way. Note that $K_N$ is normal in $G(\mb O)$. The stabilizer of any $P(\mb O)$-orbit in $G(\mb O) / K_N$ is $K_{P, N}$. We deduce from (\ref{equation-Cht-M-infty-P-O-quotient-K-U-N}) that
\begin{equation}  \label{equation-Cht-M-infty-quotient-K-P-N-K-U-N}
\begin{aligned}
\Cht_{M, N, I, W}' & = \bigsqcup_{P(\mb O) \text{-orbits in } G(\mb O) / K_N} \Cht_{M, \infty, I, W} / (K_{P, N} / K_{U, N}) \\
& = \bigsqcup_{P(\mb A) \text{-orbits in } G(\mb A) / K_N} \Cht_{M, \infty, I, W} / (K_{P, N} / K_{U, N})
\end{aligned}
\end{equation}
The second equation is because that $P(\mb O) \backslash G(\mb O) = P(\mb A) \backslash G(\mb A)$, and that in each $P(\mb A)$-orbit in $G(\mb A) / K_N$, we can choose a representative in $G(\mb O) / K_N$.

In the following, we want to generalize (\ref{equation-Cht-M-infty-quotient-K-P-N-K-U-N}) for any compact open subgroup $K \subset G(\mb O)$ (which may not be normal in $G(\mb O)$).

\sssec{}   \label{subsection-def-H-M-S-I-W-'}
Let $\mc D$ be the category of discrete sets $S$ equipped with a continuous action of $P(\mb A)$ with finitely many orbits such that the stabilizer of any point is conjugated to some open subgroup of finite index in $P(\mb O)$.
In particular, for any compact open subgroup $K \subset G(\mb O)$, the set $S = G(\mb A) / K$ is an object in $\mc D$.

For any $S \in \mc D$, we define functorially the cohomology group $H_{M, S, I, W}'$ in the following way.


When $S$ has only one orbit, choose a point $s \in S$, let $H$ be the stabilizer of $s$. Then $H$ is a subgroup of $P(\mb A)$ conjugated to some open subgroup of finite index in $P(\mb O)$. We have $S = P(\mb A) / H$. Let $R$ be a subgroup of finite index in $H \cap U(\mb A)$ and normal in $H$. 
By \ref{subsection-Cht-M-N-'-rewrite}, $\Cht_{M, \infty, I, W}$ is equipped with an action of $M(\mb A)$, thus an action of $P(\mb A)$ by the projection $P(\mb A) \twoheadrightarrow M(\mb A)$. Note that $R \subset U(\mb A)$ acts trivially on $\Cht_{M, \infty, I, W}$. We define a Deligne-Mumford stack 
$$\Cht_{M, \infty, I, W} / (H / R).$$

We define the cohomology group $H_{M, S, R, I, W}'$ as in Definition \ref{def-H-M-j} for $\Cht_{M, \infty, I, W} / (H / R)$ (instead of $\Cht_{M, N, I, W}$). 
Concretely, we have a morphism $\Cht_{M, \infty, I, W} / (H / R) \Xi \rightarrow \Cht_{M, I, W} / \Xi$, where $\Cht_{M, I, W}$ is the stack of $M$-shtukas without level structure.
Let $\mc F_{M, \infty, I, W}^{\Xi}$ be the inverse image of $\mc F_{M, I, W}^{\Xi}$. We define 
$$H_{M, S, R, I, W}^{' \, j}:=  \varinjlim_{\mu} \prod_{\nu} H_c^j(\Cht_{M, \infty, I, W}^{\leq \mu, \, \nu} / (H / R) \Xi, \mc F_{M, \infty, I, W}^{\Xi}) .$$

Let $R_1 \subset R_2$ be two subgroups of finite index in $H \cap U(\mb A)$ and normal in $H$. The projection $H / R_1 \twoheadrightarrow H / R_2$ induces a morphism
$$\mf q_{R_1, R_2}: \Cht_{M, \infty, I, W} / (H / R_1) \rightarrow \Cht_{M, \infty, I, W} / (H / R_2) .$$
It is a gerbe for the finite $q$-group $R_2 / R_1$. The counit morphism (which is equal to the trace map because $\mf q_{R_1, R_2}$ is smooth of dimension $0$) $\on{Co}(\mf q_{R_1, R_2}): (\mf q_{R_1, R_2})_! (\mf q_{R_1, R_2})^! \rightarrow \Id$ is an isomorphism. Indeed, just as in the proof (i) of Lemma \ref{lem-Tr-est-iso-mu-assez-grand}, by proper base change and the fact that $\mf q_{R_1, R_2}$ is smooth, we reduce to the case of Lemma \ref{lem-gerbe-neutre-trace-is-id} below with $\Gamma = R_2 / R_1$. 
The morphism $\on{Co}(\mf q_{R_1, R_2})$ induces an isomorphism of cohomology groups
\begin{equation}   \label{equation-H-M-S-R-1-R-2-isom}
H_{M, S, R_1, I, W}' \isom H_{M, S, R_2, I, W}' .
\end{equation}
%
We define $H_{M, S, I, W}'$ to be any $H_{M, S, R, I, W}'$, where we identify $H_{M, S, R_1, I, W}'$ and $H_{M, S, R_2, I, W}'$ by (\ref{equation-H-M-S-R-1-R-2-isom}). 

Recall that $S$ has only one orbit. $H_{M, S, I, W}'$ is independent of the choice of the point $s$ in $S$. In fact, let $s_1, s_2$ be two point of $S$ and $H_1$ (resp. $H_2$) be the stabilizer of $s_1$ (resp. $s_2$), then $H_2 = p^{-1}H_1p$ for some $p\in P(\mb A)$. The action of $p$ induces an isomorphism $\Cht_{M, \infty, I, W} / (H_1 / R)  \isom \Cht_{M, \infty, I, W} / (p^{-1}H_1p / p^{-1}Rp) $. We deduce an isomorphism of cohomology groups by the adjunction morphism.

In general, $S = \sqcup_{\alpha \in A} \alpha $ is a finite union of orbits, we define
$$H_{M, S, I, W}' := \bigoplus_{\alpha \in A} H_{M, \alpha, I, W}' .$$

When $S = G(\mb A) / K$ for some compact open subgroup $K$ in $G(\mb O)$, we write
\begin{equation}
H_{M, K, I, W}' := H_{M, S, I, W}' .
\end{equation}

\begin{lem}   \label{lem-gerbe-neutre-trace-is-id}
Let $\Gamma$ be a finite group over an algebraically closed field $k$ over $\Fq$. We denote by $B\Gamma$ the classifying stack of $\Gamma$ over $k$. Let $\mf q: B\Gamma \rightarrow \on{Spec} k$ be the structure morphism.
Then the counit morphism (equal to the trace map) $\on{Co}(\mf q): \mf q_! \mf q^! \rightarrow \Id$ of functors on $D_c(\on{Spec} k , \Ql )$ is an isomorphism.
\end{lem}
\dem
$\on{Co}(\mf q)$ is the dual of the adjunction morphism $\on{adj}(\mf q): \Id \rightarrow \mf q_* \mf q^*$. For any $\mc F \in D_c(\on{Spec} k , \Ql )$, $\mf q^*\mc F$ is a complex $F$ of $\Gamma$-modules with trivial action of $\Gamma$. Since $H^j(B\Gamma, \mf q^* \mc F) = H^j(\Gamma, F)$ (group cohomology), we have $H^{0}(B\Gamma, \mf q^* \mc F) = F^{\Gamma} = F$ and $H^{j}(B\Gamma, \mf q^* \mc F) =0$ for $j > 0$. So $\on{adj}(\mf q)$ is an isomorphism. By duality, we deduce the lemma.
\cqfd

\sssec{}   \label{subsection-def-pi-S-R}
Let $S \in \mc D$. We define 
\begin{equation}   \label{equation-Cht-P-S-I-W-'}
\Cht_{P, S, I, W}':=\Cht_{P, \infty, I, W} \overset{P(\mb A) } \times S .
\end{equation}
For each orbit $\alpha$ in $S$, choose a representative, let $H^{\alpha}$ be the stabilizer (well defined up to conjugation). Then 
$$\Cht_{P, S, I, W}' =\bigsqcup_{\alpha \in \{ P(\mb A) \text{-orbits in } S \} } \Cht_{P, \infty, I, W} / H^{\alpha} .$$

For each $\alpha$, let $R^{\alpha}$ be a subgroup of finite index in $H^{\alpha} \cap U(\mb A)$ and normal in $H^{\alpha}$. 
Let $R = (R^{\alpha})_{ \alpha \in \{ P(\mb A) \text{-orbits in } S \} } $. We define
\begin{equation}   \label{equation-def-Cht-M-S-R-'}
\Cht_{M, S, R, I, W}' := \bigsqcup_{\alpha \in \{ P(\mb A) \text{-orbits in } S \}  } \Cht_{M, \infty, I, W}  / (H^{\alpha} / R^{\alpha})  
\end{equation}

For each $\alpha$, we have morphisms of prestacks
\begin{equation}   \label{equation-Cht-P-infty-quotient-H-to-Cht-M-infty-quotient-H-R}
\Cht_{P, \infty, I, W} / H^{\alpha} \rightarrow \Cht_{M, \infty, I, W} / H^{\alpha} \rightarrow \Cht_{M, \infty, I, W}  / (H^{\alpha} / R^{\alpha}) ,
\end{equation}
where the first and third prestacks are Deligne-Mumford stacks, while the second is only a prestack.
Taking union over all the orbits, we deduce from (\ref{equation-Cht-P-infty-quotient-H-to-Cht-M-infty-quotient-H-R}) a morphism
\begin{equation}  \label{equation-pi-S-R}
\pi_{S, R}: \Cht_{P, S, I, W}' \rightarrow \Cht_{M, S, R, I, W}'
\end{equation}

In particular, when $S = G(\mb A) / K_N$, the stack $\Cht_{P, S, I, W}'$ coincides with $\Cht_{P, N, I, W}'$. For every orbit $\alpha$, we can choose a representative in $G(\mb O) / K_N$ (so that $H^{\alpha} = K_{P, N}$) and choose $R^{\alpha} = K_{U, N}$. Then $\Cht_{M, S, R, I, W}'$ coincides with $\Cht_{M, N, I, W}'$, $H_{M, S, R, I, W}'$ coincides with $H_{M, N, I, W}'$ defined in Definition \ref{def-H-M-prime-j}, and (\ref{equation-pi-S-R}) coincides with $\pi'$ defined in (\ref{equation-Cht-G-P-M-prime}).

\sssec{}   \label{subsection-def-C-G-P-S-R}
For any compact open subgroup $K \subset G(\mb O)$, let $S = G(\mb A) / K$. Note that in this case we have $\Cht_{P, K, I, W}' = \Cht_{P, S, I, W}'$. 
For any $R$ as in \ref{subsection-def-pi-S-R}, we have morphisms
\begin{equation}    \label{equation-Cht-G-P-M-K}
\xymatrixrowsep{1pc}
\xymatrixcolsep{4pc}
\xymatrix{
& \Cht_{P, S, I, W}' \ar[ld]_{i_K}   \ar[rd]^{\pi_{S, R}}  \ar[dd]^{\mf p_P} \\
\Cht_{G, K, I, W}  \ar[rd]^{\mf p_G}
& & \Cht_{M, S, R, I, W}'  \ar[ld]^{\mf p_M}  \\
& \ov{\eta^I}
}
\end{equation} 
Just as in Proposition \ref{prop-Varshavsky-proper} and Remark \ref{rem-i-proper-schematic}, the morphism $i_K$ is schematic and proper. Apply the construction in Section 3 to (\ref{equation-Cht-G-P-M-K}). Similarly to (\ref{diagram-cht-P-cht-M-pi-d}), we have
\begin{equation}  
\xymatrixrowsep{2pc}
\xymatrixcolsep{2pc}
\xymatrix{
\Cht_{P, S, I, W}'  \ar[rdd]_{\epsilon_{P, d}}    \ar[rd]|-{\pi_{S, R, d}}   \ar[rrd]^{\pi_{S, R}}
& 
&  \\
& \widetilde{\Cht}_{M, S, R, I, W}'     \ar[d]^{ \widetilde{ \epsilon_{M, d } } } \ar[r]_{ \widetilde{  \pi^0_{d}  } }
& \Cht_{M, S, R, I, W}'  \ar[d]^{\epsilon_{M, d}}  \\
&  [ P_{I, d} \backslash  \Gr_{P, I, W}]  \ar[r]_{ \ov{\pi^0_d}  }
&  [  M_{I, d} \backslash  \Gr_{M, I, W}]
}
\end{equation}
where $\pi_{S, R, d}$ is smooth. Let $\mc F_G$ be the canonical Satake sheaf on $\Cht_{G, K, I, W}$ and $\mc F_M$ be the canonical Satake sheaf on $\Cht_{M, S, R, I, W}'$.
We construct a morphism $c_{G, K}^P: (\pi_{S, R})_! (i_K)^* \mc F_G   \rightarrow  \mc F_M$ similar to (\ref{equation-base-change-with-trace}) and (\ref{equation-F-G-vers-F-M-long}). 
Namely, $c_{G, K}^P$ is the composition of some isomorphisms and the counit morphism $(\pi_{S, R, d})_! (\pi_{S, R, d})^! \rightarrow \Id$. Note that since $\pi_{S, R, d}$ is smooth, the composition $(\pi_{S, R, d})_! (\pi_{S, R, d})^* [2m](m) \isom (\pi_d)_! (\pi_d)^! \rightarrow \Id$ is the trace map in \cite{sga4}  XVIII 2, where $m$ is the dimension of $\pi_{S, R, d}$.

Similar to (\ref{morphism-TC-complex-p-G-F-G-vers-p-M-F-M}), we have a composition of morphisms of functors in $D_c^b( \ov{\eta^I} , \Ql)$:
$$(\mf p_G)_! \mc F_G \xrightarrow{\on{adj}(i_K)} (\mf p_G)_! (i_K)_* (i_K)^* \mc F_G \simeq  (\mf p_M)_! (\pi_{S, R})_! (i_K)^* \mc F_G   \xrightarrow{c_{G, K}^P} (\mf p_M)_! \mc F_M .$$

We define $$H_{P, K, I, W}' :=H_{P, S, I, W}' := (\mf p_P)_!  (i_K)^* \mc F_G .$$
The morphism $(\mf p_G)_! \mc F_G \xrightarrow{\on{adj}(i_K)} (\mf p_G)_! (i_K)_* (i_K)^* \mc F_G $ induces a morphism
\begin{equation}   \label{equation-H-G-to-H-P}
H _{G, K, I, W}^j \rightarrow H_{P, S, I, W}'
\end{equation}
The morphism $(\mf p_M)_! (\pi_{S, R})_! (i_K)^* \mc F_G   \xrightarrow{c_{G, K}^P} (\mf p_M)_! \mc F_M $ induces a morphism 
\begin{equation}  \label{equatino-H-P-to-H-M}
H_{P, S, I, W}' \rightarrow H _{M, S, R,  I, W}^{' \, j} .
\end{equation}
We define the constant term morphism to be the composition of (\ref{equation-H-G-to-H-P}) and (\ref{equatino-H-P-to-H-M})
\begin{equation}   \label{equation-def-C-G-P-S-R}
C_{G, S, R}^{P, \, j} : H _{G, K, I, W}^j \rightarrow H _{M, S, R, I, W}^{' \, j} .
\end{equation}

For $R_1 \subset R_2$ as in \ref{subsection-def-H-M-S-I-W-'}, the following diagram is commutative 
\begin{equation}    \label{diagram-CT-R-1-R-2}    
\xymatrixrowsep{2pc}
\xymatrixcolsep{5pc}
\xymatrix{
  H _{G, K, I, W}^j    \ar[r]^{ C_{G, S, R_1}^{P, \, j}    }   \ar[dr]_{ C_{G, S, R_2}^{P, \, j} } 
&   H _{M, S, R_1,  I, W}^{' \, j}   \ar[d]_{ \simeq }^{(\ref{equation-H-M-S-R-1-R-2-isom})} \\
&  H _{M, S, R_2,  I, W}^{' \, j} 
}
\end{equation}
because $C_{G, S, R_1}^{P, \, j}$, $C_{G, S, R_2}^{P, \, j} $ and (\ref{equation-H-M-S-R-1-R-2-isom}) are defined by counit morphisms (which in these cases are equal to trace maps), and by \cite{sga4} XVIII Théorème 2.9, the trace morphism is compatible with composition. 

In \ref{subsection-def-H-M-S-I-W-'} we defined $H _{M, K, I, W}^{' \, j}$. We deduce from (\ref{diagram-CT-R-1-R-2}) a morphism
\begin{equation}   \label{equation-def-C-G-P-K}
C_{G, K}^{P, \, j} : H _{G, K, I, W}^j \rightarrow H _{M, K,  I, W}^{' \, j} .
\end{equation}
which is the composition $H _{G, K, I, W}^j \rightarrow H _{P, K,  I, W}^{' \, j} \rightarrow H _{M, K,  I, W}^{' \, j}$.

\quad

\sssec{}   \label{subsection-functoriality-M-S-1-S-2}
Let $S_1, S_2 \in \mc D$ and $f: S_1 \rightarrow S_2$ be a morphism in $\mc D$. Note that $f$ is $P(\mb A)$-equivariant and it sends orbit to orbit. 
For each $P(\mb A)$-orbit $\beta$ in $S_2$, choose a representative $s^{\beta} \in \beta$ with stabilizer $H^{\beta}$.
If $f^{-1}(\beta)$ is empty, take any $R^{\beta}$ subgroup of finite index in $H_2^{\beta} \cap U(\mb A)$ and normal in $H_2^{\beta}$. 
If $f^{-1}(\beta)$ is non-empty, for every $P(\mb A)$-orbit $\alpha \in f^{-1}(\beta)$, choose a representative $s^{\alpha} \in \alpha$ such that $f(s^{\alpha}) = s^{\beta}$. Let $H_1^{\alpha}$ be the stabilizer of $s^{\alpha}$. Then $H_1^{\alpha} \subset H_2^{\alpha}$. Let $R^{\beta}$ be a subgroup of finite index in $(\cap_{\alpha \in f^{-1}(\beta)} H_1^{\alpha} ) \cap U(\mb A) \subset H_2^{\beta} \cap U(\mb A)$ and normal in $H_2^{\beta}$. 

The morphism $H_1^{\alpha} / R^{\beta} \hookrightarrow H_2^{\beta} / R^{\beta}$ for $\beta = f(\alpha)$ induces a morphism
\begin{equation}     \label{equation-Cht-M-S-1-to-S-2-each-alpha}
\mf q^M_{\alpha}: \Cht_{M, \infty, I, W} / (H_1^{\alpha} / R^{\beta})  \rightarrow \Cht_{M, \infty, I, W} / (H_2^{\beta} / R^{\beta}) 
\end{equation}

Let $R =\big( (R^{\beta})_{\beta \in \{ P(\mb A) \text{-orbits in } S_2 \} } \big)$. Similarly to (\ref{equation-def-Cht-M-S-R-'}), we define $\Cht_{M, S_1, R, I, W}'$ and $\Cht_{M, S_2, R, I, W}'$.
Then (\ref{equation-Cht-M-S-1-to-S-2-each-alpha}) for every orbit $\alpha$ induces a morphism
\begin{equation}
\mf q^M_f: \Cht_{M, S_1, R, I, W}'  \rightarrow \Cht_{M, S_2, R, I, W}' .
\end{equation} 

Similarly to \ref{subsection-adj-pr-G-K-prime-K}, the adjunction morphism $ \Id \rightarrow (\mf q_{f}^M)_* (\mf q_{f}^M)^*$ induces a morphism
\begin{equation}   \label{equation-adj-f-M}
\on{adj}(\mf q_{f}^M): H_{M, S_2, I, W}' \rightarrow H_{M, S_1, I, W}'
\end{equation}
The counit morphism $(\mf q_{f}^M)_! (\mf q_{f}^M)^! \rightarrow \Id$ induces a morphism
\begin{equation}   \label{equation-counit-f-M}
\on{Co}(\mf q_{f}^M): H_{M, S_1, I, W}' \rightarrow H_{M, S_2, I, W}' 
\end{equation}

In the following, we will apply the functoriality to the cases
\begin{itemize}
\item $K' \subset K$, $S_1 = G(\mb A) / K'$, $S_2 = G(\mb A) / K$ and $f$ is the projection $G(\mb A) / K' \twoheadrightarrow G(\mb A) / K$ 

\item $S_1=G(\mb A) / \wt K$, $S_2 = G(\mb A) / g^{-1}\wt Kg$ and $f$ is the isomorphism induced by the right multiplication by $g$: $G(\mb A) / \wt K \isom G(\mb A) / g^{-1}\wt Kg$ 
\end{itemize}

\begin{rem}
In \ref{subsection-functoriality-M-S-1-S-2}, we can also firstly define morphisms of cohomology groups for each orbit $\alpha$:
the adjunction morphism $ \Id \rightarrow (\mf q_{\alpha}^M)_* (\mf q_{\alpha}^M)^*$ induces a morphism
\begin{equation}
\on{adj}(\mf q_{\alpha}^M): H_{M, f(\alpha), I, W}' \rightarrow H_{M, \alpha, I, W}'
\end{equation}
where the orbit $\alpha$ (resp. $f(\alpha)$) is considered as subset of $S_1$ (resp. $S_2$).
The counit morphism $(\mf q_{\alpha}^M)_! (\mf q_{\alpha}^M)^! \rightarrow \Id$ induces a morphism
\begin{equation}
\on{Co}(\mf q_{\alpha}^M): H_{M, \alpha, I, W}' \rightarrow H_{M, f(\alpha), I, W}' 
\end{equation}
Then taking sum over all the orbits, we obtain (\ref{equation-adj-f-M}) and (\ref{equation-counit-f-M}).

Similarly, in \ref{subsection-functoriality-P-S-1-S-2} below, we can firstly prove the statement for cohomology groups orbit by orbit, then take the sum over all the orbits. But the notations would be more complicated. 
\end{rem}

\sssec{}  \label{subsection-functoriality-P-S-1-S-2}
Any $S_1, S_2 \in \mc D$ and $f: S_1 \rightarrow S_2$ morphism in $\mc D$ induce a morphism
$$\mf q^P_f: \Cht_{P, S_1, I, W}' \rightarrow \Cht_{P, S_2, I, W}'$$
The adjunction morphism $ \Id \rightarrow (\mf q_{f}^P)_* (\mf q_{f}^P)^*$ induces a morphism
$$
\on{adj}(\mf q_{f}^P): H_{P, S_2, I, W}' \rightarrow H_{P, S_1, I, W}'
$$
The counit morphism $(\mf q_{f}^P)_! (\mf q_{f}^P)^! \rightarrow \Id$ induces a morphism
$$ 
\on{Co}(\mf q_{f}^P): H_{P, S_1, I, W}' \rightarrow H_{P, S_2, I, W}' 
$$

For each orbit $\alpha$ in $S_1$ with $\beta = f(\alpha)$, let $H_1^{\alpha}$, $H_2^{\beta}$ and $R^{\beta}$ as in \ref{subsection-functoriality-M-S-1-S-2}. We have a Cartesian square:
\begin{equation}   \label{equation-Cht-P-M-H-1-H-2}
\xymatrixrowsep{2pc}
\xymatrixcolsep{4pc}
\xymatrix{
\Cht_{P, \infty, I, W} / H_1^{\alpha}       \ar[d]_{ (\ref{equation-Cht-P-infty-quotient-H-to-Cht-M-infty-quotient-H-R}) }   \ar[r]^{\mf q^P_{\alpha}}
& \Cht_{P, \infty, I, W} / H_2^{\beta}     \ar[d]^{ (\ref{equation-Cht-P-infty-quotient-H-to-Cht-M-infty-quotient-H-R}) }    \\
 \Cht_{M, \infty, I, W} / (H_1^{\alpha} / R^{\beta})  \ar[r]^{\mf q^M_{\alpha}} 
 & \Cht_{M, \infty, I, W} / (H_2^{\beta} / R^{\beta})
}
\end{equation}
Taking union over all the orbits, 
with the notations in \ref{subsection-def-pi-S-R} and \ref{subsection-functoriality-M-S-1-S-2}, we deduce a Cartesian square
\begin{equation}   \label{equation-Cht-P-M-S-1-S-2}
\xymatrixrowsep{2pc}
\xymatrixcolsep{4pc}
\xymatrix{
\Cht_{P, S_1, I, W}'       \ar[d]_{\pi_{S_1, R}}   \ar[r]^{\mf q^P_f}
& \Cht_{P, S_2, I, W}'     \ar[d]^{\pi_{S_2, R}}    \\
 \Cht_{M, S_1, R, I, W}'  \ar[r]^{\mf q^M_f} 
 & \Cht_{M, S_2, R, I, W}'  
}
\end{equation}

(\ref{equation-Cht-P-M-S-1-S-2}) induces a commutative diagram of cohomology groups:
\begin{equation}   \label{equation-H-P-M-S-1-S-2-adj}
\xymatrixrowsep{2pc}
\xymatrixcolsep{4pc}
\xymatrix{
H_{P, S_2, I, W}'       \ar[d]_{ (\ref{equatino-H-P-to-H-M})  }   \ar[r]^{\on{adj}(\mf q^P_f) } 
& H_{P, S_1, I, W}'     \ar[d]^{ (\ref{equatino-H-P-to-H-M})  }   \\
 H_{M, S_2, I, W}'  \ar[r]^{\on{adj}(\mf q^M_f) } 
 & H_{M, S_1, I, W}'  
}
\end{equation}
because (\ref{equation-Cht-P-M-S-1-S-2}) is Cartesian, (\ref{equatino-H-P-to-H-M}) is given by a counit morphism (equal to the trace morphism), and by \cite{sga4} XVIII Théorème 2.9, the trace morphism commutes with base change.

(\ref{equation-Cht-P-M-S-1-S-2}) induces a commutative diagram of cohomology groups:
\begin{equation}   \label{equation-H-P-M-S-1-S-2-co}
\xymatrixrowsep{2pc}
\xymatrixcolsep{4pc}
\xymatrix{
H_{P, S_1, I, W}'       \ar[d]_{ (\ref{equatino-H-P-to-H-M})  }   \ar[r]^{\on{Co}(\mf q^P_f) }
& H_{P, S_2, I, W}'     \ar[d]^{ (\ref{equatino-H-P-to-H-M})  }    \\
 H_{M, S_1, I, W}'  \ar[r]^{\on{Co}(\mf q^M_f) } 
 & H_{M, S_2, I, W}'  
}
\end{equation}
because by \cite{sga4} XVIII Théorème 2.9, the trace morphism is compatible with composition.

\begin{rem}
When $S_1 = G(\mb A) / K_{N_1}$ and $S_2 = G(\mb A) / K_{N_2}$ with $N_1 \supset N_2$, we have the projection $f: G(\mb A) / K_{N_1} \twoheadrightarrow G(\mb A) / K_{N_2}$. We have $\Cht_{M, N_1, I, W}' = \Cht_{M, S_1, R_1, I, W}'$ with $R_1^{\alpha} = K_{U, N_1}$ for each $P(\mb A)$-orbit $\alpha$ in $S_1$ and $\Cht_{M, N_2, I, W}' = \Cht_{M, S_2, R_2, I, W}'$ with $R_2^{\beta} = K_{U, N_2}$ for each $P(\mb A)$-orbit $\beta$ in $S_2$. Note that $R_1^{\alpha} \neq R_2^{f(\alpha)}$, thus the commutative diagram 
\begin{equation}  \label{equation-Cht-P-M-N-1-N-2-'}
\xymatrixrowsep{2pc}
\xymatrixcolsep{4pc}
\xymatrix{
\Cht_{P, N_1, I, W}'       \ar[d]_{\pi'}   \ar[r]
& \Cht_{P, N_2, I, W}'     \ar[d]^{\pi'}    \\
 \Cht_{M, N_1, I, W}'  \ar[r]
 & \Cht_{M, N_2, I, W}'  
}
\end{equation}
does NOT coincides with diagram (\ref{equation-Cht-P-M-S-1-S-2}). In particular, diagram (\ref{equation-Cht-P-M-N-1-N-2-'}) is not Cartesian (the morphism from $\Cht_{P, N_1, I, W}'$ to the fiber product is finite étale of degree $\sharp(K_{U, N_2} / K_{U, N_1})$ which is a power of $q$).
\end{rem}

\sssec{}   \label{subsection-functoriality-K'-K-M-adj}
Let $K' \subset K$ be two compact open subgroups of $G(\mb O)$.
Applying \ref{subsection-functoriality-M-S-1-S-2} 
to $S_1 = G(\mb A) / K'$, $S_2 = G(\mb A) / K$ and the projection $f: G(\mb A) / K' \twoheadrightarrow G(\mb A) / K$, we deduce a finite étale morphism (denoted by $\mf q_f^M$ in $loc.cit.$)
\begin{equation*}   \label{equation-pr-M-K'-K}
\on{pr}^M_{K', K}: \Cht_{M, S_1, R, I, W}' \rightarrow \Cht_{M, S_2, R, I, W}' ,
\end{equation*}
where $R$ is defined in \ref{subsection-functoriality-M-S-1-S-2}.
The adjunction morphism $\on{adj}(\on{pr}^M_{K', K}): \Id \rightarrow (\on{pr}^M_{K', K})_* (\on{pr}^M_{K', K})^*$ induces 
$$\on{adj}(\on{pr}^M_{K', K}): H _{M, K, I, W}^{' \, j}  \rightarrow H _{M, K', I, W}^{' \, j} .$$
The counit morphism $\on{Co}(\on{pr}^M_{K', K}): (\on{pr}^M_{K', K})_! (\on{pr}^M_{K', K})^! \rightarrow \Id$ induces
$$\on{Co}(\on{pr}^M_{K', K}): H _{M, K', I, W}^{' \, j}  \rightarrow  H _{M, K, I, W}^{' \, j} .$$

\begin{lem}   \label{lem-CT-commutes-with-level-augement}
For $K' \subset K$ as in \ref{subsection-functoriality-K'-K-M-adj}, the following diagram of cohomology groups commutes:
\begin{equation}     \label{equation-CT-commutes-with-adj}
\xymatrixrowsep{2pc}
\xymatrixcolsep{5pc}
\xymatrix{
  H _{G, K, I, W}^j    \ar[r]^{\on{adj}(\on{pr}^G_{K', K})}  \ar[d]_{ C_{G, K}^{P, \, j}} 
&   H _{G, K', I, W}^j  \ar[d]^{ C_{G, K'}^{P, \, j}}  \\
 H _{M, K, I, W}^{' \, j}  \ar[r]^{\on{adj}(\on{pr}^M_{K', K})}
&  H _{M, K', I, W}^{' \, j} 
}
\end{equation}
\end{lem}

\dem
(1) By (\ref{equation-Cht-P-K-'-is-fiber-prod}), we have a Cartesian square
\begin{equation}          \label{diagram-Cht-G-K-K'-P-K-K'}
\xymatrixrowsep{2pc}
\xymatrixcolsep{4pc}
\xymatrix{
\Cht_{G, K, I,  W}    
& \Cht_{G, K', I, W}  \ar[l]_{\on{pr}^G_{K', K}}  \\
\Cht_{P, K, I, W}'     \ar[u]^{i_K}  
& \Cht_{P, K', I, W}'   \ar[u]_{i_{K'}}    \ar[l]^{\on{pr}^P_{K', K}} 
}
\end{equation}
Since adjunction morphism is compatible with composition, we deduce that the following diagram is commutative:
\begin{equation*}     
\xymatrixrowsep{2pc}
\xymatrixcolsep{5pc}
\xymatrix{
  H _{G, K, I, W}^j    \ar[r]^{\on{adj}(\on{pr}^G_{K', K})}  \ar[d]_{\on{adj}(i_K)} 
&   H _{G, K', I, W}^j  \ar[d]^{\on{adj}(i_{K'})}  \\
 H _{P, K, I, W}^{' \, j}  \ar[r]^{\on{adj}(\on{pr}^P_{K', K})}
&  H _{P, K', I, W}^{' \, j} 
}
\end{equation*}

(2) Applying \ref{subsection-functoriality-P-S-1-S-2} to $S_1 = G(\mb A) / K'$, $S_2 = G(\mb A) / K$ and the projection $f: G(\mb A) / K' \twoheadrightarrow G(\mb A) / K$, we deduce from (\ref{equation-H-P-M-S-1-S-2-adj}) that the following diagram is commutative:
\begin{equation*}     
\xymatrixrowsep{2pc}
\xymatrixcolsep{5pc}
\xymatrix{
  H _{P, K, I, W}^{' \, j}    \ar[r]^{\on{adj}(\on{pr}^P_{K', K})}  \ar[d]_{(\ref{equatino-H-P-to-H-M}) }
&   H _{P, K', I, W}^{' \, j}  \ar[d]^{(\ref{equatino-H-P-to-H-M}) }  \\
 H _{M, K, I, W}^{' \, j}  \ar[r]^{\on{adj}(\on{pr}^M_{K', K})}
&  H _{M, K', I, W}^{' \, j} 
}
\end{equation*}
\cqfd

\begin{lem}    \label{lem-CT-commutes-with-level-dimunit}
For $K' \subset K$ as in \ref{subsection-functoriality-K'-K-M-adj}, the following diagram of cohomology groups commutes:
$$
\xymatrixrowsep{2pc}
\xymatrixcolsep{5pc}
\xymatrix{
  H _{G, K', I, W}^j    \ar[r]^{\on{Co}(\on{pr}^G_{K', K})}  \ar[d]_{ C_{G, K'}^{P, \, j}} 
&   H _{G, K, I, W}^j  \ar[d]^{  C_{G, K}^{P, \, j}}  \\
 H _{M, K', I, W}^{' \, j}  \ar[r]^{\on{Co}(\on{pr}^M_{K', K})}
&  H _{M, K, I, W}^{' \, j} 
}
$$
\end{lem}
\dem
(1) By \cite{sga4} XVIII Théorème 2.9, the trace morphism commutes with base change. Since (\ref{diagram-Cht-G-K-K'-P-K-K'}) is Cartesian, 
we deduce that the following diagram is commutative:
$$
\xymatrixrowsep{2pc}
\xymatrixcolsep{4pc}
\xymatrix{
  H _{G, K', I, W}^j    \ar[r]^{\on{Co}(\on{pr}^G_{K', K})}  \ar[d]_{\on{adj}(i_{K'})} 
&   H _{G, K, I, W}^j  \ar[d]^{\on{adj}(i_K)} \\
 H _{P, K', I, W}^{' \, j}  \ar[r]^{\on{Co}(\on{pr}^P_{K', K})}
&  H _{P, K, I, W}^{' \, j} 
}
$$
(2) 
Applying \ref{subsection-functoriality-P-S-1-S-2} to $S_1 = G(\mb A) / K'$, $S_2 = G(\mb A) / K$ and the projection $f: G(\mb A) / K' \twoheadrightarrow G(\mb A) / K$, we deduce from (\ref{equation-H-P-M-S-1-S-2-co}) that the following diagram is commutative:
$$
\xymatrixrowsep{2pc}
\xymatrixcolsep{4pc}
\xymatrix{
  H _{P, K', I, W}^{' \, j}    \ar[r]^{\on{Co}(\on{pr}^P_{K', K})}  \ar[d]_{(\ref{equatino-H-P-to-H-M}) } 
&   H _{P, K, I, W}^{' \, j}  \ar[d]^{(\ref{equatino-H-P-to-H-M}) } \\
 H _{M, K', I, W}^{' \, j}  \ar[r]^{\on{Co}(\on{pr}^M_{K', K})}
&  H _{M, K, I, W}^{' \, j} 
}
$$
\cqfd

\subsection{Compatibility of constant term morphisms and actions of Hecke algebras}

We first recall the action of the local Hecke algebras. The goal of this subsection is Lemma \ref{lem-TC-commute-with-Hecke-with-niveau} and Lemma \ref{lem-TC-commute-with-Hecke-without-niveau}.

\sssec{}    \label{subsection-Hecke-acts-on-H-G-infty}
Let $v$ be a place in $X$. 
Let $g \in G(F_v)$. By \ref{subsection-G-A-acts-on-Cht-G-infty}, the right action of $g$ induces an isomorphism 
\begin{equation}   \label{equation-g-acts-on-Cht-K-infty}
 \Cht_{G, \infty , I, W} \isom \Cht_{G, \infty, I, W} .
\end{equation}

Let $\wt K \subset  G(\mb O)$ be a compact open subgroup such that $g^{-1}\wt Kg \subset G(\mb O)$. 
The isomorphism (\ref{equation-g-acts-on-Cht-K-infty}) is $\wt K$-equivariant, where $k \in \wt K$ acts on the second stack by $g^{-1}kg$. 
It induces an isomorphism $$ \Cht_{G, \infty , I, W} / \wt K \isom \Cht_{G, \infty , I, W} / g^{-1} \wt Kg,$$ 
i.e. $ \Cht_{G, \wt K, I, W} \isom \Cht_{G, g^{-1} \wt Kg, I, W} .$ 
It induces (by adjunction) an isomorphism of cohomology groups 
\begin{equation}   \label{equation-adj-g-H-G-wt-K}
\on{adj}(g): H _{G, g^{-1} \wt K g, I, W}^j \isom H _{G, \wt K, I, W}^j .
\end{equation}



\sssec{}   \label{subsection-Hecke-G-acts-on-H-G}  
Let $K = K^v  K_v  \subset G(\mb O^v) G(\mc O_v) =  G(\mb O) $ be an open compact subgroup.
Let $h =  {\bf 1}_{K_v g K_v} \in C_c(K_v \backslash G(F_v) / K_v, \Ql)$ be the characteristic function of $K_v g K_v$ for some $g \in G(F_v)$. 
The action of $h$ on $H _{G, K, I, W}^j $ is given by the following composition of morphisms
\begin{equation}     \label{equation-T-h-acts-on-H-G}
T(h): H _{G, K, I, W}^j \xrightarrow{\on{adj}} H _{G, K \cap g^{-1} K g, I, W}^j \underset{\sim}{\xrightarrow{  \on{adj}(g)  }} H _{G, g Kg^{-1} \cap K, I, W}^j  \xrightarrow{\on{Co}} H _{G, K, I, W}^j ,
\end{equation}
where $\on{adj}=\on{adj}(\on{pr}^G_{K \cap g^{-1}Kg, K})$ and $\on{Co}=\on{Co}(\pr^G_{gKg^{-1} \cap K, K})$, the isomorphism $\on{adj}(g)$ is induced by (\ref{equation-adj-g-H-G-wt-K}) applied to $\wt K = g Kg^{-1} \cap K$.
Note that (\ref{equation-T-h-acts-on-H-G}) depends only on the class $K_v g K_v$ of $g$ in $G(F_v)$. 
The action of $T(h)$ is equivalent to the one constructed by Hecke correspondence (see \cite{vincent} 2.20 and 4.4).

\sssec{}   \label{subsection-functoriality-K-gKg-M}
Let $\wt K$ and $g$ as in \ref{subsection-Hecke-acts-on-H-G-infty}.
The right action of $g$ (by right multiplication by $g$) on $G(\mb A)$ induces an isomorphism
\begin{equation}   \label{equation-G-A-K-to-G-A-gKg}
G(\mb A) / \wt K \isom G(\mb A) / g^{-1}\wt Kg
\end{equation}
Applying \ref{subsection-functoriality-M-S-1-S-2} to $S_1=G(\mb A) / \wt K$, $S_2 = G(\mb A) / g^{-1}\wt Kg$ and the isomorphism (\ref{equation-G-A-K-to-G-A-gKg}), we deduce an isomorphism of cohomology groups
\begin{equation}  \label{equation-adj-g-H-M-wt-K}
\on{adj}(g): H_{M, g^{-1}\wt Kg, I, W}^{' \, j} \isom H_{M, \wt K, I, W}^{' \, j}
\end{equation}

\sssec{}   \label{subsection-Hecke-G-acts-on-H-M-'}  
Let $K$ and $h$ as in \ref{subsection-Hecke-G-acts-on-H-G}.  
The action of $h$ on $H _{M, K, I, W}^{' \, j} $ is given by the following composition of morphisms
\begin{equation}     \label{equation-T-h-acts-on-H-M-'}
T(h): H _{M, K, I, W}^{' \, j} \xrightarrow{\on{adj}} H _{M, K \cap g^{-1} K g, I, W}^{' \, j} \underset{\sim}{\xrightarrow{  \on{adj}(g)  }} H _{M, g Kg^{-1} \cap K, I, W}^{' \, j}  \xrightarrow{\on{Co}} H _{M, K, I, W}^{' \, j} ,
\end{equation}
where $\on{adj}=\on{adj}(\on{pr}^M_{K \cap g^{-1}Kg, K})$ and $\on{Co}=\on{Co}(\pr^M_{gKg^{-1} \cap K, K})$, the isomorphism $\on{adj}(g)$ is induced by (\ref{equation-adj-g-H-M-wt-K}) applied to $\wt K = g Kg^{-1} \cap K$. Note that $\wt K$ may not be normal in $G(\mb O)$.
Note that (\ref{equation-T-h-acts-on-H-M-'}) depends only on the class $K_v g K_v$ of $g$ in $G(F_v)$.

\begin{lem}    \label{lem-TC-commute-avec-Hecke-infty}
Let $\wt K$ and $g$ as in \ref{subsection-Hecke-acts-on-H-G-infty}. The following diagram of cohomology groups commutes:
\begin{equation}   \label{equation-TC-commute-avec-Hecke-infity}
\xymatrixrowsep{2pc}
\xymatrixcolsep{4pc}
\xymatrix{
  H _{G, g^{-1} \wt K g, I, W}^j   \ar[r]^{\on{adj}(g)}_{\sim}  \ar[d]_{ C_{G, g^{-1} \wt K g}^{P, \, j}} 
&   H _{G, \wt K , I, W}^j   \ar[d]^{ C_{G, \wt K}^{P, \, j}}  \\
 H _{M, g^{-1} \wt K g, I, W}^{' \, j}  \ar[r]^{\on{adj}(g)}_{\sim}
&  H _{M, \wt K, I, W}^{' \, j}  .
}
\end{equation}
\end{lem}
\dem
(1) Since the isomorphism (\ref{equation-Cht-P-infty-'-to-fiber-prod}) is $G(\mb O)$-equivariant, we deduce a Cartesian square:
\begin{equation*} 
\xymatrixrowsep{2pc}
\xymatrixcolsep{2pc}
\xymatrix{
\Cht_{P, I, W} \underset{\Cht_{G, I, W}} \times \Cht_{G, \infty, I, W} / \wt K  \ar[r]^{g \quad}_{\simeq \quad}  
& \Cht_{P, I, W} \underset{\Cht_{G, I, W}} \times \Cht_{G, \infty, I, W} / g^{-1} \wt K g \\
\Cht_{P, \infty, I,  W} \overset{P(\mb A) } \times G(\mb A) / \wt K  \ar[r]^{g \quad}_{\simeq \quad}    \ar[u]^{\simeq}  
& \Cht_{P, \infty, I,  W} \overset{P(\mb A) } \times G(\mb A) / g^{-1} \wt K g   \ar[u]^{\simeq}  
}
\end{equation*}
We deduce a Cartesian square:
\begin{equation*} 
\xymatrixrowsep{2pc}
\xymatrixcolsep{2pc}
\xymatrix{
\Cht_{G, \wt K, I, W}   \ar[r]^{g \quad}_{\simeq \quad}  
& \Cht_{G, g^{-1} \wt K g , I, W}  \\
\Cht_{P, \wt K, I, W}'   \ar[r]^{g \quad}_{\simeq \quad}    \ar[u]^{i_{\wt K}'}  
& \Cht_{P, g^{-1} \wt K g, I, W}'  \ar[u]_{i_{g^{-1} \wt K g}'}  
}
\end{equation*}
It induces a commutative diagram
$$
\xymatrixrowsep{1.5pc}
\xymatrixcolsep{4pc}
\xymatrix{
  H _{G, g^{-1} \wt K g, I, W}^j   \ar[r]^{\on{adj}(g)}_{\sim}  \ar[d]_{\on{adj}(i_{g^{-1} \wt K g}')}
&   H _{G, \wt K , I, W}^j   \ar[d]^{\on{adj}(i_{\wt K }')}  \\
 H _{P, g^{-1} \wt K g, I, W}^{' \, j}  \ar[r]^{\on{adj}(g)}_{\sim}
&  H _{P, \wt K, I, W}^{' \, j}  .
}
$$
(2) Applying \ref{subsection-functoriality-P-S-1-S-2} to $S_1=G(\mb A) / \wt K$, $S_2 = G(\mb A) / g^{-1}\wt Kg$ and $f$ the isomorphism (\ref{equation-G-A-K-to-G-A-gKg}), we deduce from (\ref{equation-H-P-M-S-1-S-2-co}) a commutative diagram
$$
\xymatrixrowsep{1.5pc}
\xymatrixcolsep{4pc}
\xymatrix{
  H _{P, g^{-1} \wt K g, I, W}^{' \, j}   \ar[r]^{\on{adj}(g)}_{\sim}  \ar[d]_{(\ref{equatino-H-P-to-H-M}) }
&   H _{P, \wt K , I, W}^{' \, j}   \ar[d]^{(\ref{equatino-H-P-to-H-M}) }  \\
 H _{M, g^{-1} \wt K g, I, W}^{' \, j}  \ar[r]^{\on{adj}(g)}_{\sim}
&  H _{M, \wt K, I, W}^{' \, j}  .
}
$$
\cqfd

\begin{lem}   \label{lem-TC-commute-with-Hecke-with-niveau}
For any place $v$ of $X$, any $K$ and $h \in C_{c}( K_v \backslash G(F_v)/ K_v, \Ql) $ as in \ref{subsection-Hecke-G-acts-on-H-G}, the following diagram of cohomology groups commutes:
\begin{equation}   \label{diagram-TC-commute-avec-Hecke}
\xymatrixrowsep{2pc}
\xymatrixcolsep{2pc}
\xymatrix{
  H _{G, K, I, W}^{j}    \ar[r]^{T(h)}  \ar[d]^{ C_{G, K}^{P, \, j}} 
&   H _{G, K, I, W}^{j}  \ar[d]^{ C_{G, K}^{P, \, j}}  \\
 H _{M, K, I, W}^{' \, j}  \ar[r]^{T(h)}
&  H _{M, K, I, W}^{' \, j} 
}
\end{equation}
where the horizontal morphisms are defined in \ref{subsection-Hecke-G-acts-on-H-G} and \ref{subsection-Hecke-G-acts-on-H-M-'}, the vertical morphism are the constant term morphisms defined in (\ref{equation-def-C-G-P-K}).
\end{lem}
\dem
By Lemma \ref{lem-CT-commutes-with-level-augement}, Lemma \ref{lem-TC-commute-avec-Hecke-infty} and Lemma \ref{lem-CT-commutes-with-level-dimunit}.
\cqfd

%

\quad

\sssec{}
From now on let $N \subset X$ be a closed subscheme and $v$ be a place in $X \sm N$. 
We have the (unnormalized) Satake transform:
\begin{equation}   \label{equation-TC-pour-op-de-Hecke}
\begin{aligned}
C_{c}(G(\mc O_v)\backslash G(F_v)/G(\mc O_v), \Ql) & \hookrightarrow C_{c}(M(\mc O_v)\backslash M(F_v) / M(\mc O_v), \Ql) \\
h \quad & \mapsto \quad  h^M: m \mapsto  \sum_{ U( F_v) / U(\mc O_v) }   h(mu) .
\end{aligned}
\end{equation}

\sssec{}
We have $K_{M, N} = K_{M, N}^vK_{M, N, v} \subset M(\mb O^v)M(\mc O_v)$.
For any $K_{M, v} \subset M(\mc O_v)$ open compact subgroup, we have $K_{M, N}^v K_{M, v} \subset M(\mb O^v)M(\mc O_v)$.
We define $H _{M, K_{M, N}^v K_{M, v}, I, W}^j$ as in Definition \ref{def-H-M-j} (replacing $\Cht_{M, N, I, W}$ by $\Cht_{M, K_{M, N}^v K_{M, v}, I, W}$).
We define $$\varinjlim_{K_{M, v}} H _{M, K_{M, N}^vK_{M, v}, I, W}^{j} . $$
As in \ref{subsection-Hecke-acts-on-H-G-infty} (by replacing $G$ by $M$), for any $m \in M(F_v)$ and $K_{M, v}$ such that $m^{-1}K_{M, v}m \subset M(\mc O_v)$, we have an isomorphism $H _{M, m^{-1}K_{M, N}^vK_{M, v}m, I, W}^{j} \isom H _{M, K_{M, N}^vK_{M, v}, I, W}^{j}$. Taking limit on $K_{M, v}$, we deduce an action of $M(F_v)$ on $\varinjlim_{K_{M, v}} H _{M, K_{M, N}^vK_{M, v}, I, W}^{j} . $

We have $K_{N} = K_{N}^vK_{N, v} \subset G(\mb O^v)G(\mc O_v)$.
For any $K_v \subset G(\mc O_v)$ open compact subgroup, we have $K_N^v K_v \subset G(\mb O^v) G(\mc O_v)$. Applying \ref{subsection-def-H-M-S-I-W-'} to $S = G(\mb A) / K_N^v K_{v} $, we define $H_{M, K_N^v K_{v}, I, W} ^{' \, j}$. 
We define $$\varinjlim_{K_v} H_{M, K_N^v K_v, I, W}^{' \, j} . $$
Note that $v$ is a place in $X \sm N$, so $K_{N, v} = G(\mc O_v)$ and $K_{M, N, v} = M(\mc O_v)$. We have 
$$H_{M, K_N^v G(\mc O_v), I, W}^{' \, j}  = H_{M, N, I, W}^{' \, j} = H_{M, N, I, W}^j \overset{ P(\mc O_N) }{ \times}  G(\mc O_N) = H_{M, K_{M, N}^vM(\mc O_v), I, W}^j \overset{ P(\mc O_N) }{ \times}  G(\mc O_N) ,$$ where $H_{M, N, I, W}^{' \, j}$ is defined in Definition \ref{def-H-M-prime-j}. 
We deduce
\begin{equation}
\varinjlim_{K_v} H_{M, K_N^v  K_v, I, W}^{' \, j}  = \on{Ind}_{P(F_v)}^{G(F_v)} \big( (\varinjlim_{K_{M, v}} H _{M, K_{M, N}^vK_{M, v}, I, W}^{j}) \overset{ P(\mc O_N) }{ \times}  G(\mc O_N)   \big),
\end{equation}
where $\on{Ind}_{P(F_v)}^{G(F_v)} $ is the (unnormalized) parabolic induction.

%

\sssec{}
Let $V$ be a $\Ql$-vector space equipped with a continuous action of $M(F_v)$, denoted by $\sigma: M(F_v) \rightarrow GL(V)$.
We recall that
$$\on{Ind}_{P(F_v)}^{G(F_v)} V = \{ f: G(F_v) \rightarrow V \text{ continuous}, f(pg)= \sigma(p)f(g), p \in P(F_v), g \in G(F_v) \}.$$
We have a morphism
\begin{equation}    \label{equation-Ind-V-to-V}
(\on{Ind}_{P(F_v)}^{G(F_v)} V )^{G(\mc O_v)}  \rightarrow  V^{M(\mc O_v)} : \quad
f \mapsto  f(1)
\end{equation}

\begin{lem}   \label{lem-action-h-equal-action-h-M}
Morphism (\ref{equation-Ind-V-to-V}) is an isomorphism.
Moreover, for $h \in C_{c}(G(\mc O_v)\backslash G(F_v)/G(\mc O_v), \Ql)$, the action of $T(h)$ on $(\on{Ind}_{P(F_v)}^{G(F_v)} V )^{G(\mc O_v)}$ coincides with the action of $T(h^M)$ on $ V^{M(\mc O_v)} $. 
\end{lem}
\dem
Morphism (\ref{equation-Ind-V-to-V}) admits an inverse $f(1) \mapsto f$ given by
$$
f(x) = f(x_P x_K) = \sigma(x_P) f(x_K) = \sigma(x_P) f(1) 
$$
where $x = x_P x_K \in G(F_v) = P(F_v) G(\mc O_v)$. Thus $ (\on{Ind}_{P(F_v)}^{G(F_v)} V )^{G(\mc O_v)} = V^{M(\mc O_v)} .$

Moreover, for $g = g_P g_K  \in G(F_v)$ and $f \in (\on{Ind}_{P(F_v)}^{G(F_v)} V)^{G(\mc O_v)}$, we have
\begin{equation}  \label{equation-g-acts-on-Ind-P-G}
g f(1) = f(g_P g_K ) = \sigma(g_P) f(g_K) = \sigma(g_P) f(1).
\end{equation}
Note that $G(F_v) = M(F_v) U(F_v) G(\mc O_v)$. Denote by $dg$ (resp. $dm, du, dk$) the Haar measure on $G(F_v)$ (resp $M(F_v), U(F_v), G(\mc O_v)$) such that the volume of $G(O_v)$ (resp. $M(\mc O_v), U(\mc O_v), G(\mc O_v)$) is $1$. We have $dg = dm du dk$. Taking the integral over $G(F_v)$ of the product by $h(g)$ of (\ref{equation-g-acts-on-Ind-P-G}), we deduce that the action of $T(h)$ on $(\on{Ind}_{P(F_v)}^{G(F_v)} V )^{G(\mc O_v)}$ coincides with the action of $T(h^M)$ on $ V^{M(\mc O_v)} $. 
\cqfd

\sssec{}
Let $V = (\varinjlim_{K_{M, v}} H _{M, K_{M, N}^vK_{M, v}, I, W}^{j}) \overset{ P(\mc O_N) }{ \times}  G(\mc O_N) $,
we have
\begin{equation*}    
V^{M(\mc O_v)} = H_{M, K_{M, N}^vM(\mc O_v), I, W}^j \overset{ P(\mc O_N) }{ \times}  G(\mc O_N) = H_{M, N, I, W}^{' \, j} ;
\end{equation*}
$$(\on{Ind}_{P(F_v)}^{G(F_v)} V )^{G(\mc O_v)} = ( \varinjlim_{K_v} H_{M, K_N^v  K_v, I, W}^{' \, j}   )^{G(\mc O_v)} = H_{M, K_N^v  G(\mc O_v), I, W}^{' \, j}  =  H_{M, N, I, W}^{' \, j} $$
By Lemma \ref{lem-action-h-equal-action-h-M}, the action of $T(h)$ on $H_{M, N, I, W}^{' \, j}$ (defined in (\ref{equation-T-h-acts-on-H-M-'})) coincides with the action of $T(h^M)$ on $H_{M, N, I, W}^{' \, j}$ (induced by the action of $T(h^M)$ on $H_{M, N, I, W}^{j}$).
Combining this fact and Lemma \ref{lem-TC-commute-with-Hecke-with-niveau}, we deduce 

\begin{lem}    \label{lem-TC-commute-with-Hecke-without-niveau}
For any place $v$ of $X \sm N$ and any $h \in C_{c}( G(\mc O_v) \backslash G(F_v)/ G(\mc O_v), \Ql) $, the following diagram of cohomology groups is commutative:
\begin{equation}   \label{diagram-TC-commute-avec-Hecke}
\xymatrixrowsep{2pc}
\xymatrixcolsep{4pc}
\xymatrix{
  H _{G, N, I, W}^j  \ar[r]^{T(h)}  \ar[d]^{ C_{G, N}^{P, \, j}} 
&   H _{G, N, I, W}^j  \ar[d]^{ C_{G, N}^{P, \, j}}  \\
 H _{M, N, I, W}^{' \, j}  \ar[r]^{T(h^M)}
&  H _{M, N, I, W}^{' \, j} 
}
\end{equation}
where the vertical morphisms are the constant term morphisms defined in Definition \ref{def-CT-cohomology}.
\cqfd
\end{lem}

\begin{rem}
For a direct proof of Lemma \ref{lem-TC-commute-with-Hecke-without-niveau}, see \cite{these} Lemme 8.1.1.
\end{rem}

\begin{rem}
We could normalize the constant term morphism $C_{G, N}^{P, \, j}$ and the Satake transform (\ref{equation-TC-pour-op-de-Hecke}) by $\delta^{\frac{1}{2}}$ as usual, where $\delta$ is the modular function of $P(F_v)$. But we do not need this normalization in this paper.
\end{rem}

\begin{rem}
When $I = \emptyset$ and $W = \bf 1$, $S = G(\mb A) / K$, we have $H_{M, S, R, I, W}^{' \, 0} \subset C(M(F)U(\mb A) \backslash G(\mb A) / K \Xi, \Ql)$. In (\ref{equation-def-C-G-P-S-R}), we defined $C_{G, S, R, I, W}^{P, \, 0}$. Commutative diagram (\ref{diagram-CT-R-1-R-2}) implies that for a given Haar measure $du$ on $U(\mb A)$, $(\int_R du) \cdot C_{G, S, R, I, W}^{P, \, 0} $ is independant on $R$. This identifies $C_{G, S, I, W}^{P, \, 0}$ with the classical constant term morphism (\ref{equation-classical-CT}) associated to $du$. 

\end{rem}

\subsection{Cuspidal cohomology and rational Hecke-finite cohomology}

\begin{defi}  \label{def-H-Hf-rat}
We define
$$H_{G, N, I, W}^{j, \; \on{Hf-rat}}:=\{  c \in H_{G, N, I, W}^j , \; \dim_{\Ql} C_{c}(K_{N}\backslash G(\mb A)/K_{N}, \Ql) \cdot c < + \infty  \}.$$
\end{defi}

Proposition \ref{prop-cusp-egal-Hfrat-section-6} will follow from Lemma \ref{lem-cusp-inclus-dans-Hf-rat} and Lemma \ref{lem-Hf-rat-inclus-dans-cusp} below.

\begin{lem}   \label{lem-cusp-inclus-dans-Hf-rat}
We have an inclusion
\begin{equation}
H_{G, N, I, W}^{j, \; \on{cusp}} \subset H_{G, N, I, W}^{j, \; \mr{Hf-rat}} .
\end{equation}
\end{lem}

\dem
By Theorem \ref{thm-cusp-dim-fini-second}, the $\Ql$-vector space $H_{G, N, I, W}^{j, \, \on{cusp}}$ has finite dimension. 
By Lemma \ref{lem-TC-commute-with-Hecke-with-niveau}, it is stable under the action of the Hecke algebra $C_{c}(K_{N}\backslash G(\mb A)/K_{N}, \Ql)$. 
We conclude by Definition \ref{def-H-Hf-rat}.
\cqfd

\begin{lem}   \label{lem-Hf-rat-inclus-dans-cusp}
We have an inclusion 
\begin{equation}
H_{G, N, I, W}^{j, \; \on{cusp}} \supset H_{G, N, I, W}^{j, \; \mr{Hf-rat}} .
\end{equation}
\end{lem}

The proof of Lemma \ref{lem-Hf-rat-inclus-dans-cusp} will use the fact that any non zero image of a constant term morphism $C_P^{G, j}$ is supported on the components $H_M^{' \, j, \, \nu}$ indexed by $\nu$ in a translated cone in $\wh{\Lambda}_{Z_G/Z_M}^{+, \Q}$.
The proof will also need the following lemma, which is for example a consequence of the Satake isomorphism:
\begin{lem} \label{lem-Hecke-M-fini-sur-Hecke-G}
Under the Satake transformation (\ref{equation-TC-pour-op-de-Hecke}), the algebra $C_{c}(M(O_v)\backslash M(F_v) / M(O_v), \Ql)$ is finite over $C_{c}(G(O_v)\backslash G(F_v)/G(O_v), \Ql)$.
\cqfd
\end{lem}

\noindent {\bf Proof of Lemma \ref{lem-Hf-rat-inclus-dans-cusp}.}  
Let $a \in H_{G, N, I, W}^{\mr{Hf-rat}}$.
We argue by contradiction. Suppose that $a \notin H_{G, N, I, W}^{j, \; \on{cusp}}$. Then there exists a maximal parabolic subgroup $P$ such that $C_G^{P, \; j}(a) \neq 0$. We denote by $M$ the Levi quotient of $P$.
Let $v$ be a place in $X \sm N$. 

(1) On the one hand, by Definition \ref{def-H-Hf-rat}, the $\Ql$-vector subspace $C_{c}(G(O_v)\backslash G(F_v)/G(O_v), \Ql) \cdot a$ has finite dimension. Then Lemma \ref{lem-TC-commute-with-Hecke-without-niveau} applied to $K=K_N$ and Lemma \ref{lem-Hecke-M-fini-sur-Hecke-G} imply that the $\Ql$-vector space $C_{c}(M(O_v)\backslash M(F_v)/M(O_v), \Ql) \cdot C^{P, \; j}_G(a)$ has finite dimension.

(2) On the other hand, since $a \in H_{G, N, I, W}^j$, there exists $\mu \in \wh{\Lambda}_{G^{ad}}^{+, \Q}$ such that $a \in \on{Im} (H_{G, N, I, W}^{j, \; \leq \mu}  \rightarrow  H_{G, N, I, W}^j )$. We deduce from (\ref{equation-TC-mu-1-mu-2}) that $C_G^{P, \; j}(a)$ is the image of an element $a^0 \in H_{M, N, I, W}^{' \; j, \; \leq \mu}$ in $ H_{M, N, I, W}^{' \; j}$. By \ref{subsection-H-M-leq-mu-support-on-cone}, $a^0$ is supported on the components $H_{M, N, I, W}^{' \; j, \, \nu }$ of $H_{M, N, I, W}^{' \; j }$ indexed by $\nu$ in the translated cone $\wh{\Lambda}_{Z_M/ Z_G}^{\mu} \subset \wh{\Lambda}_{Z_M/ Z_G}^{\Q}$. So is $C_G^{P, \; j}(a)$. 


Let $g \in Z_M(F_v)$ such that $g \notin Z_M(O_v)Z_G(F_v)$.
We denote by $\xi(g)$ the image of $g$ by the composition of morphisms
$$Z_M(F_v) \rightarrow Z_M(\mb A) \rightarrow \Bun_{M}(\Fq) \rightarrow \Bun_M \xrightarrow{(\ref{equation-Bun-M-deg-to-Lambda-Z-M-Z-G})} \wh{\Lambda}_{Z_M/ Z_G}^{\Q}.$$
The choice of $g$ implies that $\xi(g) \neq 0$. Note that $P$ is maximal, so $\wh{\Lambda}_{Z_M/ Z_G}^{\Q} \simeq \Q$.
For all $\nu \in \wh{\Lambda}_{Z_M/ Z_G}^{\Q}$, the action of $g$ on $\Cht_M /\Xi $ induces an isomorphism  $g: \Cht_M^{\nu}/\Xi \isom \Cht_M^{\nu + \xi(g)} / \Xi$ (the inverse is induced by $g^{-1}$). We denote by $T(g) \in C_{c}(M(O_v)\backslash M(F_v) / M(\mc O_v), \Ql)$ the Hecke operator associated to $g$. Then $T(g)$ induces an isomorphism $H_{M, N, I, W}^{' \,  j, \, \nu} \isom H_{M, N, I, W}^{' \,  j, \, \nu+\xi(g)}$.

Suppose that $\xi(g) > 0$ (if not, we take $g^{-1}$ in place of $g$). 
Since $C_G^{P, \; j}(a) \neq 0$, there exists $m \in \Z_{> 0}$ such that $T(g)^{m} \cdot C_G^{P, \; j}(a)$ is supported on the cone $\wh{\Lambda}_{Z_M/ Z_G}^{\mu+m\xi(g)} \supset \wh{\Lambda}_{Z_M/ Z_G}^{\mu}$, but not supported on $\wh{\Lambda}_{Z_M/ Z_G}^{\mu}$. Therefore $T(g)^{2m} \cdot C_G^{P, \; j}(a)$ is supported on the cone $\wh{\Lambda}_{Z_M/ Z_G}^{\mu+2m\xi(g)}$, but not supported on $\wh{\Lambda}_{Z_M/ Z_G}^{\mu+m\xi(g)}$, etc. We deduce that $$C_G^{P, \; j}(a), T(g)^{m} \cdot C_G^{P, \; j}(a), T(g)^{2m} \cdot C_G^{P, \; j}(a), T(g)^{3m} \cdot C_G^{P, \; j}(a), \cdots$$ are linearly independent. So the $\Ql$-vector space generated by $T(g)^{\Z} \cdot C_G^{P, \; j}(a)$ has infinite dimension.
Hence $C_{c}(M(O_v)\backslash M(F_v)/M(O_v), \Ql) \cdot C^{P, \; j}_G(a)$ has infinite dimension.

(3) We deduce from (1) and (2) a contradiction. So $a \in H_{G, N, I, W}^{j \; \on{cusp} }$.
\cqfd

\quad

\begin{defi} (\cite{vincent}, Définition 8.19)
We define $H_{G, N, I, W}^{j, \; \on{Hf}}:=$
$$\{  c \in H_{G, N, I, W}^j , \; C_{c}(K_{N}\backslash G(\mb A)/K_{N}, \Zl) \cdot c \text{ is a finitely generated } \Zl \text{-submodule}  \}.$$
\end{defi}

By definition, $H_{G, N, I, W}^{j, \; \mr{Hf}} \subset H_{G, N, I, W}^{j, \; \mr{Hf-rat}} $. Thus Proposition \ref{prop-cusp-egal-Hfrat-section-6} has the following corollary:
\begin{cor}
\begin{equation*}
H_{G, N, I, W}^{j, \; \mr{Hf}} \subset H_{G, N, I, W}^{j, \; \mr{Hf-rat}} = H_{G, N, I, W}^{j, \; \on{cusp}} .
\end{equation*}
In particular, $H_{G, N, I, W}^{j, \; \mr{Hf}} $ has finite dimension.
\end{cor}

\quad

\appendix

\section{Exact sequences associated to an open and a closed substack of the stack of shtukas}

For simplicity of the notation, we do not write the indices $N$, $I$ and $W$. 


\sssec{}    \label{subsection-appendix-replace-Lambda-by-R}
In the following, we use $\wh{\Lambda}_{G^{\mr{ad}}}^{+, \Q}$. But everything remains true if we replace it by $\frac{1}{r} \wh R_{G^{\mr{ad}}}^{+}$.

\sssec{}   \label{subsection-appendix-def-Cht-G-S-Cht-M-S}
As in \cite{DG15} 7.4.10, we equip the set $\wh{\Lambda}_{G^{\mr{ad}}}^{+, \Q}$ with the order topology, i.e. the one where a base of open subsets is formed by subsets of the form $\{ \lambda \in \wh{\Lambda}_{G^{\mr{ad}}}^{+, \Q} | \lambda \leq \lambda_0 \}$ for $\lambda_0 \in \wh{\Lambda}_{G^{\mr{ad}}}^{+, \Q}$. 
Let $S$ be a subset of $\wh{\Lambda}_{G^{\mr{ad}}}^{+, \Q}$. 
We define $$\Bun_G^S:= \underset{\lambda \in S} \bigcup \Bun_G^{=\lambda} , \quad \Cht_G^S:= \underset{\lambda \in S} \bigcup \Cht_G^{=\lambda} , \quad \Cht_M^{' \, S}:= \underset{\lambda \in S} \bigcup \Cht_M^{' \, =\lambda} .$$
where $\Cht_G^{=\lambda}$ and $\Cht_M^{' \, =\lambda}$ are defined in Definition \ref{def-Cht-G=mu-Cht-M=mu}.
If the subset $S$ is open (resp. closed) in $\wh{\Lambda}_{G^{\mr{ad}}}^{+, \Q}$, then $\Bun_G^S$ is open (resp. closed) in $\Bun_G$. So $\Cht_G^S$ is open (resp. closed) in $\Cht_G$ and $\Cht_M^{' \, S}$ is open (resp. closed) in $\Cht_M^{'}$.

If $S$ is a bounded locally closed subset of $\wh{\Lambda}_{G^{\mr{ad}}}^{+, \Q}$, then $\Cht_G^S$ and $\Cht_M^{' \, S}$ are Deligne-Mumford stacks of finite type.



\sssec{}
Let $\mu \in \wh{\Lambda}_{G^{\mr{ad}}}^{+, \Q} $. Let $S_2 = \{ \lambda \in \wh{\Lambda}_{G^{\mr{ad}}}^{+, \Q} | \lambda \leq \mu \}$. By definition it is an open subset of $\wh{\Lambda}_{G^{\mr{ad}}}^{+, \Q}$ for the order topology of $G^{\mr{ad}}$. It is also open in $\wh{\Lambda}_{\ov M}^{+, \Q}$ for the order topology of $\ov M = M / Z_G$ (because $\lambda \leq^{\ov M} \mu$ implies $\lambda \leq \mu$.). 

Let $S_1$ be an open subset of $S_2$ for the order topology of $G^{\mr{ad}}$. Thus the morphism of stacks $\Cht_G^{S_1} \xrightarrow{j_G} \Cht_G^{S_2}$ (resp. $\Cht_M^{' \, S_1} \xrightarrow{j_M} \Cht_M^{' \, S_2}$) is an open immersion. 
By definition, 
$\Cht_G^{S_2-S_1}$ (resp. $\Cht_M^{' \, S_2-S_1}$) is the closed substack in $\Cht_G^{S_2}$ (resp. $\Cht_M^{' \, S_2} $) which is the complement of $\Cht_G^{S_1}$ (resp. $\Cht_M^{' \, S_1}$).

We define $\Cht_P^{' \, S_2}$ (resp. $\Cht_P^{' \, S_1}$) to be the inverse image of $\Cht_G^{S_2}$ (resp. $\Cht_G^{S_1}$) in $\Cht_P'$. Just as in Lemma \ref{lem-Cht-P-leq-mu-to-Cht-M-leq-mu}, we have $\pi_2: \Cht_P^{' \, S_2} \rightarrow \Cht_M^{' \, S_2}$ (resp. $\pi_1: \Cht_P^{' \, S_1} \rightarrow \Cht_M^{' \, S_1}$). We have $\Cht_P^{' \, S_1} \xrightarrow{j_P} \Cht_P^{' \, S_2}$, which is an open immersion. We define $\Cht_P^{' \, S_2-S_1}:= \Cht_P^{' \, S_2 } \, \cap \, \pi^{-1} ( \Cht_M^{' \, S_2-S_1} ).$
It is a closed substack in the complement of $\Cht_P^{' \, S_1}$ in $\Cht_P^{' \, S_2}$, but may not be equal to it. 

\begin{lem}  \label{lem-ouvert-ferme-cht}
The following diagram of algebraic stacks is commutative:
\begin{equation} \label{diagram-3-3-cht}
\xymatrix{
\Cht_G^{S_2-S_1} \ar@{^{(}->}[d]_{i_G}
& \Cht_P^{' \, S_2-S_1} \ar[l]_{i_{12}}  \ar[r]^{\pi_{12}}  \ar@{^{(}->}[d]_{i_P}
& \Cht_M^{' \, S_2-S_1}  \ar@{^{(}->}[d]^{i_M} \\
\Cht_G^{S_2 } 
&  \Cht_P^{' \, S_2 }  \ar[l]_{i_2}  \ar[r]^{\pi_2}  
& \Cht_M^{' \, S_2 } \\
\Cht_G^{S_1 } \ar@{^{(}->}[u]^{j_G}
& \Cht_P^{' \, S_1 }   \ar[l]_{i_1}  \ar[r]^{\pi_1}  \ar@{^{(}->}[u]^{j_P}
& \Cht_M^{' \, S_1 }  \ar@{^{(}->}[u]_{j_M} 
}
\end{equation}
Moreover, the left bottom square and the right top square are Cartesian.
\cqfd
\end{lem}

%

\sssec{}
For any $j$, any $\nu \in \wh \Lambda_{Z_M / Z_G}^{\Q}$ 
and any bounded locally closed subset $S \subset \wh{\Lambda}_{G^{\mr{ad}}}^{+, \Q}$, 
we define
$$H _{G}^{j, \, S }  := H_c^j ( \Cht_{G, \ov{\eta^I}}^{S} / \Xi_G, \mc F_G ); \quad H _{M}^{' \, j, \, S, \, \nu }  := H_c^j ( \Cht_{M, \ov{\eta^I} }^{ ' \, S, \, \nu} / \Xi_G, \mc F_M' ) .$$

\sssec{}
By Proposition \ref{prop-Varshavsky-proper}, the restriction of morphism $i_1$ (resp. $i_2$) to $\ov{\eta^I}$ is proper. The restriction of morphism $i_{12}$ to $\ov{\eta^I}$ is also proper because $\Cht_P^{' \, S_2-S_1} \rightarrow \Cht_G^{S_2-S_1} \underset{\Cht_G^{S_2}} \times \Cht_P^{' \, S_2}$ is a closed immersion. Moreover $i_1$, $i_2$ and $i_{12}$ are schematic.
Applying the construction in Section 3 to each line in diagram (\ref{diagram-3-3-cht}), respectively, we obtain the constant term morphism $C_G^{P, \, j, \, S_1}: H _{G}^{j, \, S_1 }   \rightarrow H _{M}^{' \, j, \, S_1 }$, $C_G^{P, \, j, \, S_2}$ and $C_G^{P, \, j, \, S_2-S_1}$ (note that the morphism $\pi_{12, d}: \Cht_P^{' \, S_2-S_1} \rightarrow \wt{\Cht}_M^{' \, S_2-S_1}$ is smooth because the right top square of diagram (\ref{diagram-3-3-cht}) is Cartesian).

\sssec{}
Diagram (\ref{diagram-3-3-cht}) induces a diagram of cohomology groups with compact support for which we will study the commutativity:
\begin{equation}  \label{diagram-TC-suite-exacte-longue-general}
\xymatrix{
\cdots \ar[r]
& H _{G}^{j-1, \,  S_2-S_1 }   \ar[r]  \ar[d]^{C_G^{P, \, j-1, \, S_2-S_1}} 
& H _{G}^{j, \,  S_1}  \ar[r]   \ar[d]^{C_G^{P, \, j, \, S_1}}
& H _{G}^{j, \, S_2}   \ar[r]   \ar[d]^{C_G^{P, \, j, \, S_2}}
& H _{G}^{j, \,  S_2-S_1 }   \ar[r]  \ar[d]^{C_G^{P, \, j, \, S_2-S_1}}  
& \cdots \\
\cdots \ar[r]
& H _{M}^{' \, j-1, \, S_2-S_1 }  \ar[r] 
& H _{M}^{' \, j, \, S_1 }  \ar[r]   
& H _{M}^{' \, j, \, S_2 }   \ar[r]   
& H _{M}^{' \, j, \, S_2-S_1 }   \ar[r]
& \cdots
}
\end{equation}
The horizontal maps are the long exact sequences associated to an open substack and the complementary closed substack. The vertical maps are the constant term morphisms.

\begin{lem} \label{lem-TC-suite-exacte-longue}
For any $j$, the following diagram is commutative:
$$
\xymatrix{
 H _{G}^{j, \,  S_1}  \ar[r]   \ar[d]^{C_G^{P, \, j, \, S_1}}
& H _{G}^{j, \, S_2}   \ar[r]   \ar[d]^{C_G^{P, \, j, \, S_2}}
& H _{G}^{j, \,  S_2-S_1 }    \ar[d]^{C_G^{P, \, j, \, S_2-S_1}}   \\
 H _{M}^{' \, j, \, S_1 }  \ar[r]   
& H _{M}^{' \, j, \, S_2 }   \ar[r]   
& H _{M}^{' \, j, \, S_2-S_1 } 
}
$$
\end{lem}

\dem 
We denote by $p_G: \Cht_G^{S_2} \rightarrow \eta^I$ and $p_M: \Cht_M^{' \, S_2} \rightarrow \eta^I$. 
For $S=S_1$ or $S_2$ or $S_2-S_1$, denote $\mc F_G^{S}: = \restr{\mc F_G}{\Cht_G^{S}}$ and $\mc F_M^{S}: = \restr{\mc F_M'}{\Cht_M^{' \, S}}$. 
Note that $\mc F_G^{S_1} = (j_G)^* \mc F_G^{S_2} $ and $\mc F_G^{S_2-S_1} = (i_G)^* \mc F_G^{S_2}$. Similarly $\mc F_M^{S_1} = (j_M)^* \mc F_M^{S_2} $ and $\mc F_M^{S_2-S_1} = (i_M)^* \mc F_M^{S_2}$.
Lemma \ref{lem-TC-suite-exacte-longue} will follow from the commutativity of the following diagram of complexes in $D_c^b(\eta^I, \Ql)$:
\begin{equation}    \label{equation-exact-sequence-TC-as-complex}
\xymatrix{
\xymatrixrowsep{2pc}
\xymatrixcolsep{3pc}
(p_G)_! (j_G)_!(j_G)^* \mc F_G^{S_2}  \ar[r]   \ar[d]^{\mc C_G^{P, \, S_1}}
& (p_G)_!  \mc F_G^{S_2}   \ar[r]   \ar[d]^{\mc C_G^{P, \, S_2}}
& (p_G)_!  (i_G)_!(i_G)^* \mc F_G^{S_2}   \ar[d]^{\mc C_G^{P, \, S_2-S_1}} \\
(p_M)_!  (j_M)_!(j_M)^* \mc F_M^{S_2} \ar[r] 
& (p_M)_! \mc F_M^{S_2}   \ar[r]   
& (p_M)_!  (i_M)_!(i_M)^* \mc F_M^{S_2}   
}
\end{equation}


The commutativity of the left square is induced by (1) and (2) below. The commutativity of the right square is induced by (3) and (4) below.

\quad

We consider the left square of (\ref{equation-exact-sequence-TC-as-complex})
\begin{equation}     \label{equation-(1)-(2)-together}
\xymatrix{
\xymatrixrowsep{1pc}
\xymatrixcolsep{7pc}
(p_G)_! (j_G)_!(j_G)^* \mc F_G^{S_2} \ar[r]^{\on{Tr}_{j_G}}   \ar[d]  \ar@{}[dr]|{(1)}
& (p_G)_!  \mc F_G^{S_2}   \ar[d] \\
(p_M)_!  (j_M)_! (\pi_1)_! (i_1)^* (j_G)^* \mc F_G^{S_2} \ar[r]^{\quad \quad \on{Tr}_{j_P}}   \ar[d]^{f_1}   \ar@{}[dr]|{(2)}
& (p_M)_!  (\pi_2)_! (i_2)^* \mc F_G^{S_2}   \ar[d]^{f_2}  \\
(p_M)_!  (j_M)_!(j_M)^* \mc F_M^{S_2} \ar[r]_{\quad \quad \on{Tr}_{j_M}}   
& (p_M)_! \mc F_M^{S_2}  
}
\end{equation}
where (1) and (2) are detailed below.

(1) The following diagram of functors is commutative:
$$
\xymatrix{
(p_G)_!  (j_G)_! (j_G)^*  \ar[rr]^{\on{Tr}_{j_G}}   \ar[d]^{\on{adj}_{i_1} }    \ar[rd]^{\on{adj}_{i_2} }
& & (p_G)_!  \ar[d]^{\on{adj}_{i_2} }  \\
(p_G)_!  (j_G)_! (i_1)_! (i_1)^* (j_G)^*    \ar[r]^{\simeq}_{ (\ast) }   \ar[d]^{\simeq}
& (p_G)_! (i_2)_! (i_2)^*  (j_G)_! (j_G)^* \ar[r]^{\quad \quad \on{Tr}_{j_G} } 
&  (p_G)_! (i_2)_! (i_2)^* \ar[d]^{\simeq} \\
(p_M)_!  (j_M)_! (\pi_1)_! (i_1)^* (j_G)^* \ar[r]^{\simeq}   
& (p_M)_!  (\pi_2)_! (j_P)_! (j_P)^* (i_2)^*  \ar[r]^{\quad \quad \on{Tr}_{j_P} } 
& (p_M)_!  (\pi_2)_! (i_2)^* 
}
$$
where $(\ast)$ is given by $ (j_G)_! (i_1)_! (i_1)^*  \simeq  (i_2)_! (j_P)_! (i_1)^*  \simeq  (i_2)_! (i_2)^*  (j_G)_!$, the last isomorphism is the proper base change for the left bottom square of diagram (\ref{diagram-3-3-cht}). The commutativity of (1) follows from the fact that the adjunction morphism commutes with base change and the trace morphism commutes with base change (\cite{sga4} XVIII Théorème 2.9).

(2) Taking (\ref{diagram-cht-P-cht-M-pi-d}) into account, we have a commutative diagram, where $\pi_2$ (resp. $\pi_1$) is the composition $\wt{\pi_{2, \underline{d}}^0} \circ \pi_{2, \underline{d}}$ (resp. $\wt{\pi_{1, \underline{d}}^0} \circ \pi_{1, \underline{d}}$) for some $d$ large enough as in Proposition \ref{prop-d-assez-grand}. 
\begin{equation}  \label{diagramme-cht-TC-ouvert}
\xymatrixrowsep{2pc}
\xymatrixcolsep{3pc}
\xymatrix{
 \Cht_P^{' \, S_2 }   \ar[r]^{\pi_{2, d}}   \ar@{}[dr]|{(b)}
& \wt{\Cht}_M^{' \, S_2 }   \ar[r]^{\wt{\pi_{2, d}^0}}    \ar@{}[dr]|{(c)}
& \Cht_M^{' \, S_2 }  \\
\Cht_P^{' \, S_1 }    \ar[r]^{\pi_{1, d}}  \ar@{^{(}->}[u]^{j_P}
& \wt{\Cht}_M^{' \, S_1 }  \ar[r]^{\wt{\pi_{1, d}^0}}   \ar@{^{(}->}[u]^{\wt{j_M}}
& \Cht_M^{' \, S_1 }  \ar@{^{(}->}[u]_{j_M} 
}
\end{equation}  
The square (c) is Cartesian. The square (b) may not be Cartesian. As in Lemma \ref{lem-pi-d-est-lisse}, $\pi_{1, d}$ and $\pi_{2, d}$ are smooth. We have $\dim (\pi_{1, d} ) = \dim ( \pi_{2, d} ) = d\cdot |I| \dim U$. We denote this dimension by $m$.

By (\ref{equation-base-change-with-trace}) and (\ref{equation-F-G-vers-F-M-long}), the morphism $f_1$ (resp. $f_2$) defined in diagram (\ref{equation-(1)-(2)-together}) is the composition of $\on{Tr}_{\pi_{1, d}}: (\pi_{1, d})_! (\pi_{1, d})^* \rightarrow \Id[-2m](-m)$ (resp. $\on{Tr}_{\pi_{2, d}}: (\pi_{2, d})_! (\pi_{2, d})^* \rightarrow \Id[-2m](-m)$) with some isomorphisms. By \cite{sga4} XVIII Théorème 2.9, the trace morphism is compatible with composition, thus $$\on{Tr}_{\pi_{2, d}} \circ \on{Tr}_{j_P} \simeq \on{Tr}_{\pi_{2, d} \circ j_P} \simeq \on{Tr}_{\wt{j_M} \circ \pi_{1, d}} \simeq  \on{Tr}_{\wt{j_M}} \circ \on{Tr}_{\pi_{1, d}} $$
where the middle isomorphism is due to the commutativity of (b). Moreover, by $loc.cit.$ the trace morphism is compatible with base change, thus $$\on{Tr}_{\wt{j_M}} = (\wt{\pi_{1, d}^0})^* \on{Tr}_{j_M} .$$ We deduce that (2) is commutative.

\quad

Now we consider the right square of (\ref{equation-exact-sequence-TC-as-complex}):
\begin{equation}    \label{equation-(3)-(4)-together}
\xymatrix{
\xymatrixrowsep{1pc}
\xymatrixcolsep{7pc}
 (p_G)_!  \mc F_G^{S_2}  \ar[r]^{\on{adj}_{i_G}}   \ar[d]    \ar@{}[dr]|{(3)}
& (p_G)_!  (i_G)_!(i_G)^* \mc F_G^{S_2}  \ar[d] \\
 (p_M)_!  (\pi_2)_! (i_2)^* \mc F_G^{S_2}  \ar[r]^{\on{adj}_{i_P} \quad \quad }    \ar[d]^{f_2}    \ar@{}[dr]|{(4)}
& (p_M)_! (i_M)_!   (\pi_{12})_! (i_{12})^* (i_G)^* \mc F_G^{S_2}  \ar[d]^{f_{12}} \\
(p_M)_! \mc F_M^{S_2}  \ar[r]^{\on{adj}_{i_M}}    
& (p_M)_!  (i_M)_!(i_M)^* \mc F_M^{S_2}  
}
\end{equation}
where (3) and (4) are detailed below.

(3) The following diagram of functors is commutative:
$$
\xymatrix{
(p_G)_!    \ar[rr]^{\on{adj}_{i_G} }    \ar[d]^{\on{adj}_{i_2} } 
& & (p_G)_!  (i_G)_!(i_G)^*   \ar[d]^{\on{adj}_{i_{12}} }  \\
(p_G)_! (i_2)_! (i_2)^*  \ar[r]^{\on{adj}_{i_P} \quad \quad }    \ar[d]^{\simeq}  
&  (p_G)_! (i_2)_!(i_P)_! (i_P)^* (i_2)^*  \ar[d]^{\simeq}   \ar[r]^{\simeq}
& (p_G)_!  (i_G)_!  (i_{12})_! (i_{12})^*   (i_G)^*   \ar[d]^{\simeq}   \\
(p_M)_! (\pi_2)_! (i_2)^*  \ar[r]^{\on{adj}_{i_P} \quad \quad } 
& (p_M)_! (\pi_2)_! (i_P)_! (i_P)^* (i_2)^*  \ar[r]^{\simeq}
& (p_M)_! (i_M)_!  (\pi_{12})_! (i_{12})^* (i_G)^* 
}
$$

(4) Taking (\ref{diagram-cht-P-cht-M-pi-d}) into account, we have a commutative diagram, where $\pi_{12}$ is the composition $\wt{\pi_{12, \underline{d}}^0} \circ \pi_{12, \underline{d}}$.
$$
\xymatrixrowsep{2pc}
\xymatrixcolsep{3pc}
\xymatrix{
 \Cht_P^{' \, S_2-S_1}  \ar[r]^{\pi_{12, d}}  \ar@{^{(}->}[d]_{i_P}   \ar@{}[dr]|{(e)}
& \wt{\Cht}_M^{' \, S_2-S_1}    \ar[r]^{\wt{\pi_{12, d}^0}}  \ar@{^{(}->}[d]_{\wt{i_M}}   \ar@{}[dr]|{(f)}
& \Cht_M^{' \, S_2-S_1}  \ar@{^{(}->}[d]^{i_M} \\
  \Cht_P^{' \, S_2 }    \ar[r]^{\pi_{2, d}}  
& \wt{\Cht}_M^{' \, S_2 }   \ar[r]^{\wt{\pi_{2, d}^0}}
& \Cht_M^{' \, S_2 } 
}
$$
The squares (e) and (f) are Cartesian.

By (\ref{equation-base-change-with-trace}) and (\ref{equation-F-G-vers-F-M-long}), $f_{12}$ defined in diagram (\ref{equation-(3)-(4)-together}) is the composition of $\on{Tr}_{\pi_{12, d}}: (\pi_{12, d})_! (\pi_{12, d})^* \rightarrow \Id[-2(\dim \pi_{12, d})](-\dim \pi_{12, d})$ with some isomorphisms. By \cite{sga4} XVIII Théorème 2.9, the trace morphism is compatible with base change, thus $$\on{Tr}_{\pi_{12, d}} = (\wt{i_M})^* \on{Tr}_{\pi_{2, d}} .$$ 
We deduce that (4) is commutative.
\cqfd

\begin{rem}
We don't know if the complete diagram (\ref{diagram-TC-suite-exacte-longue-general}) is commutative.
\end{rem}

\quad

\section{Lemma of the cubic commutative diagram}

\begin{lem}   \label{lem-cube}
Let $\ms X, \ms Y, \ms Z, \ms W, \ms X', \ms Y', \ms Z', \ms W'$ be algebraic stacks. Suppose that we have two Cartesian squares:
$$
\xymatrixrowsep{1pc}
\xymatrixcolsep{1pc}
\xymatrix{
\ms Z  \ar[r]  \ar[d]
& \ms Y  \ar[d]^h
& & \ms Z'   \ar[r]  \ar[d]
& \ms Y' \ar[d]^{h'} \\
\ms X  \ar[r]^g
& \ms W
& & \ms X'  \ar[r]^{g'}
& \ms W'
}
$$
If these two squares are the front and back faces of a commutative diagram: 
\begin{equation}
\xymatrixrowsep{1pc}
\xymatrixcolsep{1pc}
\xymatrix{
\ms Z' \ar[rr] \ar[dd]  \ar[rd]^{f_{\ms Z}}
& & \ms Y' \ar@{.>}[dd]  \ar[rd]^{f_{\ms Y}}
& \\
& \ms Z \ar[rr]  \ar[dd]
&  &  \ms Y \ar[dd] \\
\ms X' \ar@{.>}[rr]  \ar[rd]^{f_{\ms X}}
&  &  \ms W'  \ar@{.>}[rd]^{f_{\ms W}}
& \\
& \ms X \ar[rr]
&  &  \ms W
}
\end{equation}
then the fibers $f_{\ms Z}$, $f_{\ms X}$, $f_{\ms Y}$ and $f_{\ms W}$ form a Cartesian square.
\end{lem}

Concretely, let $T$ be a scheme. For any morphism $T \rightarrow \ms Z$, we have the compositions of morphisms $T \rightarrow \ms Z \rightarrow \ms X$, $T \rightarrow \ms Z \rightarrow \ms Y$ and $T \rightarrow \ms Z \rightarrow \ms W$.
We denote by $\ms Z_T$ (resp. $\ms X_T$, $\ms Y_T$, $\ms W_T$) the fiber of $f_{\ms Z}$ (resp. $f_{\ms X}$, $f_{\ms Y}$, $f_{\ms W}$) over $T$. The lemma says that $\ms Z_T$ is equivalent to $\ms X_T \underset{\ms W_T} \times \ms Y_T$.

%
%
%

\dem
We will prove a more general statement: suppose that we have another Cartesian squares:
$$
\xymatrixrowsep{1pc}
\xymatrixcolsep{1pc}
\xymatrix{
\ms Z''  \ar[r]  \ar[d]
& \ms Y''  \ar[d]^{h''} \\
\ms X''  \ar[r]^{g''}
& \ms W''
}
$$
and a commutative diagram: 
\begin{equation}
\xymatrixrowsep{1pc}
\xymatrixcolsep{1pc}
\xymatrix{
\ms Z'' \ar[rr] \ar[dd]  \ar[rd]
& & \ms Y'' \ar@{.>}[dd]  \ar[rd]
& \\
& \ms Z \ar[rr]  \ar[dd]
&  &  \ms Y \ar[dd] \\
\ms X'' \ar@{.>}[rr]  \ar[rd]
&  &  \ms W''  \ar@{.>}[rd]
& \\
& \ms X \ar[rr]
&  &  \ms W
}
\end{equation}
Then we have a canonical isomorphism:
\begin{equation}    \label{equation-prod-fibre-Z-prod-fibre-X-Y-W}
\ms Z' \times_{\ms Z} \ms Z'' \isom (  \ms X' \times_{\ms X} \ms X''   ) \underset{\ms W' \times_{\ms W} \ms W''}{\times} (\ms Y' \times_{\ms Y} \ms Y'') .
\end{equation}

In fact, by definition, we have $$\ms Z' \times_{\ms Z} \ms Z'' \simeq (  \ms X' \times_{\ms W'} \ms Y'   ) \underset{\ms X \times_{\ms W} \ms Y}{\times} (\ms X'' \times_{\ms W''} \ms Y'') .$$
For any scheme $S$, the $S$-points of both sides of (\ref{equation-prod-fibre-Z-prod-fibre-X-Y-W}) classify the data of $S$-points $x'$ in $\ms X'$, $x'' $ in $\ms X''$, $y'$ in $\ms Y'$, $y'' $ in $ \ms Y''$, an isomorphism between the images of $x'$ and $x''$ in $\ms X$, an isomorphism between the images of $y'$ and $y''$ in $\ms Y$, an isomorphism between the images of $x'$ and $y'$ in $\ms W'$, an isomorphism between the images of $x''$ and $y''$ in $\ms W''$, such that the diagram deduced from these four isomorphisms between the images of $x'$, $x''$, $y'$, $y''$ in $\ms W$ is commutative. 
We deduce (\ref{equation-prod-fibre-Z-prod-fibre-X-Y-W}).

The lemma is the special case when $\ms X'' = \ms Y'' = \ms W'' = \ms Z'' = T$.
\cqfd

\quad

\end{document}